\documentclass{article}
\usepackage{amssymb,amsmath,titlesec}
\usepackage{verbatim,amsthm,slashed}
\usepackage[all]{xy}

\theoremstyle{plain}
\newtheorem{thm}[equation]{Theorem}
\newtheorem{prop}[equation]{Proposition}
\newtheorem{lemma}[equation]{Lemma}
\renewcommand{\qedsymbol}{\textsquare}

\theoremstyle{definition}
\newtheorem{defn}[equation]{Definition}

\renewcommand{\theequation}{\thesubsection.\arabic{equation}}

\def\beq{\begin{eqnarray}}
\def\eeq{\end{eqnarray}}
\def\beqq{\begin{eqnarray*}}
\def\eeqq{\end{eqnarray*}}
\def\t{\textrm}
\def\bb{\mathbb}
\def\mc{\mathcal}
\def\mf{\mathfrak}
\def\wh{\widehat}

\def\M{\widehat{\mf{M}}}
\def\ev{\t{even}}
\def\od{\t{odd}}
\def\r{\t{red}}

\title{\sf\bfseries SUPERSYMMETRIC FIELD THEORIES AND COHOMOLOGY}
\author{\sf Pokman Cheung}
\date{\sf October 30, 2006}

\begin{document}

\maketitle

\begin{abstract}
This is the Ph.D.~dissertation of the author.
The project has been motivated by the conjecture that 
the Hopkins-Miller tmf spectrum can be described in terms of 
`spaces' of conformal field theories.
In this dissertation, spaces of field theories are constructed as 
classifying spaces of categories whose objects are certain types of 
field theories. 
If such a category has a symmetric monoidal structure and its components 
form a group, by work of Segal, its classifying space is 
an infinite loop space and defines a cohomology theory. 
This has been carried out for two classes of field theories:
(i) For each $n\in\bb{Z}$, there is a category $\mc{SEFT}_n$ 
whose objects are the Stolz-Teichner $(1|1)$-dimensional 
super Euclidean field theories of degree $n$. 
It is proved that the classifying space $|\mc{SEFT}_n|$ represents 
the degree-$n$ $K$ or $KO$ cohomology, depending on the coefficients 
of the field theories. 
(ii) For each $n\in\bb{Z}$, there is a category $\mc{AFT}_n$ whose objects
are a kind of $(2|1)$-dimensional field theories called 
`annular field theories,' defined using supergeometric versions of 
circles and annuli only. 
It is proved that the classifying space $|\mc{AFT}_n|$ represents the 
degree-$n$ elliptic cohomology associated with the Tate curve. 
To the author's knowledge, this is the first time the definitions of 
low-dimensional supersymmetric field theories are given in full detail.
\end{abstract}

\setcounter{section}{-1}

\titleformat{\section}[block]
  {\large\sffamily\filcenter\bfseries}
  {}{.5em}{}

\section{INTRODUCTION}

{\bf Background.} 
One of the first connections between (supersymmetric) quantum field 
theories (QFTs) and topology appeared in Witten's approach to the 
Atiyah-Singer index theorem. 
He argued that the circle partition function of 
a canonical $1$-dimensional QFT associated to 
a compact spin manifold $M^{2n}$ equals $\hat{A}(M)$.~\cite{Witten.indD} 
(See also \cite{LM}, III, \S 17.)
Later, Ochanine discovered the elliptic genus, 
which is a bordism invariant of spin manifolds and takes values of 
integral modular forms.~\cite{Ochanine}
Witten identified the elliptic genus as the torus partition function 
of a $2$-dimensional conformal field theory (CFT); 
furthermore, variations of CFTs give rise to other modular form valued 
bordism invariants, e.g. the Witten genus.~\cite{Witten.ell.genus}
The field-theoretic interpretation of these bordism invariants 
provides intuitive explanations of some of their properties, 
e.g. integrality, modularity and rigidity, 
which can also be proved rigorously.~\cite{Segal.Ell,BottTaubes}

Elliptic cohomology theories are cohomology theories where 
the elliptic genus, the Witten genus, etc., of families take value. 
Such theories were first constructed by Landweber, Ravenel and Stong, 
using the Landweber exact functor theorem.~\cite{LRS}
There are plenty of elliptic cohomology theories, each associated to an 
elliptic curve. 
A lot of work has gone into the search of a `universal' theory. 
Ando, Hopkins and Strickland defined and studied elliptic spectra; 
then Hopkins and Miller constructed the tmf spectrum, 
which is a suitable inverse limit in the category of 
elliptic spectra.~\cite{AHS,Hopkins.ICM}
The families Witten genus, for instance, can now be described as 
a multiplicative map from a Thom spectrum to tmf. 
However, these cohomology theories and spectra have been defined 
using homotopy theory and lack a geometric interpretation. 

In \cite{Segal.Ell}, Segal proposed that an elliptic cohomology class 
is represented by something he called an elliptic object. 
An elliptic object over a point is a CFT, in the sense of \cite{Segal.CFT},
and an elliptic object over a general space $X$ can be thought of as 
`a CFT parametrized by $X$.' 
It has also been conjectured that the tmf spectrum can be recovered 
in terms of spaces of CFTs. 
Stolz and Teichner have significantly modified and elaborated upon 
the precise notions of CFTs and elliptic objects that should be related 
to tmf; 
they have also constructed canonical elliptic objects 
which are expected to represent certain tmf Euler classes.~\cite{ST}

\vspace{0.1in}
\noindent
{\bf Goal and results.}
This work focuses on constructing `spaces' of $1$- and $2$-dimensional 
QFTs and relating them to cohomology theories.  
Following the axiomatic point of view of Atiyah, Segal and Witten, 
we think of a QFT as a vector space representation of a bordism category. 
Our approach is to define symmetric monoidal categories 
whose objects are certain kinds of QFTs, 
and study their classifying spaces. 
By \cite{Segal.Gamma}, 
these classifying spaces are infinite loop spaces, 
\footnote{
The components of each of our symmetric monoidal categories 
form a group.
}
hence defining cohomology theories. 

The main results are as follows. 

{\renewcommand{\theequation}{\ref{SEFT.Fred}}
\thm 
\renewcommand{\theequation}{\thechapter.\arabic{equation}}
For each $n\in\bb{Z}$, there is a (small) category $\mc{SEFT}_n$ 
whose objects are 
`$1$-dimensional super Euclidean field theories of degree $n$,' 
such that 
\beqq
|\mc{SEFT}_n|\simeq KO_n\t{ or }K_n\;,
\eeqq
where $KO_n$ (resp. $K_n$) is the $n$-th space of the periodic $KO$-theory 
(resp. $K$-theory) spectrum.
Whether we have $KO_n$ or $K_n$ depends on the coefficients in the field 
theories.
}

{\renewcommand{\theequation}{\ref{AFT.TateK}}
\thm 
\renewcommand{\theequation}{\thechapter.\arabic{equation}}
For each $n\in\bb{Z}$, there is a (small) category $\mc{AFT}_n$ 
whose objects are 
`super annular field theories of degree $n$,' such that 
\beqq
|\mc{AFT}_n|\simeq (K_{\t{Tate}})_n\;,
\eeqq
where the right hand side is the $n$-th space of the elliptic spectrum 
$K_{\t{Tate}}$ associated with the Tate curve.
}

Here is a section-by-section overview. 
In \S 2.1, we define Stolz and Teichner's 
$1$-dimensional super Euclidean field theories 
as functors between super categories. 
The use of super categories has been explained and emphasized by 
Stolz and Teichner. 
However, to the author's knowledge, its details have never appeared 
in the literature. 
\S 2.2 concerns the definition of $\mc{SEFT}_n$. 
We digress in \S 2.3 to study another sequence of categories $\mc{V}_n$, 
which are related to $\mc{SEFT}_n$, 
before proving theorem \ref{SEFT.Fred} in \S 2.4. 
In \S 3.1, we construct certain bordism (super) categories 
using supergeometric versions of $S^1$ and $S^1\times[0,t]$, 
and define `super annular field theories.' 
Finally, in \S 3.2, we define $\mc{AFT}_n$ and prove theorem 
\ref{AFT.TateK}. 
The proofs of two technical lemmas needed in \S 2.1 are given in 
an appendix.

\vspace{0.1in}
\noindent
{\bf Acknowledgements.}
The author is indebted to his Ph.D. advisor, Ralph Cohen, who initially 
pointed him to this very exciting research topic, and has provided 
generous support, guidance and encouragement.
He would also like to thank Stephan Stolz and Peter Teichner for their 
interest in this project and for explaining various aspects of their 
groundbreaking work, which largely motivated this project; 
and S\o ren Galatius for his suggestion on the proof of 
lemma \ref{C2.Fred}.

\newpage

\titleformat{\section}[block]
  {\large\sffamily\filcenter\bfseries}
  {CHAPTER \thesection:}{.5em}{}
\titleformat{\subsection}[block]
  {\large\sffamily\bfseries}
  {\thesubsection.}{.5em}{}

\section{ONE-DIMENSIONAL FIELD THEORIES AND $KO$-THEORY}
\label{chap.SEFT}

Throughout this chapter, one can replace real coefficients 
with complex ones and obtain parallel results for $K$-theory. 
For simplicity, we focus on $KO$.
\emph{All gradings are $\bb{Z}/2$-gradings} and all tensor 
products of graded vector spaces or algebras are graded tensor products. 
The even and odd parts of a graded vector space $V$ are denoted $V^\ev$ and 
$V^\od$. 


\subsection{Stolz and Teichner's One-Dimensional Field Theories} 
\label{sec.1SEFT}
{\bf Super manifolds of dimensions $(0|1)$ and $(1|1)$.}
We first recall the definition of (smooth) super manifolds. \cite{QFS.susy} 
A \emph{$(p|q)$-dimensional super manifold} $M=(M_\r,\mc{O}_M)$ consists of a 
smooth $p$-dimensional manifold $M_\r$, and a sheaf $\mc{O}_M$ of real graded 
commutative algebras on $M_\r$. 
Any point in $M_\r$ has a neighborhood $U$ such that 
$\mc{O}_M|_U\cong \mc{O}_{M_\r}|_U\otimes\Lambda^*\bb{R}^q$, 
where $\mc{O}_{M_\r}$ is the structure sheaf of $M_\r$ and 
$\Lambda^*\bb{R}^q$ is the exterior algebra generated by $\bb{R}^q$. 
Elements in $\mc{O}(M):=\mc{O}_M(M_\r)$ are refered to as \emph{(global) 
functions} on $M$. 
Examples of super manifolds are 
$\bb{R}^{p|q}=(\bb{R}^p,\mc{O}_{\bb{R}^p}\otimes\Lambda^*\bb{R}^q)$. 
A map $f=(f_\r,f^*):M\rightarrow N$ of two super manifolds consists of a map 
$f_\r:M_\r\rightarrow N_\r$ of smooth manifolds and a morphism 
$f^*:\mc{O}_N\rightarrow (f_\r)_*\mc{O}_M$ of sheaves of graded commutative 
algebras extending $f_\r^*$. 

In this chapter, by a \emph{Euclidean $(0|1)$-manifold} $(Z,\xi)$, we mean 
a $(0|1)$-manifold $Z$ together with a $1$-form $\xi$ on $Z$ of 
the form $\lambda d\lambda$ for some $\lambda\in\mc{O}(Z)^\od$ that is 
nowhere vanishing. 
\footnote{
The notion of differential form extends to smooth super manifolds. 
See \cite{QFS.susy}.
}
On the other hand, a \emph{Euclidean $(1|1)$-manifold} $(Y,\omega)$ is a 
$(1|1)$-manifold $Y$ with a `metric' $\omega$. 
Following \cite{ST}, a $1$-form $\omega$ on $Y$ is called a \emph{metric} 
if every point in $Y_\r$ has a neighborhood $U$ and 
$s_U\in\mc{O}(Y_U)^\ev$, $\lambda_U\in\mc{O}(Y_U)^\od$, 
\footnote{ \label{restrict.smfld}
Given a super manifold $M$ and an open subset $U\subset M_\r$, we denote by 
$M_U$ the super submanifold $(U,\mc{O}_M|_U)$.
}
such that $ds_U$, $\lambda_U$ are nowhere vanishing and 
$\omega|_{Y_U}=ds_U+\lambda_U d\lambda_U$. 
For any Euclidean $(0|1)$- or $(1|1)$-manifolds $(M_i,\chi_1)$, $i=1,2$, we 
write $\phi:(M_1,\chi_1)\rightarrow(M_2,\chi_2)$ for a map 
$\phi:M_1\rightarrow M_2$ of super manifolds that satisfies 
$\chi_1=\phi^*\chi_2$; the set of these maps is denoted 
$\t{Hom}(M_1,\chi_1;\,M_2,\chi_2)$. 

Consider the Euclidean $(0|1)$-manifold $(\bb{R}^{0|1},\lambda d\lambda)$ and 
the Euclidean $(1|1)$-manifold $(\bb{R}^{1|1},ds+\lambda d\lambda)$, where 
$\lambda$ denotes both the element 
$1\in\bb{R}=\Lambda^1\bb{R}=\mc{O}(\bb{R}^{0|1})^\od$ and its pullback in 
$\mc{O}(\bb{R}^{1|1})^\od$ along the canonical projection 
$\bb{R}^{1|1}\rightarrow\bb{R}^{0|1}$, while $s\in\mc{O}(\bb{R}^{1|1})^\ev$ 
is the pullback of the identity function on $\bb{R}$ along the other 
canonical projection $\bb{R}^{1|1}\rightarrow\bb{R}$. 
Define the following maps 
\beqq
&& \epsilon:\bb{R}^{0|1}\rightarrow\bb{R}^{0|1}, \quad 
  \epsilon^*(\lambda)=-\lambda, \\
\t{or} && \epsilon:\bb{R}^{1|1}\rightarrow\bb{R}^{1|1}, \quad
  \epsilon^*(s,\lambda)=(s,-\lambda), \\ 
&& \gamma_t:\bb{R}^{0|1}\rightarrow\bb{R}^{1|1}, \quad
  \gamma_t^*(s,\lambda)=(t,\lambda),\quad t\in\bb{R} \\ 
&& \tau_t:\bb{R}^{1|1}\rightarrow\bb{R}^{1|1}, \quad
  \tau_t^*(s,\lambda)=(s+t,\lambda),\quad t\in\bb{R}. 
\eeqq
Notice that $\gamma_t^*(s)=t$ means 
$(\gamma_t)_\r$ maps $\bb{R}^0$ to $\{t\}$; 
$\tau_t^*(s)=s+t$ means 
$(\tau_t)_\r$ is translation by $t$ on $\bb{R}$. 
We have 
\beq
&& \t{Hom}(\bb{R}^{0|1},\lambda d\lambda;\,\bb{R}^{0|1},\lambda d\lambda)=
  \{1,\epsilon\}, \label{Hom.00} \\
&& \t{Hom}(\bb{R}^{0|1},\lambda d\lambda;\,\bb{R}^{1|1},ds+\lambda d\lambda)=
  \{\gamma_t,\gamma_t\epsilon\}, \label{Hom.01} \\
&& \t{Hom}(\bb{R}^{1|1},ds+\lambda d\lambda;\,\bb{R}^{1|1},
  ds+\lambda d\lambda)=\{\tau_t,\tau_t\epsilon\}. \label{Hom.11}
\eeq
For example, a map $\tau:\bb{R}^{1|1}\rightarrow\bb{R}^{1|1}$ is determined 
by the pair $\tau^*(s,\lambda)=(a(s),b(s)\lambda)$. 
The condition $\tau^*(ds+\lambda d\lambda)=ds+\lambda d\lambda$ is then 
equivalent to $a'(s)=1$ and $b(s)^2=1$. 
This explains (\ref{Hom.11}). 

{\it Note.}
The definition of a Euclidean $(1|1)$-manifold $(Y,\omega)$ can be restated 
as follows: every point in $Y_\r$ has a neighhorhood $U$ and a map 
$\phi:(Y_U,\omega|_{Y_U})\rightarrow(\bb{R}^{1|1},ds+\lambda d\lambda)$. 
Indeed, suppose $\omega|_{Y_U}=ds_U+\lambda_U d\lambda_U$ for some 
$s_U\in\mc{O}(Y_U)^\ev$, $\lambda_U\in\mc{O}(Y_U)^\od$. 
The map $\phi=(\phi_\r,\phi^*)$ can then be defined by 
$\phi_\r=s_U:U\rightarrow\bb{R}$ and $\phi^*(\lambda)=\lambda_U$. 

Euclidean $(0|1)$- and $(1|1)$-manifolds have the following alternative 
descriptions:

{\lemma \label{01.spin}
A Euclidean $(0|1)$-manifold $(Z,\xi)$ with an orientation on $Z_\r$ is 
equivalent to a $0$-manifold $Z_\r$ with a spin structure. 
}

\proof
Suppose $Z_\r$ has a spin structure, which consists of an orientation 
and a real line bundle $S\rightarrow Z_\r$ with metric. 
We then obtain a Euclidean $(0|1)$-manifold $(Z,\xi)$, with   
$Z=(Z_\r,\Gamma(-,\Lambda^*S))$ and $\xi=\lambda d\lambda$ where 
$\lambda\in\mc{O}(Z)^\od=\Gamma(Z_\r,S)$ is a section of constant length $1$. 
Notice $\xi$ depends only on the metric of $S$ but not on the choice of 
$\lambda$. 
Conversely, given a Euclidean $(0|1)$-manifold $(Z,\xi)$ with 
$\xi=\lambda d\lambda$ for some $\lambda\in\mc{O}(Z)^\od$, we obtain a line 
bundle $S=\coprod_{x\in Z_\r}\mc{O}_Z(\{x\})^\od$ over $Z_\r$ with a metric 
defined by requiring $\lambda$ be of constant length $1$. 
This, together with the orientation on $Z_\r$, gives a spin structure on 
$Z_\r$. 
The two constructions above are inverse to each other. 
\qedsymbol

{\lemma \label{11.spin}
A Euclidean $(1|1)$-manifold $(Y,\omega)$ is equivalent to a $1$-manifold 
$Y_\r$ with a spin structure. 
}

\proof
Suppose $Y_\r$ has a spin structure, which consists of a Riemannian 
metric, an orientation and a $\bb{Z}/2$-bundle $P\rightarrow Y_\r$. 
Let $S^+=P\times_{\bb{Z}/2}\bb{R}$ be the associated positive spinor bundle 
with the canonical metric. 
We first obtain a $(1|1)$-manifold $Y=(Y_\r,\Gamma(-,\Lambda^*S^+))$; 
a metric $\omega$ on $Y$ is then defined by demanding 
$\omega|_{Y_U}=ds_U+\lambda_U d\lambda_U$ in any open set $U\subset Y_\r$ that 
admits a function $s_U$ measuring the arc length in the positive direction and 
a section $\lambda_U\in\Gamma(U,S^+)=\mc{O}(Y_U)^\od$ of constant length $1$. 
Such $\omega$ exists and is unique. 
Conversely, given a Euclidean $(1|1)$-manifold $(Y,\omega)$, pick an open 
cover $\{U_\alpha\}$ of $Y_\r$ and maps 
$f_\alpha:(Y_{U_\alpha},\omega|_{Y_\alpha})\rightarrow 
(\bb{R}^{1|1},ds+\lambda d\lambda)$. 
If $U_\alpha\cap U_\beta\neq\emptyset$, the change of coordinates 
$\phi_{\beta\alpha}=f_\beta f_\alpha^{-1}$ preserves $ds+\lambda d\lambda$. 
By (\ref{Hom.11}), $\phi_{\beta\alpha}^*(ds)=ds$ and 
$\phi_{\beta\alpha}^*(\lambda)=\sigma_{\beta\alpha}\lambda$ for some 
$\sigma_{\beta\alpha}=\pm 1$. 
The nowhere vanishing $1$-forms $\{f_\alpha^*(ds) \}$ consistently 
induce a Riemannian metric and an orientation on $Y_\r$; the cocycle 
$\{\sigma_{\beta\alpha}\}$ defines a $\bb{Z}/2$-bundle over $Y_\r$. 
These constitute a spin structure on $Y_\r$. 
The two constructions above are inverse to each other. 
\qedsymbol

{\it Remark.}
Despite these equivalences, in most of the discussions involving 
Euclidean $(0|1)$- and $(1|1)$-manifolds, we view them as super manifolds. 

Let us interpret functions and $1$-forms on a super manifold $M$ in 
terms of its functor of points: super manifold 
$B\mapsto M(B):=\{\t{maps from }B\t{ to }M\}$. 
Associated to any $a\in\mc{O}(M)$ is a natural transformation 
\beq \label{fcn.B}
a_B:M(B)\rightarrow\mc{O}(B), \quad f\mapsto f^* a. 
\eeq
Also, each $1$-form on $M$ of the form $adb$, where 
$a,b\in\mc{O}(M)$, induces a natural transformation as follows 
\beq \label{1form.B}
(adb)_B:TM(B)\rightarrow\mc{O}(B), \quad 
(f;\delta f)\mapsto(f^*a)\left.\left(\frac{d}{dt}f_t^* b\right)\right|_{t=0}, 
\eeq
where $\{f_t\}\subset M(B)$ is a smooth one-parameter family whose tangent 
at $t=0$ equals $(f;\delta f)$. 
\footnote{ \label{MB.smooth}
In our situations below, $M(B)$ always has the structure of a smooth, 
finite dimensional vector bundle over a smooth manifold, 
hence in particular has a smooth structure. 
}
For instance, the elements 
$\lambda\in\mc{O}(\bb{R}^{0|1})^\od\subset\mc{O}(\bb{R}^{1|1})^\od$, 
$s\in\mc{O}(\bb{R}^{1|1})^\ev$ defined above induce the following bijections 
\beq \label{R01B.R11B}
\lambda_B:\bb{R}^{0|1}(B)\stackrel{\sim}{\rightarrow}\mc{O}(B)^\od, \qquad
(s_B,\lambda_B):\bb{R}^{1|1}(B)\stackrel{\sim}{\rightarrow}
  \mc{O}(B)^\ev\times\mc{O}(B)^\od. 
\eeq
From now on, we will identify $\bb{R}^{0|1}(B)$ with $\mc{O}(B)^\od$ and 
$\bb{R}^{1|1}(B)$ with $\mc{O}(B)^\ev\times\mc{O}(B)^\od$ 
using these bijections. 
The $1$-form $\lambda d\lambda$ on $\bb{R}^{0|1}$ defines the maps
\beqq
(\lambda d\lambda)_B:T\bb{R}^{0|1}(B)\rightarrow\mc{O}(B), \quad 
  (\eta;\delta\eta)\mapsto\eta\delta\eta. 
\eeqq
We have used (\ref{R01B.R11B}), so that $\eta\in\mc{O}(B)^\od$, 
$\delta\eta\in T_\eta\mc{O}(B)^\od\cong\mc{O}(B)^\od$ and 
$\lambda_B(\eta)=\eta$. 
Similarly, the $1$-form $ds+\lambda d\lambda$ on $\bb{R}^{1|1}$ 
defines the maps 
\beqq
(ds+\lambda d\lambda)_B:T\bb{R}^{1|1}(B)\rightarrow\mc{O}(B), \quad 
  (w,\eta;\delta w,\delta\eta)\mapsto\delta w+\eta\delta\eta, 
\eeqq
where $w,\delta w\in\mc{O}(B)^\ev$. 

Given Euclidean $(0|1)$- or $(1|1)$-manifolds $(M_i,\chi_i)$, $i=1,2$, 
we use the notation 
$\phi:\big(M_1(B),(\chi_1)_B\big)\rightarrow\big(M_2(B),(\chi_2)_B\big)$
for a smooth map $\phi:M_1(B)\rightarrow M_2(B)$ that satisfies 
$(\chi_1)_B=(\chi_2)_B\circ T\phi$; 
we write $\t{Hom}\big(M_1(B),(\chi_1)_B;M_2(B),(\chi_2)_B\big)$ 
for the set of such maps. 
For example, keeping (\ref{R01B.R11B}) in mind, we define 
\beqq
&& \epsilon:\bb{R}^{0|1}(B)\rightarrow\bb{R}^{0|1}(B), \quad 
  \epsilon(\eta)=-\eta, \\
\t{or} && \epsilon:\bb{R}^{1|1}(B)\rightarrow\bb{R}^{1|1}(B), \quad
  \epsilon(w,\eta)=(w,-\eta), \\ 
&& \gamma_{z,\theta}:\bb{R}^{0|1}(B)\rightarrow\bb{R}^{1|1}(B), \quad
  \gamma_{z,\theta}(\eta)=(z-\theta\eta,\theta+\eta), \\ 
&& \tau_{z,\theta}:\bb{R}^{1|1}(B)\rightarrow\bb{R}^{1|1}(B), \quad
  \tau_{z,\theta}(w,\eta)=(w+z-\theta\eta,\theta+\eta), 
\eeqq
where $(z,\theta)\in\mc{O}(B)^\ev\times\mc{O}(B)^\od$. 
Then, we have 
\beq
&& \t{Hom}\big(\bb{R}^{0|1}(B),(\lambda d\lambda)_B;\,\bb{R}^{0|1}(B),
  (\lambda d\lambda)_B\big)=
  \{1,\epsilon\},\label{BHom.00} \\
&& \t{Hom}\big(\bb{R}^{0|1}(B),(\lambda d\lambda)_B;\,\bb{R}^{1|1}(B),
  (ds+\lambda d\lambda)_B\big)=
  \{\gamma_{z,\theta},\gamma_{z,\theta}\epsilon\},\label{BHom.01} \\
&& \t{Hom}\big(\bb{R}^{1|1}(B),(ds+\lambda d\lambda)_B;\,\bb{R}^{1|1}(B),
  (ds+\lambda d\lambda)_B\big)=
  \{\tau_{z,\theta},\tau_{z,\theta}\epsilon\}.\label{BHom.11}
\eeq
To obtain, say (\ref{BHom.11}), we suppose a map 
$\tau:\bb{R}^{1|1}(B)\rightarrow\bb{R}^{1|1}(B)$ has the form 
$\tau(w,\eta)=(a(w)+b(w)\eta,c(w)+d(w)\eta)$, apply both sides of 
$(ds+\lambda d\lambda)_B=(ds+\lambda d\lambda)_B\circ T\tau$ to 
$(w,\eta;\delta w,\delta\eta)\in T\bb{R}^{0|1}(B)$, and solve the resulting 
equation for $a(w)$, $b(w)$, $c(w)$ and $d(w)$. 
For any $B$, there is a natural map
\beq \label{Hom.BHom}
\begin{array}{ccc}
\t{Hom}(M_1,\chi_1;M_2,\chi_2)&\rightarrow&
\t{Hom}\big(M_1(B),(\chi_1)_B;M_2(B),(\chi_2)_B\big) \\
\phi &\mapsto& (\phi_B:f\mapsto\phi\circ f).
\end{array}
\eeq 
In particular, (\ref{Hom.00}) $\rightarrow$ (\ref{BHom.00}) is a bijection, 
while (\ref{Hom.01}) $\rightarrow$ (\ref{BHom.01}) is given by 
$\gamma_t\mapsto\gamma_{t,0}$, $\gamma_t\epsilon\mapsto\gamma_{t,0}\epsilon$, 
and (\ref{Hom.11}) $\rightarrow$ (\ref{BHom.11}) by $\tau_t\mapsto\tau_{t,0}$, 
$\tau_t\epsilon\mapsto\tau_{t,0}\epsilon$. 
Although (\ref{Hom.BHom}) is not always surjective, we have the following 
result, as well as lemma \ref{BHom.Hom}. 

{\lemma \label{BHom.red}
Suppose $(M_i,\chi_i)$, $i=1,2$, are Euclidean $(0|1)$- or $(1|1)$-manifolds, 
and $B$ is a super manifold with $B_\r=\{b\}$. 
Given a map
\beqq
\phi:\big(M_1(B),(\chi_1)_B\big)\rightarrow\big(M_2(B),(\chi_2)_B\big), 
\eeqq
there exists a map $\phi_\r:(M_1)_\r\rightarrow(M_2)_\r$ of ordinary 
manifolds so that the diagram 
\beqq
\xymatrix{
  \big(M_1(B),(\chi_1)_B\big)\ar[r]^\phi\ar[d] & 
    \big(M_2(B),(\chi_2)_B\big)\ar[d] \\
  (M_1)_\r\ar[r]^{\phi_\r} & (M_2)_\r \; , 
} 
\eeqq
with the vertical arrows given by $f\mapsto f_\r(b)$, commutes. 
In other words, $\phi(f)_\r=\phi_\r\circ f_\r$, $\forall f\in M_1(B)$. 
Furthermore, $(\phi'\phi)_\r=\phi'_\r\circ\phi_\r$ whenever 
$\phi$ and $\phi'$ are composable.
}

\proof
See appendix \ref{lemmas}.
\qedsymbol

\vspace{0.1in}
\noindent
{\bf Super Euclidean field theories of degree zero.} 
Let $\mc{S}$ be the category whose objects are super manifolds $B$ 
with $B_\r=$ a point, and whose morphisms are maps of super manifolds. 
In this chapter, a \emph{super category} $\mc{C}$ is a functor 
$\mc{C}:\mc{S}^{\t{op}}\rightarrow\t{\bf CAT}$, where $\t{\bf CAT}$ is the 
category of categories and functors. 
\footnote{
A super category may also be defined as a contravariant functor from the 
category of \emph{all} super manifolds to $\t{\bf CAT}$. 
The use of super categories in defining field theories is emphasized by Stolz 
and Teichner. \cite{ST}
}
A \emph{functor of super categories} $F:\mc{C}_1\rightarrow\mc{C}_2$ is a 
natural transformation between $\mc{C}_1$ and $\mc{C}_2$. 
We will now give two examples of super categories. 

The first super category we define is $\mc{S}\t{Hilb}$. 
Consider the category $\mc{S}\t{Hilb}(B)$ associated with an object
$B$ of $\mc{S}$.
The objects of $\mc{S}\t{Hilb}(B)$, which are independent of $B$, are 
real graded Hilbert spaces.
Let $\mc{H}_1$ and $\mc{H}_2$ be two of them. 
The set of morphisms from $\mc{H}_1$ to $\mc{H}_2$ in $\mc{S}\t{Hilb}(B)$ 
is the graded tensor product $\mc{O}(B)\otimes B(\mc{H}_1,\mc{H}_2)$, 
where the grading on $B(\mc{H}_1,\mc{H}_2)$ is induced by those of $\mc{H}_1$ 
and $\mc{H}_2$. 
Composition is defined by composition of operators and 
the algebra structure on $\mc{O}(B)$. 
Given a morphism $f:B\rightarrow B'$ in $\mc{S}$, the functor 
$\mc{S}\t{Hilb}(f):\mc{S}\t{Hilb}(B')\rightarrow
\mc{S}\t{Hilb}(B)$ is defined to be identity on objects and 
$f^*\otimes 1:\mc{O}(B')\otimes B(\mc{H}_1,\mc{H}_2)\rightarrow 
\mc{O}(B)\otimes B(\mc{H}_1,\mc{H}_2)$ on morphisms. 

The second super category $\mc{SEB}^1$ is a `super' version of a 
$(0+1)$-dimensional bordism category. 
Fix an object $B$ of $\mc{S}$ and let $\Lambda=\mc{O}(B)$. 
We now describe the category $\mc{SEB}^1(B)$. 
The objects of $\mc{SEB}^1(B)$ are independent of $B$. 
For each pair of nonnegative integers $m$ and $n$, there is an object 
$Z_{m,n}$, which is the super manifold 
$\underbrace{\bb{R}^{0|1}\sqcup\cdots\sqcup\bb{R}^{0|1}}_{m+n}$ with its first 
$m$ underlying points positively oriented and the rest negative oriented. 
Let $\lambda\in\mc{O}(Z_{m,n})^\od$ be the element whose restriction to each 
$\bb{R}^{0|1}$ equals what we also denote by $\lambda$ in e.g. (\ref{Hom.00}). 
There are two types of morphisms in $\mc{SEB}^1(B)$, 
`$B$-isomorphisms' and equivalence classes of `$B$-bordisms.' 

A \emph{$B$-isomorphism}, always a morphism from some $Z_{m,n}$ to 
itself, is a map 
\beqq
\xymatrix{
 \big(Z_{m,n}(B),(\lambda d\lambda)_{B}\big)\ar[r]^\mu & 
 \big(Z_{m,n}(B),(\lambda d\lambda)_{B}\big), 
}
\eeqq
such that $\mu_\r$ (see lemma \ref{BHom.red}) is bijective and 
orientation-preserving. 
For example, according to (\ref{BHom.00}), the $B$-isomorphisms on 
$Z_{1,0}$ are $1$ and $\epsilon$. 

A \emph{$B$-bordism} from $Z_{m,n}$ to $Z_{p,q}$ is a $4$-tuple 
$(Y,\omega,\alpha,\beta)$, consisting of a Euclidean $(1|1)$-manifold 
$(Y,\omega)$ with $Y_\r$ compact, and two maps 
\beqq
\xymatrix{
 & \big(Y(B),\omega_{B}\big) & \\
 \big(Z_{m,n}(B),(\lambda d\lambda)_{B}\big)\ar[ru]^\alpha & & 
 \big(Z_{p,q}(B),(\lambda d\lambda)_{B}\big). \ar[lu]_\beta 
} 
\eeqq 
By lemma \ref{11.spin}, $\omega$ induces an orientation on $Y_\r$, which in 
turn induces an orientation on $\partial Y_\r$. 
We require $\alpha_\r\cup\beta_\r$ map $(Z_1)_\r\cup(Z_2)_\r$ bijectively 
onto $\partial Y_\r$ with $\alpha_\r$ orientation-reversing and $\beta_\r$ 
orientation-preserving. 
(See lemma \ref{BHom.red}.) 
Two $B$-bordisms $(Y,\omega,\alpha,\beta)$ and 
$(Y',\omega',\alpha',\beta')$ are \emph{equivalent} if there is 
an invertible map $\tau:\big(Y(B),\omega_{B}\big)\rightarrow
\big(Y'(B),\omega'_{B}\big)$ making the following diagram 
\beq \label{equiv.sbordisms}
\xymatrix{
 & Y(B)\ar[dd]_{\tau} & \\
 Z_{m,n}(B)\ar[ru]^\alpha \ar[rd]_{\alpha'} & & 
  Z_{p,q}(B)\ar[lu]_\beta \ar[ld]^{\beta'} \\
 & Y'(B) &
}
\eeq
commute. 
We denote the equivalence class of $(Y,\omega,\alpha,\beta)$ by 
$[Y,\omega,\alpha,\beta]$ and \emph{identify equivalent $B$-bordisms 
as the same morphism} in $\mc{SEB}^1(B)$. 

As an example, let us find all $B$-bordisms between $Z_{1,0}$ and 
itself of the form $(\bb{R}^{1|1}_{[0,t]},ds+\lambda d\lambda,\alpha,\beta)$, 
where $t>0$, and classify them up to equivalence. 
Notice that 
\beqq
\bb{R}^{1|1}_{[0,t]}(B)=\Lambda^\ev_{[0,t]}\times\Lambda^\od, 
\eeqq
where, for any $S\subset\bb{R}$, we denote by 
$\Lambda^\ev_S\subset\Lambda^\ev$ the subset of elements 
whose $\Lambda^0$-parts lie in $S$. 
By (\ref{BHom.01}), we have $\alpha=\gamma_{z_1,\theta_1}$ or 
$\gamma_{z_1,\theta_1}\epsilon$ and 
$\beta=\gamma_{z_2,\theta_2}$ or $\gamma_{z_2,\theta_2}\epsilon$, 
where $z_1\in\Lambda^\ev_{\{0\}}$, $z_2\in\Lambda^\ev_{\{t\}}$, 
$\theta_1,\theta_2\in\Lambda^\od$. 
Let $\tau$ be as in (\ref{equiv.sbordisms}) with 
$(Y,\omega)=(Y',\omega')=(\bb{R}^{1|1}_{[0,t]},ds+\lambda d\lambda)$. 
By (\ref{BHom.11}), $\tau=\tau_{z,\theta}$ or $\tau_{z,\theta}\epsilon$ 
with $(z,\theta)\in\Lambda^\ev_{\{0\}}\times\Lambda^\od$. 
Observe that, given $\alpha=\gamma_{z_1,\theta_1}$ (resp. 
$\gamma_{z_1,\theta_1}\epsilon$), there is a unique $\tau$, namely 
$\tau_{-z_1,-\theta_1}$ (resp. $\tau_{-z_1,-\theta_1}\epsilon$), such that 
$\alpha'=\tau\alpha=\gamma_{0,0}$. 
Therefore, up to equivalence, we may assume $\alpha=\gamma_{0,0}$, while 
$\beta$ can be any $\gamma_{z_2,\theta_2}$ or $\gamma_{z_2,\theta_2}\epsilon$. 
We define 
\beq \label{I.z.theta}
I_{z,\theta} &=&
  [\,\bb{R}^{1|1}_{[0,z^0]},ds+\lambda d\lambda,
  \gamma_{0,0},\gamma_{z,\theta}\,] \\ 
\epsilon I_{z,\theta} &=&
  [\,\bb{R}^{1|1}_{[0,z^0]},ds+\lambda d\lambda,
  \gamma_{0,0},\gamma_{z,\theta}\epsilon\,] \nonumber
\eeq
where now $z\in\Lambda^\ev_{(0,\infty)}$, $\theta\in\Lambda^\od$ and 
$z^0$ is the $\Lambda^0$-part of $z$.
According to the definition of composition to be given below, 
$\epsilon I_{z,\theta}$ is indeed the composition of 
$I_{z,\theta}$ and $\epsilon$. 

To finish the definition of $\mc{SEB}^1(B)$, it remains to describe 
how morphisms compose ($\mu$, $\mu'$ are $B$-isomorphisms):
\beqq
\begin{array}{lcl}
\mu'\circ \mu &=& \mu'\mu \\
\mu'\circ[Y,\omega,\alpha,\beta] &=& [Y,\omega,\alpha,\beta{\mu'}^{-1}] \\ 
\left[Y',\omega',\alpha',\beta'\right]\circ\mu &=& 
  [Y',\omega',\alpha'\mu,\beta'] \\ 
\left[Y',\omega',\alpha',\beta'\right]\circ[Y,\omega,\alpha,\beta] &=& 
  [Y'',\omega'',\iota\alpha,\iota'\beta']. 
\end{array}
\eeqq
The last case is described in a commutative diagram (the $1$-forms 
suppressed) 
\beq \label{sbordism.compose.diagram}
\xymatrix{
 && Y''(B) && \\
 & Y(B)\ar[ru]^\iota && Y'(B)\ar[lu]_{\iota'} & \\
 Z_{m,n}(B)\ar[ru]^\alpha && Z_{p,q}(B)\ar[lu]_\beta\ar[ru]^{\alpha'} && 
  Z_{r,s}(B)\ar[lu]_{\beta'},
} 
\eeq
where $\iota$ and $\iota'$ are injections making the square in the above 
diagram a pushout. 
It is straightforward to check that 
$(Y'',\omega'',\iota\alpha,\iota'\beta')$ 
is well-defined up to equivalence. 

Let us compute the relations among various endomorphisms of $Z_{1,0}$. 
Clearly, $\epsilon^2=1$. 
Let $Y=\bb{R}^{1|1}_{[0,z^0]}$ and $\omega=ds+\lambda d\lambda$. 
Consider
\beq \label{epsilon.I.epsilon}
\epsilon I_{z,\theta}\epsilon
&=& \epsilon\circ[\,Y,\omega,\gamma_{0,0},\gamma_{z,\theta}\,]\circ\epsilon 
  = [\,Y,\omega,\gamma_{0,0}\epsilon,\gamma_{z,\theta}\epsilon\,] \nonumber\\
&=& [\,Y,\omega,\epsilon\gamma_{0,0}\epsilon,
  \epsilon\gamma_{z,\theta}\epsilon\,]
  =[\,Y,\omega,\gamma_{0,0},\gamma_{z,-\theta}\,] \nonumber\\
&=& I_{z,-\theta}\; ,
\eeq
where the third equality follows from the definition of equivalent 
$B$-bordisms and the fourth equality is an easy computation. 
On the other hand, the commutative diagram 
\beqq 
\xymatrix{
 && \Lambda^\ev_{[0,z^0+z'^0]}\times\Lambda^\od && \\
 & \Lambda^\ev_{[0,z^0]}\times\Lambda^\od \ar[ru]^{\tau_{0,0}} 
   && \Lambda^\ev_{[0,z'^0]}\times\Lambda^\od \ar[lu]_{\tau_{z,\theta}} & \\
 \Lambda^\od \ar[ru]^{\gamma_{0,0}} 
   && \Lambda^\od \ar[lu]_{\gamma_{z,\theta}} \ar[ru]^{\gamma_{0,0}} 
   && \Lambda^\od \ar[lu]_{\gamma_{z',\theta'}} \; ,
}
\eeqq
together with the facts $\tau_{0,0}\gamma_{0,0}=\gamma_{0,0}$ and 
$\tau_{z,\theta}\gamma_{z',\theta'}= 
\gamma_{z+z'-\theta\theta',\theta+\theta'}$, shows 
\beq \label{sbordism.compose}
I_{z',\theta'}I_{z,\theta}=I_{z+z'-\theta\theta',\theta+\theta'}.
\eeq
This defines the structure of a semigroup on 
$\{I_{z,\theta}\}\cong\Lambda^\ev_{(0,\infty)}\times\Lambda^\od$. 
\footnote{\label{ssemigroup}
Furthermore, as $B$ varies, we obtain a contravariant functor from $\mc{S}$ 
to the category of semigroups, i.e. a \emph{super semigroup}. 
\cite{QFS.susy}
}

Finally, given a morphism $f:B'\rightarrow B$ in $\mc{S}$, we describe 
the functor 
$\mc{SEB}^1(f):\mc{SEB}^1(B)\rightarrow\mc{SEB}^1(B')$. 
It is identity on objects and $B$-isomorphisms. 
Notice that any $B$-isomorphism is an orientation-preserving 
permutation of points together with a choice of $1$ or $\epsilon$ for each 
point. 
The set of $B$-isomorphisms is thus `the same' for all $B$ 
and it makes sense to say $\mc{SEB}^1(f)$ is identity on them. 
As to $B$-bordisms, we first define for any 
$(z,\theta)\in\Lambda^\ev_{(0,\infty)}\times\Lambda^\od$ the following: 
\beq
\bar{I}_{z,\theta}
  &=& [\,\bb{R}^{1|1}_{[0,z^0]},ds+\lambda d\lambda,
  \gamma_{z,\theta},\gamma_{0,0}\,] \; , \nonumber \\
LI_{z,\theta}
  &=& [\,\bb{R}^{1|1}_{[0,z^0]},ds+\lambda d\lambda,
  \gamma_{0,0}\sqcup\gamma_{z,\theta},-\,] \; , \label{Ibar.LI.RI} \\
RI_{z,\theta}
  &=& [\,\bb{R}^{1|1}_{[0,z^0]},ds+\lambda d\lambda,
  -,\gamma_{z,\theta}\sqcup\gamma_{0,0}\,] \; , \nonumber
\eeq
which are respectively $B$-bordisms from $Z_{0,1}$ to $Z_{0,1}$, from 
$Z_{1,1}$ to $Z_{0,0}$, and from $Z_{0,0}$ to $Z_{1,1}$. 
A general $B$-bordism is a disjoint union of ones of the form 
$I_{z,\theta}$, $\bar{I}_{z,\theta}$, $LI_{z,\theta}$, $RI_{z,\theta}$ and 
$[Y,\omega,-,-]$ where $Y_\r=S^1$. 
The map $\mc{SEB}^1(f)$ on $B$-bordisms is defined by 
$I_{z,\theta}\mapsto I_{f^*z,f^*\theta}$ (and similarly for 
$\bar{I}_{z,\theta}$, $LI_{z,\theta}$ and $RI_{z,\theta}$), 
$[Y,\omega,-,-]\mapsto [Y,\omega,-,-]$, and the requirement that it preserves 
disjoint unions. 

\defn \label{defn.SEFT0}
A \emph{one-dimensional super Euclidean field theory}, or \emph{$1$-SEFT}, of 
degree $0$ is a functor $E:\mc{SEB}^1\rightarrow\mc{S}\t{Hilb}$ of super 
categories. 
Let $E_B:\mc{SEB}^1(B)\rightarrow\mc{S}\t{Hilb}(B)$ be the (ordinary) functor
associated with each object $B$ of $\mc{S}$. 
Since on objects $E_B$ does not depend on $B$, we simply write $E(Z_{m,n})$. 
Let $\mc{H}=E(Z_{1,0})$. 
A priori $E_B(\epsilon)\in\mc{O}(B)\otimes B(\mc{H})$, 
but since $(f^*\otimes 1)E_{B'}(\epsilon)=E_B(\epsilon)$ 
for every morphism $f:B\rightarrow B'$ in $\mc{S}$, each 
$E_B(\epsilon)$ must in fact be an element of $B(\mc{H})$ independent of $B$. 
Abusing notations, we denote that element also by $\epsilon$. 
Now, fix $B$ and let $\Lambda=\mc{O}(B)$. 
The functor $E_B$ is required to satisfy the following assumptions:
\begin{itemize}
\item[(i)]
$E_{B}$ sends $Z_{0,0}$ to $\bb{R}$ and disjoint unions of objects 
(resp. morphisms) to tensor products of graded Hilbert spaces 
(resp. operators). 
More precisely, tensor products of operators are maps of the form 
given by multiplication in $\Lambda$: 
\beqq
\big(\Lambda\otimes B(\mc{H}_1,\mc{H}_2)\big)\otimes
\big(\Lambda\otimes B(\mc{H}'_1,\mc{H}'_2)\big)\rightarrow 
\Lambda\otimes B(\mc{H}_1\otimes\mc{H}'_1,\mc{H}_2\otimes\mc{H}'_2).
\eeqq
\item[(ii)] 
$E_{B}$ sends objects (resp. $B$-bordisms) with opposite orientations to 
the same graded Hilbert space (resp. operator), 
e.g. $E(Z_{n,m})=E(Z_{m,n})$, 
$E_{B}(\bar{I}_{z,\theta})=E_{B}(I_{z,\theta})$. 
\footnote{  \label{cplx.coeff}
For complex coefficients, opposite orientation corresponds to complex 
conjugation. 
}
\item[(iii)]
If $(Y,\omega,\alpha,\beta)$ is a $B$-bordism from $Z_{m,n}$ to $Z_{p,q}$, 
then $(Y,\omega,-,\alpha\sqcup\beta)$ is a $B$-bordism from $Z_{0,0}$ to 
$Z_{n+p,m+q}$. 
The map 
\beqq
\t{Ad}_{E(Z_{m,n})}:\Lambda\otimes E(Z_{n,m})\otimes E(Z_{p,q})\rightarrow 
\Lambda\otimes B\big(E(Z_{m,n}),E(Z_{p,q})\big)
\eeqq
given by taking adjoints sends $E_{B}([Y,\omega,-,\alpha\sqcup\beta])$ to 
$E_{B}([Y,\omega,\alpha,\beta])$. 
We have used the fact $E(Z_{n,m})$ is canonically isomorphic to the dual 
of $E(Z_{m,n})$. 
\footnote{
This is the case for both real and complex coefficients. 
See footnote \ref{cplx.coeff}. 
}
\item[(iv)]
The operator $\epsilon=E_{B}(\epsilon)$ on $\mc{H}$ is the grading 
homomorphism, i.e $\epsilon=1$ on $\mc{H}^\ev$ and $\epsilon=-1$ on 
$\mc{H}^\od$. 
\item[(v)]
$E_{B}(I_{z,\theta})$ is analytic in 
$(z,\theta)\in\Lambda^\ev_{(0,\infty)}\times\Lambda^\od$. 
\end{itemize}

Let us examine the implications of assumptions (i)-(v). 
For objects, we have $E(Z_{m,n})=\mc{H}^{\otimes(m+n)}$ by (i) and (ii). 
Also, (i), (ii) and (iv) determine the action of $E_B$ on $B$-isomorphisms. 
Since every $B$-bordism is a disjoint union of ones of the forms 
$I_{z,\theta}$, $\bar{I}_{z,\theta}$, $LI_{z,\theta}$, $RI_{z,\theta}$ and 
$[Y,\omega,-,-]$ with $Y_\r=S^1$, and $[Y,\omega,-,-]$ is a composition of 
the others, it suffices by (i) and functoriality to understand 
$E_{B}(I_{z,\theta})$, $E_{B}(\bar{I}_{z,\theta})$, 
$E_{B}(LI_{z,\theta})$ and $E_{B}(RI_{z,\theta})$. 
According to (iii), the three adjoint maps 
\beqq
\xymatrix{
  & \Lambda\otimes\mc{H}^{\otimes 2} \ar[ld]_{\t{Ad}^{(1)}_{\mc{H}}} 
    \ar[d]_{\t{Ad}^{(2)}_{\mc{H}}} \ar[rd]^{\t{Ad}_{\mc{H}^{\otimes 2}}} & \\
  \Lambda\otimes B(\mc{H}) \ar[r]^{(\cdot)^T} 
    & \Lambda\otimes B(\mc{H}) \ar[l] 
    & \Lambda\otimes(\mc{H}^{\otimes 2})^\vee 
}
\eeqq
send $E_{B}(RI_{z,\theta})$ to $E_{B}(I_{z,\theta})$, 
$E_{B}(\bar{I}_{z,\theta})$ and $E_{B}(LI_{z,\theta})$, where 
$\t{Ad}^{(1)}_{\mc{H}}$ and $\t{Ad}^{(2)}_{\mc{H}}$ correspond to the 
two different factors of $\mc{H}^{\otimes 2}$. 
Any one of the four operators, say $E_{B}(I_{z,\theta})$, determines 
the others because the adjoint maps are injective. 
Being in the image of $\t{Ad}^{(1)}_{\mc{H}}$ precisely means 
$E_{B}(I_{z,\theta})$ is Hilbert-Schmidt. 
Also, the images under $\t{Ad}^{(1)}_{\mc{H}}$ and $\t{Ad}^{(2)}_{\mc{H}}$ of 
the same element are transpose of one another, hence 
$E_{B}(\bar{I}_{z,\theta})=E_{B}(I_{z,\theta})^T$. 
It then follows by (ii) that $E_{B}(I_{z,\theta})$ is self-adjoint. 

The relations among the morphisms in $\mc{SEB}^1(B)$ impose further 
restrictions on the operators $E_{B}(I_{z,\theta})$. 
Let $E_{B}(I_{z,\theta})=X(z)+\theta Y(z)$, where 
$X(z),Y(z)\in\Lambda\otimes B(\mc{H})$ are self-adjoint and Hilbert-Schmidt 
(hence compact). 
The relation in (\ref{epsilon.I.epsilon}) and (iv) imply $X(z)$ are even and 
$Y(z)$ are odd. 
The composition law (\ref{sbordism.compose}) translates into the equations 
\beq
& X(z)X(z')-\theta\theta'Y(z)Y(z')=X(z+z'+\theta\theta'), 
\label{XY1} \\
& \theta'X(z)Y(z')+\theta Y(z)X(z')=(\theta+\theta')Y(z+z'+\theta\theta').
\label{XY2}
\eeq
By (v), $X(z)$ and $Y(z)$ are determined by their values for, say positive 
real $z$. 
Let $t,t'$ denote positive real numbers. 
Arguing in the same way we did for $\epsilon$, we see that $X(t)$ and $Y(t)$ 
are elements of $B(\mc{H})$ independent of $B$. 
We claim that all $X(t)$ and $Y(t)$ commute. 
This is true among $X(t)$ by (\ref{XY1}) with $\theta=0$. 
Using (\ref{XY1}) again and the commutativity among $X(t)$, we see that 
$Y(t)Y(t')-Y(t')Y(t)\in B(\mc{H})$ is annihilated by any $\theta\theta'$, 
and so must itself vanish. 
Finally, using (\ref{XY2}) with $(z,\theta,z',\theta')=(t,0,t',\theta')$ and 
$(t',\theta,t,0)$, we obtain $X(t)Y(t')=Y(t')X(t)$. 
This proves our claim. 
By the spectral theorem, the commuting family of self-adjoint compact 
operators $\{X(t),Y(t)\}$ has a simultaneous eigenspace decomposition. 
Since $X(t)X(t')=X(t+t')$ for all $t,t'>0$, there is an even, self-adjoint 
operator $H$ with compact resolvent defined on a graded subspace 
$\mc{H}'\subset\mc{H}$, such that $X(t)=e^{-tH}$ on $\mc{H}'$ and $X(t)=0$ on 
$\mc{H}'^\perp$. 
Let $Q=Y(1)e^H$, which is an odd, self-adjoint operator on $\mc{H}'$ commuting 
with $H$. 
For $t>1$, (\ref{XY2}) implies $Y(t)=Y(1)X(t-1)=Qe^{-tH}$ on $\mc{H}'$, and 
$Y(t)=0$ on $\mc{H}'^\perp$; this in fact holds for all $t>0$ by analyticity. 
Equation (\ref{XY1}), restricted to $\mc{H}'$, now becomes 
$e^{-(z+z')H}(1-\theta\theta'Q^2)=e^{-(z+z'+\theta\theta')H}$, which implies 
$Q^2=H$. 
Therefore, we have 
\beq \label{E.I.z.theta}
E_{B}(I_{z,\theta})=X(z)+\theta Y(z)=
\begin{cases}
  e^{-z Q^2}+\theta Qe^{-zQ^2} & \t{on }\mc{H}' \\
  0 & \t{on }\mc{H}'^\perp
\end{cases}.
\eeq
\footnote{
More conceptually, the operators $E_{B}(I_{z,\theta})$ are determined by a 
single one because the Lie algebra of the super semigroup 
described in footnote \ref{ssemigroup} is generated by a single element. 
\cite{QFS.susy}
}
Conversely, (\ref{E.I.z.theta}) always satisfies 
(\ref{XY1}) and (\ref{XY2}). 

To summarize, a $1$-SEFT $E$ of degree $0$ is determined by a real graded 
Hilbert space $\mc{H}=E(Z_{1,0})$ and an odd, self-adjoint operator $Q$ with 
compact resolvent defined on a graded subspace of $\mc{H}$. 
This operator is referred to as the \emph{generator} of $E$. 

\vspace{0.1in}
\noindent
{\bf Super Euclidean field theories of general degrees.} 
Given a Euclidean $(0|1)$-manifold $(Z,\xi)$ with $Z_\r$ oriented, 
let $C(Z,\xi)$ be the Clifford algebra on $\mc{O}(Z)^\od$ 
with respect to the quadratic form 
$\lambda'\mapsto\sum_{x\in Z_\r}\t{sgn}(x)[\lambda'(x)/\lambda(x)]^2$, 
where $\lambda\in\mc{O}(Z)^\od$ such that $\xi=\lambda d\lambda$, and 
$\t{sgn}(x)=\pm 1$ is the orientation of $x$. 
Notice the quadratic form depends only on $\xi$ but not on the choice of 
$\lambda$. 
A map $\mu:(Z,\xi)\rightarrow(Z',\xi')$ between two such $(0|1)$-manifolds, 
with $\mu_\r$ orientation-preserving, induces an algebra homomorphism 
$\mu^*:C(Z',\xi')\rightarrow C(Z,\xi)$. 
By lemma \ref{01.spin}, $(Z,\xi)$ induces a spin structure on $Z_\r$, which 
includes a real line bundle $S\rightarrow Z_\r$ with metric $||\cdot||$; 
$C(Z,\xi)$ can then be described as in \cite{ST}, i.e. the Clifford algebra 
on $\Gamma(Z_\r,S)$ with respect to the quadratic form 
$\lambda'\mapsto\sum_{x\in Z_\r}\t{sgn}(x)||\lambda'(x)||^2$. 


For any Euclidean $(1|1)$-manifold $(Y,\omega)$, define a real vector space 
$L(Y,\omega)$ as follows. 
It is a subspace of $\mc{O}(Y)^\od$. 
An element $\psi\in\mc{O}(Y)^\od$ is in $L(Y,\omega)$ iff, for any open set 
$U\subset Y_\r$ with $\omega|_{Y_U}=ds_U+\lambda_U d\lambda_U$ for some 
$s_U\in\mc{O}(Y_U)^\ev$, $\lambda_U\in\mc{O}(Y_U)^\od$, 
the restriction $\psi|_U$ is a constant multiple of $\lambda_U$. 
A map $\tau:(Y,\omega)\rightarrow(Y',\omega')$ induces a linear map 
$\tau^*:L(Y',\omega')\rightarrow L(Y,\omega)$. 
By lemma \ref{11.spin}, $(Y,\omega)$ induces a spin structure on $Y_\r$; 
$L(Y,\omega)$ can then be described as in \cite{ST}, i.e. the space of 
sections of the positive spinor bundle $S^+$ satisfying the Dirac equation. 

The \emph{fermionic Fock space} of a Euclidean $(1|1)$-manifold $(Y,\omega)$ 
is 
\beqq
F(Y,\omega)=\Lambda^{\t{top}}L(Y_1,\omega|_{Y_1})^\vee\otimes
  \Lambda^*\big(L(Y_2,\omega|_{Y_2})^\vee\big), 
\eeqq
where $Y_1$ is the union of the closed components of $Y$ and $Y_2=Y-Y_1$. 
The Fock space depends covariantly with $(Y,\omega)$. 
Let us recall a number of its properties from \cite{ST}. 

Suppose $(Z_i,\xi_i)$, $i=1,2$, are two Euclidean $(0|1)$-manifolds with 
$(Z_i)_\r$ oriented. 
A \emph{bordism} $(Y,\omega,\alpha,\beta)$ from $(Z_1,\xi_1)$ to 
$(Z_2,\xi_2)$ consists of a Euclidean $(1|1)$-manifold $(Y,\omega)$ with 
$Y_\r$ compact, and two maps 
\beqq
\xymatrix{
 & (Y,\omega) & \\
 (Z_1,\xi_1)\ar[ru]^\alpha & & (Z_2,\xi_2), \ar[lu]_\beta 
} 
\eeqq
such that $\alpha_\r\cup\beta_\r$ maps $(Z_1)_\r\cup(Z_2)_\r$ bijectively 
onto $\partial Y_\r$ with $\alpha_\r$ orientation-reversing and $\beta_\r$ 
orientation-preserving. 
(Recall, by lemma \ref{11.spin}, the metric induces an orientation on 
$Y_\r$, and hence on $\partial Y_\r$.) 
\emph{Each bordism $(Y,\omega,\alpha,\beta)$ induces 
a $C(Z_2,\xi_2)$-$C(Z_1,\xi_1)$-action on $F(Y,\omega)$}. 
\footnote{
This action is an example of the fermionic Fock representation. 
}
This action has the following properties. 
Suppose $(Z_3,\xi_3)$ is another Euclidean $(0|1)$-manifold with $(Z_3)_\r$ 
oriented and $\mu:(Z_3,\xi_3)\rightarrow(Z_2,\xi_2)$ is 
an invertible, orientation-preserving map. 
The $C(Z_3,\xi_3)$-$C(Z_1,\xi_1)$-action on $F(Y,\omega)$ induced by the 
new bordism $(Y,\omega,\alpha,\beta\mu)$ pulls back to the above 
$C(Z_2,\xi_2)$-$C(Z_1,\xi_1)$-action along the algebra isomorphism 
$\mu^*:C(Z_2,\xi_2)\rightarrow C(Z_3,\xi_3)$. 
On the other hand, 
any invertible map $\tau:(Y,\omega)\rightarrow(Y',\omega')$ defines 
another bordism $(Y',\omega',\tau\alpha,\tau\beta)$, as well as 
an isomorphism $\tau_*:F(Y,\omega)\rightarrow F(Y',\omega')$ of the induced 
$C(Z_2,\xi_2)$-$C(Z_1,\xi_1)$-bimodules. 


\emph{Composition of bordisms gives rise to tensor product of bimodules}. 
More precisely, suppose we have a diagram 
\beqq
\xymatrix{
 && (Y'',\omega'') && \\
 & (Y,\omega)\ar[ru]^\iota && (Y',\omega')\ar[lu]_{\iota'} & \\
 (Z_1,\xi_1)\ar[ru]^\alpha && 
   (Z_2,\xi_2)\ar[lu]_\beta\ar[ru]^{\alpha'} && 
   (Z_3,\xi_3)\ar[lu]_{\beta'},
}
\eeqq
where $(Y,\omega,\alpha,\beta)$ and $(Y',\omega',\alpha',\beta')$ are bordisms 
and the square is a pushout. 
Then $(Y'',\omega'',\iota\alpha,\iota'\beta')$ is also a bordism, and there 
is a canonical isomorphism 
\beq \label{Fock.gluing}
F(Y',\omega')\otimes_{C(Z_2,\xi_2)}F(Y,\omega)\cong F(Y'',\omega'')
\eeq
of $C(Z_3,\xi_3)$-$C(Z_1,\xi_1)$-bimodules. 

For each $p,q\geq 0$, let $\lambda\in\mc{O}(Z_{p,q})^\od$ be the canonical 
element defined earlier. 
Notice that 
\beqq
C(Z_{p,q}):=C(Z_{p,q},\lambda d\lambda)\cong Cl_{p,q}.
\eeqq 
If $\epsilon:(Z_{p,q},\lambda d\lambda)\rightarrow(Z_{p,q},\lambda d\lambda)$ 
is given by $\lambda\mapsto-\lambda$, 
the induced map $\epsilon^*$ on $C(Z_{p,q})$ is precisely the grading 
homomorphism. 
For $t>0$, we use the notations
\beqq
\begin{array}{c}
F(\bb{R}^{1|1}_{[0,t]}):=F(\bb{R}^{1|1}_{[0,t]},ds+\lambda d\lambda), \\
\Omega:=1\in
  \Lambda^0\big(L(\bb{R}^{1|1}_{[0,t]},ds+\lambda d\lambda)^\vee\big)
  \subset F(\bb{R}^{1|1}_{[0,t]}).
\end{array}
\eeqq 
\emph{The element $\Omega$ has the following properties}. 
The self-maps of $(\bb{R}^{1|1}_{[0,t]},ds+\lambda d\lambda)$, which consist 
of $\tau_0=1$ and $\tau_0\epsilon=\epsilon$ in (\ref{Hom.11}), preserve 
$\Omega$. 
For any appropriate bordism that induces the structure of a 
$C(Z_{1,0})$-$C(Z_{1,0})$-bimodule on $F(\bb{R}^{1|1}_{[0,t]})$, the 
bimodule is generated by $\Omega$ (hence irreducible) and has the property 
that $\lambda\Omega=\Omega\lambda$, where 
$\lambda\in\mc{O}(Z_{1,0})^\od\subset C(Z_{1,0})$. 
Also, any isomorphism 
$F(\bb{R}^{1|1}_{[0,s]})\otimes_{C(Z_{1,0})}F(\bb{R}^{1|1}_{[0,t]})
\cong F(\bb{R}^{1|1}_{[0,s+t]})$ 
corresponding to gluing bordisms sends $\Omega\otimes\Omega$ to $\Omega$.

{\lemma \label{BHom.Hom}
Suppose $(M_i,\chi_i)$, $i=1,2$, are Euclidean $(0|1)$- or $(1|1)$-manifolds 
and $B$ is a super manifold with $B_\r=$ a point. 
There is a map 
\beqq
\ell:\t{Hom}\big(M_1(B),(\chi_1)_B;M_2(B),(\chi_2)_B\big)\rightarrow
\t{Hom}(M_1,\chi_1;M_2,\chi_2),
\eeqq
such that $\ell(\phi)_\r=\phi_\r$ (see lemma \ref{BHom.red}) and 
it is a left inverse of (\ref{Hom.BHom}). 
Also, $\ell(\phi'\phi)=\ell(\phi')\ell(\phi)$ whenever $\phi$ and $\phi'$ are 
composable. 
}

\proof 
See appendix \ref{lemmas}.
\qedsymbol

This technical result allows us to replace bordisms with 
$B$-bordisms in the above discussion of the Clifford algebras $C(Z,\xi)$ 
and the Fock spaces $F(Y,\omega)$. 
\footnote{
A more desirable approach would be to find appropriate functor-of-points 
interpretations of $C(Z,\xi)$ and $F(Y,\omega)$. 
However, the author has not been able to achieve that, and resorts to 
this ad hoc approach. 
}
For any $B$-isomorphism $\mu$, say on $Z_{p,q}$, the map $\ell(\mu)$ 
induces an automorphism $\ell(\mu)^*$ on $C(Z_{p,q})$. 
Given a $B$-bordism $(Y,\omega,\alpha,\beta)$, 
it follows from $\ell(\alpha)_\r=\alpha_\r$ and $\ell(\beta)_\r=\beta_\r$ 
that $(Y,\omega,\ell(\alpha),\ell(\beta))$ is a 
bordism, which induces a bimodule structure on $F(Y,\omega)$. 
That $\ell$ respects compositions implies that 
an equivalence of $B$-bordisms yields an isomorphism of bimodules. 
Similarly, composition of $B$-bordisms gives rise to 
composition of bordisms, and hence tensor product of bimodules. 

Now, we define a super category $\mc{SEB}^1_n$ for each $n\in\bb{Z}$. 
Let us first fix an object $B$ of $\mc{S}$ and describe the category 
$\mc{SEB}^1_n(B)$. 
The objects of $\mc{SEB}^1_n(B)$ are $Z_{p,q}$, $p,q\geq 0$. 
There are two types of morphisms, `decorated $B$-isomorphisms' and 
equivalence classes of `decorated $B$-bordisms.' 
A \emph{decorated $B$-isomorphism} $(\mu,c)$ consists of a 
$B$-isomorphism $\mu$ from some $Z_{p,q}$ to itself, and an element 
$c\in C(Z_{p,q})^{\otimes -n}$. 
\footnote{ \label{neg.power.tensor}
If $n>0$, we use the notation $C(Z_{p,q})^{\otimes -n}$ for the opposite 
algebra of $C(Z_{p,q})^{\otimes n}$ and also $F(Y,\omega)^{\otimes -n}$ for 
the opposite bimodule of $F(Y,\omega)^{\otimes n}$. 
See \cite{ST}. 
}
One should think of $c$ as coming from the first $Z_{p,q}$. 
A \emph{decorated $B$-bordism} is a $5$-tuple 
$(Y,\omega,\alpha,\beta,\Psi)$, where the first four components comprise a 
$B$-bordism and $\Psi\in F(Y,\omega)^{\otimes -n}$. 
This $5$-tuple is \emph{equivalent} to precisely those of the form 
$\big(Y',\omega',\tau\alpha,\tau\beta,\ell(\tau)_*\Psi\big)$, for some 
invertible map $\tau:(Y(B),\omega_B)\rightarrow(Y'(B),\omega'_B)$. 
Denote the equivalence class by $[Y,\omega,\alpha,\beta,\Psi]$. 
Composition in $\mc{SEB}^1_n(B)$ is defined as follows:
\beqq
\begin{array}{lcl}
(\mu',c')\circ(\mu,c) &=& (\mu'\mu,\ell(\mu)^*c'\cdot c), \\
(\mu',c')\circ[Y,\omega,\alpha,\beta,\Psi] &=& 
  [Y,\omega,\alpha,\beta{\mu'}^{-1},(\mu'^{-1})^*c'\Psi], \\ 
\left[Y',\omega',\alpha',\beta',\Psi'\right]\circ(\mu,c) &=& 
  [Y',\omega',\alpha'\mu,\beta',\Psi'c], \\ 
\left[Y',\omega',\alpha',\beta',\Psi'\right]\circ[Y,\omega,\alpha,\beta,\Psi] 
  &=& [Y'',\omega'',\iota\alpha,\iota'\beta',\Psi'']. 
\end{array}
\eeqq
For the last case, $(Y'',\omega'',\iota\alpha,\iota'\beta')$ is 
as in (\ref{sbordism.compose.diagram}), and $\Psi''$ is the image of 
$\Psi'\otimes\Psi$ under the induced isomorphism (\ref{Fock.gluing}). 
Also, if a decorated $B$-bordism is written as 
$(Y,\omega,\alpha,\beta,c_1\Psi c_2)$, where $\Psi\in F(Y,\omega)$ and 
$c_1$, $c_2$ are Clifford algebra elements, 
$c_1\Psi c_2$ is defined using the bimodule structure on $F(Y,\omega)$ 
induced by $(Y,\omega,\alpha,\beta)$. 

For any morphism $f:B'\rightarrow B$ in $\mc{S}$, 
the functor $\mc{SEB}^1_n(f)$ is 
identity on objects and on decorated $B$-isomorphisms. 
As to decorated $B$-bordisms, we first define, 
for $(z,\theta)\in\mc{O}(B)^\ev\times\mc{O}(B)^\od$ and 
$\Psi\in F(\bb{R}^{1|1}_{[0,z^0]},ds+\lambda d\lambda)$, 
\beqq
I_{z,\theta,\Psi}=[\bb{R}^{1|1}_{[0,z^0]},ds+\lambda d\lambda,\gamma_{0,0}, 
  \gamma_{z,\theta},\Psi],\qquad 
I_{z,\theta}=I_{z,\theta,\Omega^{\otimes -n}} 
\eeqq
and the related morphisms $\bar{I}_{z,\theta,\Psi}$, 
$LI_{z,\theta,\Psi}$, $RI_{z,\theta,\Psi}$ as in (\ref{Ibar.LI.RI}). 
Every equivalence class of decorated $B$-bordisms is a disjoint union 
of ones of the form $I_{z,\theta,\Psi}$, $\bar{I}_{z,\theta,\Psi}$, 
$LI_{z,\theta,\Psi}$, $RI_{z,\theta,\Psi}$ and $[Y,\omega,-,-,\Psi]$, 
where $Y_\r=S^1$. 
We define $\mc{SEB}^1_n(f)$ to 
send $I_{z,\theta,\Psi}$ to $I_{f^*z,f^*\theta,\Psi}$ 
(and similarly for $\bar{I}_{z,\theta,\Psi}$, $LI_{z,\theta,\Psi}$, 
$RI_{z,\theta,\Psi}$), 
$[Y,\omega,-,-,\Psi]$ to itself, 
and preserve disjoint unions. 
This finishes the definition of the super category $\mc{SEB}^1_n$.

Let us compute the relations among certain endomorphisms of $Z_{1,0}$ in 
$\mc{SEB}^1_n(B)$. 
For any $c\in C(Z_{1,0})^{\otimes -n}$, we have 
\beq \label{Clifford.epsilon}
(1,c)\circ(\epsilon,1)
=(\epsilon,\epsilon^*c)
=(\epsilon,1)\circ(1,\epsilon^*c). 
\eeq
Notice that, for the $B$-isomorphism $\epsilon$, $\ell(\epsilon)$ is also 
denoted $\epsilon$ (\ref{ell.basic}), which induces the grading 
homomorphism $\epsilon^*$ on $C(Z_{1,0})^{\otimes -n}\cong Cl_{-n}$. 
Since $\lambda\Omega=\Omega\lambda$, where $\lambda\in\mc{O}(Z_{1,0})^\od$ 
generates $C(Z_{1,0})$, we have 
\beq \label{bordism.Clifford}
(1,c)\circ I_{z,\theta}=I_{z,\theta}\circ(1,c).
\eeq
The following analogues of (\ref{epsilon.I.epsilon}) and 
(\ref{sbordism.compose}),
\beq 
&(\epsilon,1)\circ I_{z,\theta,\Psi}\circ(\epsilon,1)=I_{z,-\theta,\Psi},&
  \label{epsilon.I.epsilon.n}  \\
&I_{z_1,\theta_1}\circ I_{z_2,\theta_2}=
  I_{z_1+z_2+\theta_1\theta_2,\theta_1+\theta_2},& 
  \label{Bbordism.compose.n}
\eeq
are proved in the same way. 
For (\ref{Bbordism.compose.n}), we also need the fact that $\Omega$ is 
preserved in tensor products associated with composing $B$-bordisms. 

\defn \label{defn.SEFTn}
A \emph{$1$-SEFT of degree $n$} is a functor 
$E:\mc{SEB}^1_n\rightarrow\mc{S}\t{Hilb}$ of super categories. 
Let $E_B:\mc{SEB}^1_n(B)\rightarrow\mc{S}\t{Hilb}(B)$ be the (ordinary) functor
associated with each object $B$ of $\mc{S}$. 
Since the values of $E_B$ on objects and on decorated $B$-isomorphisms 
does not depend on $B$, 
\footnote{
See definition \ref{defn.SEFT0}.
}
we simply write $E(Z_{p,q})$ and $E(\mu,c)$. 
Now, fix $B$ and let $\Lambda=\mc{O}(B)$. 
The functor $E_B$ is required to satisfy the following assumptions: 
\begin{itemize}
\item[(i)]
$E_{B}$ sends disjoint unions of objects (resp. morphisms) to tensor 
products of graded Hilbert spaces (resp. operators). 
\item[(ii)] 
$E_{B}$ sends objects (resp. $B$-bordisms) with opposite orientations to 
the same graded Hilbert space (resp. operator). 
\item[(iii)]
Given a $B$-bordism $(Y,\omega,\alpha,\beta)$ from $Z_{p,q}$ to $Z_{r,s}$, 
the map 
\beqq
\t{Ad}_{E(Z_{p,q})}:\Lambda\otimes E(Z_{q,p})\otimes E(Z_{r,s})\rightarrow 
\Lambda\otimes B\big(E(Z_{p,q}),E(Z_{r,s})\big), 
\eeqq
sends $E_B([Y,\omega,-,\alpha\sqcup\beta,\Psi])$ to 
$E_B([Y,\omega,\alpha,\beta,\Psi])$, for any $\Psi\in F(Y,\omega)$. 
\item[(iv)]
The operator $\epsilon=E(\epsilon,1)$ on $E(Z_{1,0})$ is the grading 
homomorphism. 
\item[(v)]
$E_B(I_{z,\theta})$ is analytic in $(z,\theta)$. 
\item[(vi)]
$E(\mu,c)$ is linear in $c$; 
$E_B([Y,\omega,\alpha,\beta,\Psi])$ is linear in $\Psi$.
\end{itemize}

{\it Remark.}
The two versions of degree zero $1$-SEFT, 
i.e definition \ref{defn.SEFT0} and the case $n=0$ of definition 
\ref{defn.SEFTn}, are equivalent. 
Indeed, since $C(Z_{p,q})^{\otimes 0}=\bb{R}=F(Y,\omega)^{\otimes 0}$, 
we may regard $\mc{SEB}^1$ as a subcategory of $\mc{SEB}^1_0$ via
$Z_{p,q}\mapsto Z_{p,q}$, $\mu\mapsto(\mu,1)$ and 
$[Y,\omega,\alpha,\beta]\mapsto[Y,\omega,\alpha,\beta,1]$.
Clearly, any functor $\mc{SEB}^1_0\rightarrow\mc{S}\t{Hilb}$ satisfying 
assumption (vi) of \ref{defn.SEFTn} is uniquely determined by 
its restriction to $\mc{SEB}^1$, and any functor 
$\mc{SEB}^1\rightarrow\mc{S}\t{Hilb}$ extends linearly to $\mc{SEB}^1_0$. 

Given a degree $n$ $1$-SEFT $E$, (i) and (ii) imply 
$E(Z_{p,q})=\mc{H}^{\otimes(p+q)}$, where $\mc{H}=E(Z_{1,0})$. 
For decorated $B$-isomorphisms $(\mu,c)$, it suffices by (i) and (iv) to 
consider the cases $\mu=1$ and $c\in C(Z_{1,0})^{\otimes -n}\cong Cl_{-n}$. 
The corresponding operators $E(1,c)$ define a $Cl_{-n}$-action on $\mc{H}$. 
By (\ref{Clifford.epsilon}), it is in fact a graded $Cl_{-n}$-action. 
Now, fix an object $B$ of $\mc{S}$. 
The same arguments following definition \ref{defn.SEFT0} 
show that the values of $E_B$ on all decorated $B$-bordisms 
are determined by the values on $I_{z,\theta,\Psi}$. 
Since $\Psi\in F(\bb{R}^{1|1}_{[0,z^0]})^{\otimes -n}$, which is generated by 
$\Omega^{\otimes -n}$ as a 
$C(Z_{1,0})^{\otimes -n}$-$C(Z_{1,0})^{\otimes -n}$-bimodule, 
the operators $E_B(I_{z,\theta,\Psi})$ are in turn determined by 
$E_B(I_{z,\theta})$ and the $Cl_{-n}$-action on $\mc{H}$. 
(Recall: $I_{z,\theta}:=I_{z,\theta,\Omega^{\otimes -n}}$.)
The morphisms $I_{z,\theta}$ satisfy (\ref{epsilon.I.epsilon.n}) and 
(\ref{Bbordism.compose.n}), hence the same analysis following definition 
\ref{defn.SEFT0} applies, again yielding (\ref{E.I.z.theta}):
\beqq
E_{B}(I_{z,\theta})=X(z)+\theta Y(z)=
\begin{cases}
  e^{-z Q^2}+\theta Qe^{-zQ^2} & \t{on }\mc{H}' \\
  0 & \t{on }\mc{H}'^\perp
\end{cases},
\eeqq
where $Q$ is an odd, self-adjoint operator with compact resolvent on a 
graded subspace $\mc{H}'\subset\mc{H}$. 
Furthermore, by (\ref{bordism.Clifford}), 
all $E_B(I_{z,\theta})$ are $Cl_{-n}$-linear, hence 
$\mc{H}'$ is in fact a graded $Cl_{-n}$-submodule and $Q$ is 
$Cl_{-n}$-linear. 
The operator $Q$ is called \emph{the generator of $E$}. 


\setcounter{equation}{0}
\subsection{Categories $\mc{SEFT}_n$}\label{sec.SEFT}
Throughout this section, we fix an integer $n$ and 
a real Hilbert space $\mc{H}$ with a stable, orthogonal 
\footnote{ \label{stable}
A graded $Cl_{-n}$-action on a vector space is \emph{stable} 
if each isomorphism class of irreducible graded $Cl_{-n}$-module 
appears with infinite muliplicity; 
it is \emph{orthogonal} if $\mc{H}^\ev\perp\mc{H}^\od$ and 
the standard basis elements $e_1,\ldots,e_{|n|}$ of $\bb{R}^{|n|}$, 
regarded as elements of $Cl_{-n}$, act orthogonally on $\mc{H}$. 
}
graded $Cl_{-n}$-action. 

\defn \label{cat.SEFTn}
$\mc{SEFT}_n(\mc{H})$ is a topological category whose objects are 
$1$-SEFTs $E$ such that $E(Z_{1,0})=\mc{H}$. 
We take the discrete topology on the set of objects. 
A morphism from $E$ to $E'$ is a `deformation' from the spectral 
decomposition of the generator $Q_E$ of $E$ to the spectral decomposition 
of the generator $Q_{E'}$ of $E'$.
To explain what this means, let us first describe the spectral decomposition 
of the generator of a $1$-SEFT of degree $n$.

\begin{figure}[t] 
\setlength{\unitlength}{0.037cm}
\begin{picture}(400,40)
  \put(0,20){\line(1,0){400}}
  \put(5,0){$\cdots$}
  \put(50,20){\circle*{3}}
  \put(28,0){\scriptsize $V^E_{-\lambda'}=\epsilon V^E_{\lambda'}$}
  \put(150,20){\circle*{3}}
  \put(132,0){\scriptsize $V^E_{-\lambda}=\epsilon V^E_\lambda$}
  \put(200,20){\circle*{3}}
  \put(198,0){\scriptsize $V^E_0$}
  \put(250,20){\circle*{3}}
  \put(248,0){\scriptsize $V^E_\lambda$}
  \put(350,20){\circle*{3}}
  \put(348,0){\scriptsize $V^E_{\lambda'}$}
  \put(380,0){$\cdots$}
\end{picture}
\caption{The spectral decomposition of $Q_E$ for a $1$-SEFT $E$. 
\label{sp.Q}}
\end{figure}
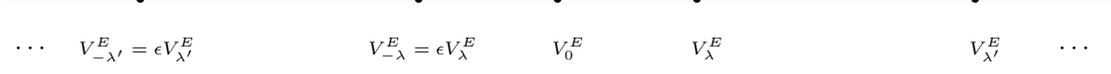

Let $E$ be a $1$-SEFT. 
Recall that its generator $Q_E$ is an (i) odd, (ii) self-adjoint, 
(iii) $Cl_{-n}$-linear operator (defined on a graded $Cl_{-n}$-submodule of 
$\mc{H}$) (iv) with compact resolvent. 
Because of (ii), $Q_E$ has real eigenvalues and its domain decomposes into 
eigenspaces. 
(The orthogonal complement of the domain of $Q_E$, if nontrivial, 
should be regarded as the eigenspace associated with the eigenvalue 
$\infty$.) 
It is useful to have in mind the picture in figure \ref{sp.Q}. 
The line represents the real numbers. 
The points on the real line are the eigenvalues of $Q_E$ and they are 
labeled by the corresponding eigenspaces. 
By (iii), the eigenspaces are $Cl_{-n}$-submodules.
Then, (iv) implies that the eigenvalues have no accumulation point on 
the real line and each eigenspace is finite dimensional. 
Finally, (i) has the following implication. 
If $\epsilon$ denotes the grading homomorphism on $\mc{H}$, (i) means 
$\epsilon Q_E=-Q_E\epsilon$, and hence $\epsilon$ maps the eigenspace 
$V^E_\lambda$ to $V^E_{-\lambda}$. 
This gives rise to a symmetry about the origin in figure \ref{sp.Q}. 
Notice that $V^E_0$ is the only eigenspace that is a 
\emph{graded} $Cl_{-n}$-module. 

We now continue with the definition of $\mc{SEFT}_n(\mc{H})$. 
Given two objects $E,E'$, let $V^E_\lambda$ (resp. $V^{E'}_\lambda$) 
denote the eigenspace of $Q_E$ (resp. $Q_{E'}$) with eigenvalue $\lambda$. 
A deformation $(\alpha,f,A)$ from the spectral decomposition of $Q_E$ to 
that of $Q_{E'}$ consists of the following data:
\begin{itemize}
\item[-]
a map $\alpha:\t{sp}\,Q_E\rightarrow\t{sp}\,Q_{E'}$, which is odd, proper 
and order preserving,
\item[-]
an even linear map $f:\bigoplus_\lambda V^E_\lambda\rightarrow
\bigoplus_{\lambda'}V^{E'}_{\lambda'}$, which embeds each $V^E_\lambda$ 
isometrically into $V^{E'}_{\alpha(\lambda)}$, and 
\item[-]
an (ungraded) $Cl_{-n}$-submodule $A$ of $\mc{H}$, such that $Q_{E'}$ is 
nonnegative on $A$, $A\perp\epsilon A$ and 
$\bigoplus_{\lambda'}V^{E'}_{\lambda'}=f(\bigoplus_\lambda V^E_\lambda)
\oplus A\oplus\epsilon A$. 
\end{itemize}
In particular, fixing $\alpha$, 
\beqq
f\in\prod_{\lambda\in\t{sp}\,Q_E}
  \t{Hom}(V^E_\lambda,V^{E'}_{\alpha(\lambda)}), 
\eeqq
where $\t{Hom}(\cdot,\cdot)$ denotes the space of linear maps between two 
finite dimensional vector spaces, with the usual topology. 
Furthermore, given $\alpha$ and $f$, $A$ is determined by 
\beqq
A\cap V^{E'}_0\in\t{Gr}(V^{E'}_0),
\eeqq
where $\t{Gr}(\cdot)$ denotes the space of linear subspaces in a finite 
dimensional vector space, again with the usual topology. 
The set of morphisms $\{(\alpha,f,A)\}$ from $E$ to $E'$ is topologized 
as a subspace of 
\beqq
\coprod_\alpha\left(
  \prod_{\lambda\in\t{sp}\,Q_E}
  \t{Hom}(V^E_\lambda,V^{E'}_{\alpha(\lambda)})
  \times\t{Gr}(V^{E'}_0)
\right),
\eeqq
which has one component for each $\alpha$. 
Composition of morphisms is defined by 
\beqq
(\alpha',f',A')\circ(\alpha,f,A)=(\alpha'\alpha,f'f,f'(A)\oplus A').
\eeqq 
This finishes the definition of $\mc{SEFT}_n(\mc{H})$.
Sometimes, we simply write $\mc{SEFT}_n$, with the 
underlying Hilbert space $\mc{H}$ unspecified but understood. 

\begin{figure}[p] 
\setlength{\unitlength}{0.037cm}
\begin{picture}(400,480)


  \put(0,440){\line(1,0){400}}
  \put(50,440){\circle*{3}}
  \put(48,450){\scriptsize $V^E_{-\lambda'}$}
  \put(150,440){\circle*{3}}
  \put(148,450){\scriptsize $V^E_{-\lambda}$}
  \put(200,440){\circle*{3}}
  \put(198,450){\scriptsize $V^E_0$}
  \put(250,440){\circle*{3}}
  \put(248,450){\scriptsize $V^E_\lambda$}
  \put(350,440){\circle*{3}}
  \put(348,450){\scriptsize $V^E_{\lambda'}$}

  \put(5,395){$f$}
  \put(50,420){\vector(0,-1){40}}
  \put(55,395){\scriptsize $\cong$}
  \put(150,420){\vector(0,-1){40}}
  \put(155,395){\scriptsize $\cong$}
  \put(200,420){\vector(0,-1){40}}
  \put(205,395){\scriptsize $\cong$}
  \put(250,420){\vector(0,-1){40}}
  \put(255,395){\scriptsize $\cong$}
  \put(350,420){\vector(0,-1){40}}
  \put(355,395){\scriptsize $\cong$}

  \put(0,360){\line(1,0){400}}
  \put(50,360){\circle*{3}}
  \put(48,340){\scriptsize $V^{E'}_{-\lambda'}$}
  \put(150,360){\circle*{3}}
  \put(148,340){\scriptsize $V^{E'}_{-\lambda}$}
  \put(200,360){\circle*{3}}
  \put(198,340){\scriptsize $V^{E'}_0$}
  \put(250,360){\circle*{3}}
  \put(248,340){\scriptsize $V^{E'}_\lambda$}
  \put(350,360){\circle*{3}}
  \put(348,340){\scriptsize $V^{E'}_{\lambda'}$}


  \put(0,280){\line(1,0){400}}
  \put(50,280){\circle*{3}}
  \put(48,290){\scriptsize $V^E_{-\lambda'}$}
  \put(160,280){\circle*{3}}
  \put(158,290){\scriptsize $V^E_{-\lambda}$}
  \put(200,280){\circle*{3}}
  \put(198,290){\scriptsize $V^E_0$}
  \put(240,280){\circle*{3}}
  \put(238,290){\scriptsize $V^E_\lambda$}
  \put(350,280){\circle*{3}}
  \put(348,290){\scriptsize $V^E_{\lambda'}$}

  \put(5,235){$f$}
  \put(60,260){\vector(1,-2){20}}
  \put(60,235){\scriptsize $=$}
  \put(170,260){\vector(1,-2){20}}
  \put(170,235){\scriptsize $\cap$}
  \put(200,260){\vector(0,-1){40}}
  \put(205,235){\scriptsize $\cap$}
  \put(230,260){\vector(-1,-2){20}}
  \put(225,235){\scriptsize $\cap$}
  \put(340,260){\vector(-1,-2){20}}
  \put(335,235){\scriptsize $=$}

  \put(0,200){\line(1,0){400}}
  \put(90,200){\circle*{3}}
  \put(88,180){\scriptsize $V^{E'}_{-\mu}$}
  \put(200,200){\circle*{3}}
  \put(198,180){\scriptsize $V^{E'}_0$}
  \put(310,200){\circle*{3}}
  \put(308,180){\scriptsize $V^{E'}_\mu$}


  \put(0,120){\line(1,0){400}}
  \put(50,120){\circle*{3}}
  \put(48,130){\scriptsize $V^E_{-\lambda}$}
  \put(200,120){\circle*{3}}
  \put(198,130){\scriptsize $V^E_0$}
  \put(350,120){\circle*{3}}
  \put(348,130){\scriptsize $V^E_\lambda$}

  \put(5,75){$f$}
  \put(50,100){\vector(0,-1){40}}
  \put(55,75){\scriptsize $=$}
  \put(200,100){\vector(0,-1){40}}
  \put(205,75){\scriptsize $=$}
  \put(350,100){\vector(0,-1){40}}
  \put(355,75){\scriptsize $=$}

  \put(0,40){\line(1,0){400}}
  \put(50,40){\circle*{3}}
  \put(48,20){\scriptsize $V^{E'}_{-\lambda}$}
  \put(130,40){\circle*{3}}
  \put(128,20){\scriptsize $V^{E'}_{-\lambda'}\subset\epsilon A$}
  \put(200,40){\circle*{3}}
  \put(198,20){\scriptsize $V^{E'}_0$}
  \put(270,40){\circle*{3}}
  \put(268,20){\scriptsize $V^{E'}_{\lambda'}\subset A$}
  \put(350,40){\circle*{3}}
  \put(348,20){\scriptsize $V^{E'}_\lambda$}

\end{picture}
\caption{Three basic types of morphisms in $\mc{SEFT}_n$.
\label{SEFTn.mor}}
\end{figure}
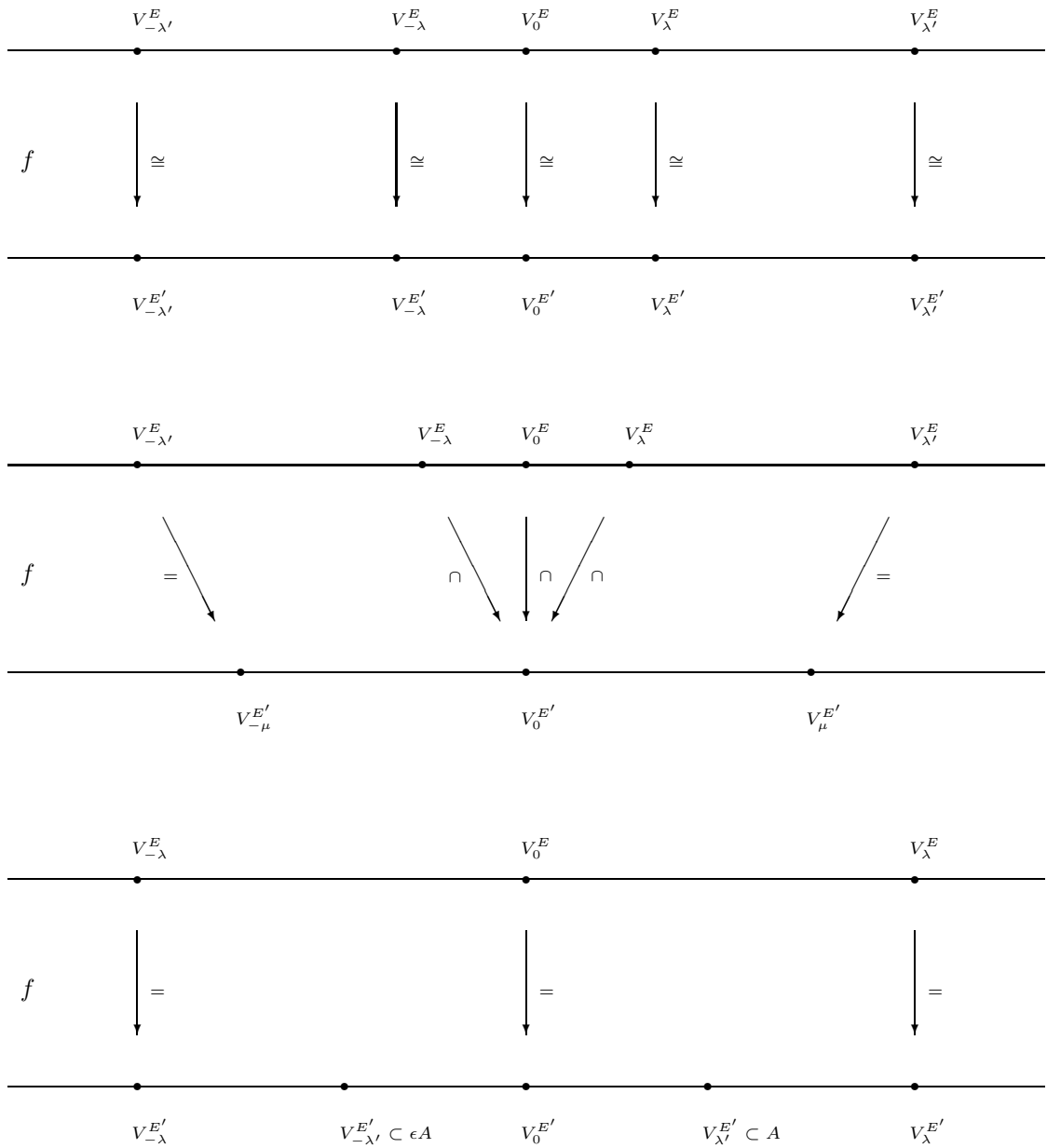

{\it Remarks.}
(a) The assumptions that $\alpha$ is odd and $f$ is even imply that 
$f|_{V^E_{-\lambda}}=\epsilon\circ f|_{V^E_\lambda}\circ\epsilon$. 
In this sense, the symmetry about the origin in the spectral decompositions 
associated with the objects is respected by the morphisms. 
(b) There are three special types of morphisms. (See figure \ref{SEFTn.mor}.)
For the first type, $\t{sp}\,Q_E=\t{sp}\,Q_{E'}$, $\alpha$ is 
the identity map and $A=0$, which means the eigenvalues are the same 
but the eigenspaces change by isomorphisms.
For the second type, $A=0$ and $f$ is identity, which means the eigenspaces 
remain the same but the eigenvalues may have shifted, with the possibility 
of a finite number of eigenvalues merging into one and the corresponding 
eigenspaces combining together. 
For the third type, $\t{sp}\,Q_E\subset\t{sp}\,Q_{E'}$, 
$\alpha$ is the inclusion and $f$ is also the inclusion, which means 
the original eigenvalues and eigenspaces are unchanged but there are 
new ones that `emerge from $+\infty$ (those in $A$) and 
$-\infty$ (those in $\epsilon A$).'
All morphisms are generated by these three special types. 
Indeed, any $(\alpha,f,A)$ can be factored as 
\beqq
(\t{inclusion},\t{inclusion},A)\circ(\t{identity},f,0)\circ
(\alpha,\t{identity},0).
\eeqq

The result below concerns the homotopy types of $|\mc{SEFT}_n(\mc{H})|$. 
To state the result, we need to first define a space of operators 
$Fred_n(\mc{H})$. 
Let $F_n(\mc{H})$ denote the space (with the norm topology) of Fredholm 
operators on $\mc{H}$ which are odd, self-adjoint and $Cl_{-n}$-linear.
If $n\not\equiv 1\pmod{4}$, we define $Fred_n(\mc{H})=F_n(\mc{H})$. 
If $n\equiv 1\pmod{4}$, consider an orthonormal basis $e_1,\ldots,e_{|n|}$ of 
$\bb{R}^{|n|}$, regarded as elements of $Cl_{-n}$ and in turn as orthogonal 
operators on $\mc{H}$. 
For each $F\in F_n(\mc{H})$, the operator $e_1\cdots e_{|n|}F$ is even and 
self-adjoint. 
$Fred_n(\mc{H})$ is defined to be the subspace of $F_n(\mc{H})$ consisting of 
those $F$ such that 
\beq \label{ess.pos.neg}
(e_1\cdots e_{|n|}F)|_{\mc{H}^\ev} \t{ is neither essentially positive nor 
essentially negative.}
\eeq 
\footnote{
A self-adjoint operator is essentially positive (resp. negative) if it is 
positive (resp. negative) on an invariant subspace of a finite codimension.
}
Sometimes we write $Fred_n$ instead of $Fred_n(\mc{H})$ when 
the underlying Hilbert space is understood.

{\thm \label{SEFT.Fred}
For each $n\in\bb{Z}$, there is a homotopy equivalence 
\beqq
|\mc{SEFT}_n(\mc{H})|\simeq Fred_n(\mc{H}).
\eeqq
}

{\it Remark.}
It follows from \cite{AS.Fred} (see also \cite{LM}, III, \S 10) that, 
for each $n\leq 0$, $Fred_n$ represents the functor $KO^n$. 
We claim that this is also true for $n>0$. 
Suppose $n=4k-\ell$, where $k>0$ and $\ell\geq 0$ are integers. 
Let $\mc{H}$ continue to be a Hilbert space with a stable graded 
$Cl_{-n}$-action and $V$ be an irreducible graded $Cl_{\ell,\ell}$-module. 
Since $Cl_{0,n}\otimes Cl_{\ell,\ell}\cong Cl_{\ell,4k}\cong Cl_{4k+\ell}$ 
as a graded tensor product, $\mc{H}\otimes V$ admits a stable graded 
$Cl_{4k+\ell}$-action with the even and odd parts being 
$(\mc{H}^\ev\otimes V^\ev)\oplus(\mc{H}^\od\otimes V^\od)$ and 
$(\mc{H}^\ev\otimes V^\od)\oplus(\mc{H}^\od\otimes V^\ev)$. 
Let us express every operator on a graded vector space as a $2\times 2$ block 
with respect to the parity decomposition of the vector space. 
The map 
\beq 
\left(\begin{array}{cc} 0 & F^1 \\ F^0 & 0 \end{array}\right) \mapsto 
\left(\begin{array}{cc} 
  0 & -F^0\otimes 1_{V^\od}+F^1\otimes 1_{V^\ev} \\ 
  F^0\otimes 1_{V^\ev}-F^1\otimes 1_{V^\od} & 0 
\end{array}\right)
\label{map.btx.Fred} 
\eeq
defines a homeomorphism from $Fred_n(\mc{H})$ to 
$Fred_{-4k-\ell}(\mc{H}\otimes V)$. 
\footnote{
Had we chosen to use skew-adjoint, Clifford anitlinear operators to define 
$Fred_n$, the map in (\ref{map.btx.Fred}) would simply be 
$F\mapsto F\otimes\t{id}_V$.
(See \cite{AS.Fred}, Theorem 5.1.)
}
Since $n\equiv -4k-\ell\pmod{8}$, this shows that $Fred_n(\mc{H})$ represents 
$KO^{-4k-\ell}=KO^n$.

By the above remark, theorem \ref{SEFT.Fred} implies that the space 
$|\mc{SEFT}_n|$ represents $KO^n$. 
The proof of the theorem, which is given in section \ref{sec.proof}, 
consists of two steps
\beqq
|\mc{SEFT}_n(\mc{H})|\stackrel{\sim}{\rightarrow}|\mc{V}_n(\mc{H})|\simeq 
Fred_n(\mc{H}),
\eeqq 
where $\mc{V}_n(\mc{H})$ is a category we define below. 


\setcounter{equation}{0}
\subsection{Categories $\mc{V}_n$} \label{sec.Vn}

\defn
Let $n$ be an integer and $\mc{H}$ a real Hilbert space with a stable, 
orthogonal (see footnote \ref{stable}) graded $Cl_{-n}$-action. 
$\mc{V}_n(\mc{H})$ is a topological category whose objects are 
finite dimensional graded $Cl_{-n}$-submodules of $\mc{H}$, and the set 
of objects is discrete. 
A morphism from an object $V$ to another object $V'$ is a pair 
$(f,A)$, where $f:V\hookrightarrow V'$ is an injective map of graded 
$Cl_{-n}$-modules, and $A\subset V'$ is an (ungraded) 
$Cl_{-n}$-submodule such that $A\perp\epsilon A$, $A\perp f(V)$ and 
$V'=f(V)\oplus A\oplus\epsilon A$. 
The set of morphisms from $V$ to $V'$ is topologized as a subspace of 
$\t{Hom}(V,V')\times \t{Gr}(V')$. 
Composition is defined by 
\beqq
(f',A')\circ(f,A)=(f'f,f'(A)\oplus A').
\eeqq

For the proof of theorem \ref{SEFT.Fred}, the reader can go directly to 
section \ref{sec.proof}. 
The rest of this section is a study of the categories $\mc{V}_n(\mc{H})$. 

\vspace{0.1in}
\noindent
{\bf Comparing with Quillen's categories $\wh{Vect}$ and $QVect$.}
Each $\mc{V}_n$ is designed to admit a functor from $\mc{SEFT}_n$, 
as we will see in the next section. 
Curiously, $\mc{V}_0$ and $\mc{V}_1$ coincide with two 
categories studied by Quillen \cite{Quillen.algK} and Segal 
\cite{Segal.Khlgy}.

Let us first recall $\wh{Vect}$, the category of `virtual vector spaces.' 
To make explicit the dependence on two underlying infinite dimensional 
Hilbert spaces $\mc{H}^0$ and $\mc{H}^1$, we use the notation 
$\wh{Vect}(\mc{H}^0,\mc{H}^1)$. 
An object of $\wh{Vect}(\mc{H}^0,\mc{H}^1)$ is a pair $(V^0,V^1)$, where 
$V^0\subset\mc{H}^0$ and $V^1\subset\mc{H}^1$ are finite dimensional 
subspaces. 
A morphism from an object $(V^0,V^1)$ to another object $(V'^0,V'^1)$ is a 
triple $(f^0,f^1;\phi)$, where $f^i:V^i\rightarrow V'^i$ are isometric 
injections and $\phi:U^0\stackrel{\sim}{\rightarrow}U^1$ is an isometric 
isomorphism, $U^i$ being the orthogonal complement of $f^i(V^i)$ in $V'^i$, 
$i=0,1$. 
The composition of $(f^0,f^1;\phi)$ with a morphism $(f'^0,f'^1;\phi')$ 
from $(V'^0,V'^1)$ to $(V''^0,V''^1)$ is defined by 
\beqq
(f'^0,f'^1;\phi')\circ(f^0,f^1;\phi)=
(f'^0 f^0,f'^1 f^1; f'^1\phi(f'^0)^{-1}\oplus\phi').
\eeqq
More precisely, the map $({f'^0})^{-1}$ appearing here is defined on 
$f'^0(U^0)$, so that $f'^1\phi{f'^0}^{-1}\oplus\phi'$ is indeed defined on 
the orthogonal complement of $(f'^0 f^0)(V^0)$ in $V''^0$. 

Let $Vect$ denote the category of finite dimensional vector spaces (say in 
$\mc{H}^0$) and isomorphisms. 
Taking direct sums defines a symmetric monoidal structure on $Vect$. 
It is known \cite{Segal.Khlgy} that
\beqq
|\wh{Vect}(\mc{H}^0,\mc{H}^1)|\simeq\Omega B|Vect|
\simeq \Omega B\left(\coprod_{n\geq 0} BGL_n\right),
\eeqq 
i.e. $\wh{Vect}(\mc{H}^0,\mc{H}^1)$ is a `category-level group completion' 
of $Vect$. 

Suppose $\mc{H}$ is a graded Hilbert space such that 
$\mc{H}^\ev=\mc{H}^0$, $\mc{H}^\od=\mc{H}^1$ 
are both infinite dimensional.
We have the following

{\prop \label{V0.virVect}
The categories $\mc{V}_0(\mc{H})$ and $\wh{Vect}(\mc{H}^0,\mc{H}^1)$ are 
isomorphic.}

\proof 
We define a functor
\beqq
F:\mc{V}_0(\mc{H})\rightarrow\wh{Vect}(\mc{H}^0,\mc{H}^1).
\eeqq
For an object $V$ of $\mc{V}_0(\mc{H})$, i.e. a finite dimensional graded 
subspace of $\mc{H}$, we have 
\beqq
F(V)=(V\cap\mc{H}^0,V\cap\mc{H}^1).
\eeqq
Suppose $(f,A)$ is a morphism from $V$ to $V'$ in $\mc{V}_0(\mc{H})$. 
By definition, the subspaces $A$, $\epsilon A$, and $f(V)$ of $V'$ are 
pairwise orthogonal such that $V'=f(V)\oplus A\oplus\epsilon A$. 
Define a map $\phi_A:(A\oplus\epsilon A)\cap\mc{H}^0\rightarrow
(A\oplus\epsilon A)\cap\mc{H}^1$ as follows. 
Let $u=a+\epsilon a'$ be an element of $(A\oplus\epsilon A)\cap \mc{H}^0$, 
where $a,a'\in A$.  
Since $\epsilon u=u$, we in fact have $a=a'$, hence 
$(A\oplus\epsilon A)\cap\mc{H}^0$ consists precisely of elements of the form 
$a+\epsilon a$, $a\in A$.
Similarly, $(A\oplus\epsilon A)\cap\mc{H}^1$ consists precisely of elements 
of the form $a-\epsilon a$, $a\in A$.
The map $\phi_A$ is defined by $\phi_A(a+\epsilon a)=a-\epsilon a$ for 
any $a\in A$. 
Now, the definition of the functor $F$ on morphisms is given by
\beqq
F(f,A)=(f|_{V\cap\mc{H}^0},f|_{V\cap\mc{H}^1};\phi_A).
\eeqq
Indeed, $f|_{V\cap\mc{H}^i}$ embeds $V\cap\mc{H}^i$ isometrically 
into $V'\cap\mc{H}^i$, and $(A\oplus\epsilon A)\cap\mc{H}^i$ is the 
orthogonal complement of $f(V\cap\mc{H}^i)$ in $V'\cap\mc{H}^i$, where 
$i=0,1$. 
Also, since $A\perp\epsilon A$, we have $||a+\epsilon a||=||a-\epsilon a||$, 
so that $\phi_A$ is indeed an isometry.
It is straightforward to verify that $F$ respects compositions. 

To show that $F$ is an isomorphism of categories, we define a functor
\beqq
G:\wh{Vect}(\mc{H}^0,\mc{H}^1)\rightarrow\mc{V}_0(\mc{H})
\eeqq
and show that it is the inverse of $F$. 
For an object $(V^0,V^1)$ of $\wh{Vect}(\mc{H}^0,\mc{H}^1)$, we have 
\beqq
G(V^0,V^1)=V^0\oplus V^1,
\eeqq
regarded as a graded subspace of $\mc{H}$, hence an object of 
$\mc{V}_0(\mc{H})$. 
Suppose $(f^0,f^1;\phi)$ is a morphism from $(V^0,V^1)$ to $({V'}^0,{V'}^1)$ 
in $\wh{Vect}(\mc{H}^0,\mc{H}^1)$. 
In particular, $\phi:U^0\rightarrow U^1$ is an isometric isomorphism, where 
$U^i$ is the orthogonal complement of $f^i(V^i)$ in ${V'}^i$, $i=0,1$. 
Let $A_\phi=\{u+\phi(u)|u\in U^0\}$, and we define
\beqq
G(f^0,f^1;\phi)=(f^0\oplus f^1,A_\phi).
\eeqq
Notice that $\epsilon A_\phi=\{u-\phi(u)|u\in U^0\}$ and  
$U^0\oplus U^1=A_\phi \oplus\epsilon A_\phi$. 
Also, the facts that $\phi$ is an isometry and $U^0\perp U^1$ imply 
$A_\phi\perp\epsilon A_\phi$.
These show that $G(f^0,f^1;\phi)$ just defined is indeed a morphism from 
$G(V^0,V^1)$ to $G({V'}^0,{V'}^1)$ in $\mc{V}_0(\mc{H})$. 
Again, it is straightforward to verify that $G$ respects compositions. 

Let us check that $GF$ and $FG$ are identity functors. 
It is clear on the objects. 
Suppose $(f,A)$ is a morphism from an object $V$ to another object $V'$ in 
$\mc{V}_0(\mc{H})$. 
We have $GF(f,A)=(f,A_{\phi_A})$, where $A_{\phi_A}$ consists precisely of 
the elements of the form $(a+\epsilon a)+\phi_A(a+\epsilon a)$, $a\in A$. 
Since $\phi_A(a+\epsilon a)=a-\epsilon a$, the elements of $A_{\phi_A}$ are 
in fact of the form $2a$, where $a\in A$.
Hence, $A_{\phi_A}=A$ and $GF$ is indeed the identity on $\mc{V}_0(\mc{H})$. 
On the other hand, consider a morphism $(f^0,f^1;\phi)$ from an object 
$(V^0,V^1)$ to another object $(V'^0,V'^1)$ in $\wh{Vect}(\mc{H}^0,\mc{H}^1)$. 
Let $U^i$ again denote the orthogonal complement of $f^i(V^i)$ in $V'^i$, 
$i=0,1$. 
We have $FG(f^0,f^1;\phi)=(f^0,f^1;\phi_{A_\phi})$.
Since $A_\phi=\{u+\phi(u)|u\in U^0\}$, the map 
$\phi_{A_\phi}:U^0\rightarrow U^1$ is given by 
\beqq
\phi_{A_\phi}(u)
&=& \phi_{A_\phi}\left[\left(\frac{u}{2}+\frac{\phi(u)}{2}\right)+
\epsilon\left(\frac{u}{2}+\frac{\phi(u)}{2}\right)\right] \\
&=& \left(\frac{u}{2}+\frac{\phi(u)}{2}\right)-
\epsilon\left(\frac{u}{2}+\frac{\phi(u)}{2}\right) \\
&=& \phi(u).
\eeqq
Thus, $\phi_{A_\phi}=\phi$ and $FG$ is indeed the identity on 
$\wh{Vect}(\mc{H}^0,\mc{H}_1)$. 
\qedsymbol

Now we recall $QVect$, the $Q$-construction of $Vect$. 
The make explicit its dependence on an underlying infinite dimensional 
Hilbert space $\mc{H}'$, we use the notation $QVect(\mc{H}')$. 
An object of $QVect(\mc{H}')$ is a finite dimensional subspace of $\mc{H}'$. 
A morphism from an object $U$ to another object $U'$ is a triple 
$(g;W_1,W_2)$, where $g:U\rightarrow U'$ is an isometric embedding, and 
$W_1,W_2\subset U'$ are subspaces such that $U'=g(U)\oplus W_1\oplus W_2$, 
and $g(U)$, $W_1$ and $W_2$ are pairwise orthogonal. 
The composition of $(g;W_1,W_2)$ with a morphism $(g';W'_1,W'_2)$ from 
$U'$ to $U''$ is given by
\beqq
(g';W'_1,W'_2)\circ(g;W_1,W_2)=(g'g;g'(W_1)\oplus W'_1,g'(W_2)\oplus W'_2).
\eeqq
It is known \cite{Segal.Khlgy} that 
\beqq
\Omega|QVect(\mc{H}')|\simeq\Omega B|Vect| 
\simeq \Omega B\left(\coprod_{n\geq 0} BGL_n\right).
\eeqq

Suppose $\mc{H}$ is a Hilbert space with a stable, orthogonal graded 
$Cl_{-1}$-action. 
Let $p\in Cl_{-1}=\bb{R}\oplus\bb{R}$ denote the projection onto the 
first factor.
Notice that, as operators on a graded $Cl_{-1}$-module, we have 
$\epsilon p=(1-p)\epsilon$. 
We have the following

{\prop \label{V1.QVect}
The categories $\mc{V}_1(\mc{H})$ and $QVect(p\mc{H})$ are isomorphic.
}

\proof
We first define a functor 
\beqq
F:\mc{V}_1(\mc{H})\rightarrow QVect(p\mc{H}).
\eeqq
For an object $V$ of $\mc{V}_1(\mc{H})$, i.e. a finite dimensional 
$Cl_{-1}$-submodule of $\mc{H}$, we have
\beqq
F(V)=pV.
\eeqq
Given a morphism $(f,A)$ from an object $V$ to another object $V'$, we
define 
\beqq
F(f,A)=(f|_{pV};pA,p\epsilon A).
\eeqq
Here, since $f$ embeds $V$ isometrically into $V'$ and 
$f$ is $Cl_{-1}$-linear, 
$f|_{pV}$ embeds $pV$ isometrically into $pV'$. 
Also, we have $pV'=f(pV)\oplus pA\oplus p\epsilon A$, hence 
$F(f,A)$ defined above is indeed a morphism from $F(V)$ 
to $F(V')$ in $QVect(p\mc{H})$. 

Then we define a functor 
\beqq
G:QVect(p\mc{H})\rightarrow\mc{V}_1(\mc{H})
\eeqq
that is the inverse of $F$. 
Given an object $U$ of $QVect(p\mc{H})$, let
\beqq
G(U)=U\oplus\epsilon U.
\eeqq
Notice that $\epsilon U=\epsilon pU=(1-p)\epsilon U\subset (1-p)\mc{H}$ is 
orthogonal to $U\subset p\mc{H}$, and $U\oplus\epsilon U$ is indeed a 
graded $Cl_{-1}$-submodule. 
Given a morphism $(g;W_1,W_2)$ from $U$ to $U'$ in $QVect(p\mc{H})$, define
\beqq
G(g;W_1,W_2)=(g\oplus\epsilon g\epsilon,W_1\oplus\epsilon W_2).
\eeqq
To see that the right hand side is indeed a morphism from $U\oplus\epsilon U$ 
to $U'\oplus\epsilon U'$ in $\mc{V}_1(\mc{H})$, observe that
\beqq
U'\oplus\epsilon U'
&=& \big(g(U)\oplus W_1\oplus W_2\big)\oplus
\epsilon\big(g(U)\oplus W_1\oplus W_2\big) \\
&=& (g\oplus\epsilon g\epsilon)(U\oplus\epsilon U)\oplus
(W_1\oplus\epsilon W_2)\oplus\epsilon(W_1\oplus\epsilon W_2).
\eeqq
It is clear that $g\oplus\epsilon g\epsilon$ commutes with $p$ and 
$\epsilon$, hence is graded and $Cl_{-1}$-linear; 
also, $p(W_1\oplus\epsilon W_2)\subset(W_1\oplus\epsilon W_2)$, so that 
$W_1\oplus\epsilon W_2$ is indeed a $Cl_{-1}$-submodule. 

Finally, we check that $GF$ and $FG$ are identity functors. 
For objects, we observe 
\beqq
GF(V)&=&pV\oplus\epsilon pV=pV\oplus (1-p)\epsilon V=pV\oplus(1-p)V=V \\
FG(U)&=&p(U\oplus\epsilon U)=pU\oplus\epsilon(1-p)U=U\oplus 0=U.
\eeqq
For a morphism $(f,A)$ from $V$ to $V'$ in $\mc{V}_1(\mc{H})$, 
we have by definition 
$GF(f,A)=\big(f|_{pV}\oplus\epsilon(f|_{pV})\epsilon,
pA\oplus\epsilon p\epsilon A\big)$.
It follows from $\epsilon p=(1-p)\epsilon$ that 
this is equal to $(f,A)$. 
Therefore, $GF$ is the identity on $\mc{V}_1(\mc{H})$. 
For a morphism $(g;W_1,W_2)$ from an object $U$ to another object $U'$ in 
$QVect(p\mc{H})$, we have  
\beqq
FG(g;W_1,W_2)=\big((g\oplus\epsilon g\epsilon)|_{p(U\oplus\epsilon U)};
p(W_1\oplus\epsilon W_2),p\epsilon(W_1\oplus\epsilon W_2)\big).
\eeqq
Since $p(U\oplus\epsilon U)=U$, $p(W_1\oplus\epsilon W_2)=W_1$ and 
$p(\epsilon W_1\oplus W_2)=W_2$, the right hand side above equals 
$(g;W_1,W_2)$, proving that $FG$ is the identity on $QVect(p\mc{H})$.
\qedsymbol

\vspace{0.1in}
\noindent
{\bf Components of $\mc{V}_n$.}
For any $n\in\bb{Z}$, let $\M_n$ denote the set of isomorphism classes 
of finite dimensional graded $Cl_n$-modules.
Given such a module $V$, we denote its isomorphism class by $[V]$.
Taking direct sums of graded $Cl_n$-modules induces the structure of an 
abelian monoid on $\M_n$. 
Recall that (for $n\neq 0$) $Cl_n$ is generated by a set of elements 
$e_1,\ldots,e_{|n|}$ satisfying the relations $e_i^2=-\t{sgn}(n)$ for any 
$i$, and $e_i e_j=-e_j e_i$ for any $i\neq j$. 

For each $n$, we define a monoid homomorphism 
\beqq
i_n:\M_n\rightarrow\M_{n-1}.
\eeqq
For $n\geq 1$, $i_n$ is induced by restricting a graded $Cl_n$-action 
to the action generated by $e_1,\ldots,e_{n-1}$. 
For $n\leq 0$, let $k=-n$ and $\eta\in\M_{-k}$. 
Suppose $V^e$ and $V^o$ are two graded $Cl_{-k}$-modules representing $\eta$, 
and $\gamma:V^e\rightarrow V^o$ is a module isomorphism. 
Denote by $\epsilon$ the grading homomorphisms on both $V^e$ and $V^o$. 
Let $V=V^e\oplus V^o$ be the graded vector space whose even and odd parts 
are $V^e$ and $V^o$ respectively. 
The following operators
\beq \label{in.neg}
e'_i=\begin{cases}
  \epsilon e_i\gamma & \t{on }V^e \\
  -\epsilon e_i\gamma^{-1} & \t{on }V^o
\end{cases}\quad 1\leq i\leq k, 
\qquad
e'_{k+1}=\begin{cases}
  \gamma & \t{on } V^e \\
  \gamma^{-1} & \t{on } V^o
\end{cases}
\eeq
reverse the parity on $V$, and satisfy ${e'_i}^2=1$ for 
$1\leq i\leq k+1$, $e'_i e'_j=-e'_j e'_i$ for $1\leq i<j\leq k+1$.
Hence, they generate a graded $Cl_{-k-1}$-action on $V$. 
We define $i_{-k}\eta$ to be $[V]\in\M_{-k-1}$.

{\prop \label{prop.pi0.Vn}
For each $n\in\bb{Z}$, there is a bijective correspondance
\beq \label{pi0.Vn}
\pi_0\mc{V}_{-n}\cong\M_n/i_{n+1}\M_{n+1}.
\eeq
Furthermore, the right hand side above is in fact an abelian group.
}

{\it Remarks.}
(a) If $M$ is an abelian monoid and $N\subset M$ is a submonoid, the set 
$M/N$ is defined to be $M/\sim$, where $m\sim m'$ if and only if 
$m+n_1=m'+n_2$ for some $n_1,n_2\in N$.
(b) Together with subsequent results on the categories $\mc{V}_{-n}$, this
proposition recovers the well known interpretation of $KO^{-n}(\t{pt})$ in 
terms of Clifford modules \cite{ABS}, at least for $n\geq 0$.
The construction of the right hand side of (\ref{pi0.Vn}) for $n<0$ is new to 
the author. 

\proof[Proof of Proposition \ref{prop.pi0.Vn}]
Recall that the objects of $\mc{V}_{-n}$ are finite dimensional graded 
$Cl_n$-submodules of a stable graded $Cl_n$-module $\mc{H}$. 
Since submodules of every isomorphism class appear in $\mc{H}$ and two 
isomorphic ones clearly lie in the same component of $\mc{V}_{-n}$, there 
is a surjective map
\beqq
\phi:\M_n\rightarrow\pi_0\mc{V}_{-n}.
\eeqq

Now, we show that $\phi$ factors through the projection 
$\M_n\rightarrow\M_n/i_{n+1}\M_{n+1}$. 
It suffices to show that, for any $\eta\in\M_n$ and $\xi\in i_{n+1}\M_{n+1}$, 
$\phi(\eta)=\phi(\eta+\xi)$. 
Suppose $V,W\subset\mc{H}$ are two orthogonal graded 
$Cl_n$-submodules such that $[V]=\eta$, $[W]=\xi$. 
We show that there is a morphism from $V$ to $V\oplus W$ in $\mc{V}_{-n}$. 

First consider the case $n\geq 0$. 
The assumption that $[W]\in i_{n+1}\M_{n+1}$ implies there is an operator 
$e_{n+1}$ on $W$ extending the graded $Cl_n$-action on $W$ to a graded 
$Cl_{n+1}$-action. 
Let $A=\{w+e_{n+1}w|w\in W^\ev\}\subset W$. 
Since for any $w\in W^\ev$ and $1\leq i\leq n$, we have 
$e_i(w+e_{n+1}w)=(e_i e_{n+1}w)+e_{n+1}(e_i e_{n+1}w)$ and 
$e_i e_{n+1}w\in W^\ev$, $A$ is a $Cl_n$-submodule.
Also, $\epsilon A=\{w-e_{n+1}w|w\in W^\ev\}$ is orthogonal 
to $A$, and the fact that $e_{n+1}$ maps $W^\ev$ isomorphically to 
$W^\od$ implies $W=A\oplus\epsilon A$. 
Therefore, $(\t{inc},A)$ is a morphism from $V$ to $V\oplus W$ in the 
category $\mc{V}_{-n}$, where $\t{inc}$ is the inclusion of 
$V$ into $V\oplus W$. 

Then we consider the case $n\leq -1$.
Let $k=-n$. 
By definition of $i_{-k+1}$, there are isomorphic graded 
$Cl_{-(k-1)}$-actions on $W^\ev$ and $W^\od$, 
whose relation with the graded $Cl_{-k}$-action on $W$ is described in 
(\ref{in.neg}). 
Let us use the same notations as in (\ref{in.neg}). 
Let $A=\{w+\epsilon\gamma w|w\in W^\ev\}\subset W$. 
For $1\leq i\leq k-1$ and $w\in W^\ev$, we have
\beqq
e'_i(w+\epsilon\gamma w)
=\epsilon e_i\gamma w-\epsilon e_i\gamma^{-1}\epsilon\gamma w
=(e_i w)+\epsilon\gamma (e_i w), 
\eeqq
and 
\beqq
e'_k(w+\epsilon\gamma w)
=\gamma w+\gamma^{-1}\epsilon\gamma w
=(\epsilon w)+\epsilon\gamma(\epsilon w), 
\eeqq
showing that $A$ is a $Cl_{-k}$ submodule of $W$. 
Since $\epsilon'A=\{w-\epsilon\gamma w|w\in W^0\}$ 
is orthogonal to $A$, 
$(\t{inc},A)$ is a morphism from $V$ to $V\oplus W$ in $\mc{V}_k$, 
where $\t{inc}$ is again the inclusion of $V$ into $V\oplus W$. 
This finishes the proof that $\phi$ induces a map
\beqq
\bar{\phi}:\M_n/i_{n+1}\M_{n+1}\rightarrow\pi_0\mc{V}_{-n}.
\eeqq

It remains to show that $\bar{\phi}$ is one to one.
Consider a morphism $(f,A)$ from $V$ to $V'$ in $\mc{V}_{-n}$. 
Since $V'=f(V)\oplus A\oplus\epsilon A$ and $f$ is a graded, 
$Cl_n$-linear embedding, we have $[V']=[V]+[A\oplus\epsilon A]\in\M_n$. 
Therefore, it suffices to show that 
\beq \label{A+epsilonA}
[A\oplus\epsilon A]\in i_{n+1}\M_{n+1}. 
\eeq 

If $n\geq 0$, the operator $e_{n+1}$ on $A\oplus\epsilon A$ defined by 
$e_{n+1}|_A=-\epsilon$, $e_{n+1}|_{\epsilon A}=\epsilon$ 
extends the graded $Cl_n$-action on $A\oplus\epsilon A$ to a graded 
$Cl_{n+1}$-action, proving (\ref{A+epsilonA}) in this case. 
For $n\leq -1$, let $k=-n$ and $W=A\oplus\epsilon A$. 
As in the proof of proposition \ref{V0.virVect}, 
$W^\ev=\{a+\epsilon a|a\in A\}$ and $W^\od=\{a-\epsilon a|a\in A\}$. 
Let $\phi_A:W^\ev\rightarrow W^\od$ be the map defined by 
$\phi_A(a+\epsilon a)=a-\epsilon a$. 
We denote by $e'_1,\ldots,e'_k$ the generators (as above) of the 
$Cl_{-k}$-action on $W$. 
Notice that for $i=1,\ldots,k$, $e'_i\phi_A=\phi_A^{-1}e'_i$ as operators 
on $W^\ev$. 
It then follows that the following operators on $W^\ev$ 
\beqq
\epsilon=\phi_A^{-1}e'_k,\qquad
e_i=\phi_A^{-1}e'_i,\quad i=1,\ldots,k-1
\eeqq
define a graded $Cl_{-(k-1)}$-action. 
Similarly, the following operators on $W^\od$
\beqq
\epsilon=\phi_A e'_k,\qquad
e_i=-\phi_A e'_i,\quad i=1,\ldots,k-1
\eeqq
also define a graded $Cl_{-(k-1)}$-action.
Furthermore, $\gamma=e'_k|_{W^\ev}:W^\ev\rightarrow W^\od$ 
commutes with these actions. 
According to (\ref{in.neg}), one may use these data to define a graded 
$Cl_{-k}$-action on $W^\ev\oplus W^\od=W$, 
which is easily checked to coincide with the graded $Cl_{-k}$-action 
we started with. 
This proves $[W]\in i_{-k+1}\M_{-k+1}$. 
We have now established the bijection
$\pi_0\mc{V}_{-n}\cong\M_n/i_{n+1}\M_{n+1}$. 

Finally, we show that the monoid structure on $\M_n/i_{n+1}\M_{n+1}$ induced 
from taking direct sums admits an inverse. 
Given $\eta\in\M_n$, let $V$ be a graded $Cl_n$-module representing $\eta$. 
Reversing the parity in $V$ yields another graded $Cl_n$-module $\bar{V}$. 
If $\epsilon$ denotes the grading homomorphism of $V$, the grading 
homomorphism of $V\oplus\bar{V}$ is given by $\epsilon\oplus -\epsilon$. 
Let $\Delta=\{(v,v)|v\in V\}\subset V\oplus\bar{V}$. 
Since $V\oplus\bar{V}=\Delta\oplus(\epsilon\oplus -\epsilon)\Delta$ 
and $\Delta$ is a $Cl_n$-submodule, 
by (\ref{A+epsilonA}), $[V]+[\bar{V}]\in i_{n+1}\M_{n+1}$. 
Therefore, the images of $[V]=\eta$ and $[\bar{V}]$ in 
$\M_n/i_{n+1}\M_{n+1}$ are inverses of each other. 
\qedsymbol


\setcounter{equation}{0}
\subsection{Proof of $|\mc{SEFT}_n|\simeq Fred_n$} \label{sec.proof}
Fix an integer $n$ and a Hilbert space $\mc{H}$ with a stable, orthogonal 
(see footnote \ref{stable}) graded $Cl_{-n}$-action.
Consider the categories $\mc{SEFT}_n=\mc{SEFT}_n(\mc{H})$, 
$\mc{V}_n=\mc{V}_n(\mc{H})$ and the space $Fred_n=Fred_n(\mc{H})$ 
defined in the previous sections. 

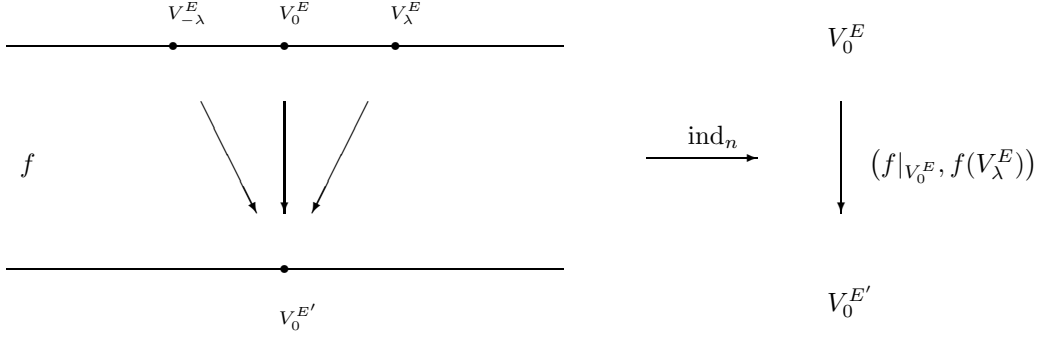
\begin{figure}[t]
\setlength{\unitlength}{0.037cm}
\begin{picture}(400,150)
  \put(0,120){\line(1,0){200}}
  \put(60,120){\circle*{3}}
  \put(58,130){\scriptsize $V^E_{-\lambda}$}
  \put(100,120){\circle*{3}}
  \put(98,130){\scriptsize $V^E_0$}
  \put(140,120){\circle*{3}}
  \put(138,130){\scriptsize $V^E_\lambda$}

  \put(5,75){$f$}
  \put(70,100){\vector(1,-2){20}}
  \put(100,100){\vector(0,-1){40}}
  \put(130,100){\vector(-1,-2){20}}

  \put(0,40){\line(1,0){200}}
  \put(100,40){\circle*{3}}
  \put(98,20){\scriptsize $V^{E'}_0$}

  \put(230,80){\vector(1,0){40}}
  \put(245,85){$\t{ind}_n$}

  \put(295,120){$V^E_0$}
  \put(295,25){$V^{E'}_0$}
  \put(300,100){\vector(0,-1){40}}
  \put(310,75){$\big(f|_{V^E_0},f(V^E_\lambda)\big)$}
\end{picture}
\caption{The functor $\t{ind}_n$ on a particular morphism in $\mc{SEFT}_n$.
\label{ind}}
\end{figure}

\vspace{0.1in}
\noindent
{\bf Comparing $\mc{SEFT}_n$ and $\mc{V}_n$.}
We define a functor
\beqq
\t{ind}_n:\mc{SEFT}_n\rightarrow\mc{V}_n.
\eeqq
Let $E$ be an object of $\mc{SEFT}_n$, 
and $V^E_\lambda=\ker(Q_E-\lambda)$ the eigenspaces of the generator 
$Q_E$ of $E$ (see section \ref{sec.SEFT}). 
On objects, $\t{ind}_n$ is given by 
\beqq
\t{ind}_n E=V^E_0.
\eeqq
For a morphism $(\alpha,f,A)$ from an object $E$ to another object $E'$ in 
$\mc{SEFT}_n$, we define
\beqq
\t{ind}_n(\alpha,f,A)=
\Big(f|_{V^E_0},
f(\bigoplus_{\substack{\lambda>0 \\ \alpha(\lambda)=0}} V^E_\lambda)
\oplus (A\cap V^{E'}_0)\Big).
\eeqq
The right hand side is indeed a morphism from $\t{ind}_n E$ to $\t{ind}_n E'$ 
in $\mc{V}_n$ because
\beqq
\t{ind}_n E'
&=&
V^{E'}_0 =
f(\bigoplus_{\alpha(\lambda)=0}V^E_\lambda) \oplus
\big((A\oplus\epsilon A)\cap V^{E'}_0\big) \\
&=& 
f(V^E_0)\oplus 
f(\bigoplus_{\substack{\lambda>0\\ \alpha(\lambda)=0}} V^E_\lambda)\oplus
f(\bigoplus_{\substack{\lambda<0\\ \alpha(\lambda)=0}} V^E_\lambda)\oplus 
(A\cap V^{E'}_0) \oplus 
(\epsilon A\cap V^{E'}_0) \\
&=& 
f(\t{ind}_n E)\oplus 
f(\bigoplus_{\substack{\lambda>0\\ \alpha(\lambda)=0}} V^E_\lambda)\oplus
\epsilon f(\bigoplus_{\substack{\lambda>0\\ \alpha(\lambda)=0}} V^E_\lambda)
\oplus 
\big(A\cap V^{E'}_0\big) \oplus 
\epsilon\big(A\cap V^{E'}_0\big).
\eeqq
An example (with $A=0$) is depicted in figure \ref{ind}. 
It is straightforward to check that $\t{ind}_n$ respects compositions. 

{\prop \label{SEFTn.Vn}
The functor $\t{ind}_n$ induces a homotopy equivalence 
\beqq
|\t{ind}_n|:|\mc{SEFT}_n|\stackrel{\sim}{\longrightarrow}|\mc{V}_n|.
\eeqq
}

\proof
We first define a functor $\mc{E}_n:\mc{V}_n\hookrightarrow\mc{SEFT}_n$
which identifies $\mc{V}_n$ with a subcategory of $\mc{SEFT}_n$. 
Given a graded $Cl_{-n}$-submodule $V\subset\mc{H}$, the $1$-SEFT 
$\mc{E}_n(V)$ is defined by $\mc{E}_n(V)(Z_{1,0})=\mc{H}$ and 
\beqq
\mc{E}_n(V)_B(I_{z,\theta})=\begin{cases}
  1 & \t{on  }V \\
  0 & \t{on  }V^\perp\subset\mc{H}
\end{cases}, 
\eeqq
for any object $B$ of $\mc{S}$ and 
$(z,\theta)\in\mc{O}(B)^\ev_{(0,\infty)}\times\mc{O}(B)^\od$. 
If $V\neq 0$, $\t{sp}\,Q_{\mc{E}_n(V)}=\{0\}$ and $V^{\mc{E}_n(V)}_0=V$; 
if $V=0$, $\t{sp}\,Q_{\mc{E}_n(V)}$ is empty.
For a morphism $(f,A)$ from an object $V$ to another object $V'$ in 
$\mc{V}_n$, we have 
\beqq
\mc{E}_n(f,A)=(\alpha_{V,V'},f,A),
\eeqq
where $\alpha_{V,V'}$ is the inclusion of $\t{sp}\,Q_{\mc{E}_n(V)}$ 
into $\t{sp}\,Q_{\mc{E}_n(V')}$.  

The composition $\t{ind}_n\mc{E}_n$ is the identity on $\mc{V}_n$, as 
it is clear that
\beqq
\begin{array}{lcl} 
\t{ind}_n\mc{E}_n(V)&=&V^{\mc{E}_n(V)}_0=V, \\
\t{ind}_n\mc{E}_n(f,A)&=&\t{ind}_n(\alpha_{V,V'},f,A)=(f,A).
\end{array}
\eeqq
On the other hand, there is a natural transformation $N$ from 
$\mc{E}_n\t{ind}_n$ to the identity on $\mc{SEFT}_n$, 
\beqq
N(E)=\Big(\alpha_E,i_E,\bigoplus_{\lambda>0}V^E_\lambda\Big),
\eeqq
where $\alpha_E$ is the inclusion of $\{0\}$ or $\emptyset$ into 
$\t{sp}\,Q_E$, depending on whether $0\in\t{sp}\,Q_E$ or not, 
and $i_E$ is the inclusion of $V^E_0$ 
(possibly $0$) into $\bigoplus_\lambda V^E_\lambda$. 
It is clear the right hand side above is indeed a morphism from 
$\mc{E}_n\t{ind}_n(E)=\mc{E}_n(V^E_0)$ to $E$. 
It remains to check the commutativity of the diagram 
\beqq
\xymatrix{
 \mc{E}_n\t{ind}_n(E)\ar[rr]^{N(E)}
   \ar[dd]_{\mc{E}_n\t{ind}_n(\alpha,f,A)} &&
 E \ar[dd]^{(\alpha,f,A)} \\ \\
 \mc{E}_n\t{ind}_n(E') \ar[rr]^{N(E')} && E'
}
\eeqq
for each morphism $(\alpha,f,A)$ in $\mc{SEFT}_n$. 
On one hand, we have 
\beqq
(\alpha,f,A)\circ N(E)=
\Big(\beta,f|_{V^E_0},f(\bigoplus_{\lambda>0}V^E_\lambda)\oplus A\Big),
\eeqq
where $\beta$ is the inclusion of $\{0\}$ or $\emptyset$ into 
$\t{sp}\,Q_{E'}$, depending on whether $0\in\t{sp}\,Q_E$ or not. 
On the other hand, 
\beqq
N(E')\circ\mc{E}_n\t{ind}_n(\alpha,f,A) 
= \Big(\beta, f|_{V^E_0},
f(\bigoplus_{\substack{\lambda>0 \\ \alpha(\lambda)=0}} V^E_\lambda)
\oplus (A\cap V^{E'}_0)
\oplus\bigoplus_{\lambda'>0}V^{E'}_{\lambda'}\Big).
\eeqq
Because $A=\bigoplus_{\lambda'\geq 0} A\cap V^E_{\lambda'}$ and 
$f(\bigoplus_{\alpha(\lambda)=\lambda'}V^E_\lambda)\oplus 
(A\cap V^{E'}_{\lambda'})=V^{E'}_{\lambda'}$ for each $\lambda'>0$, 
the two compositions are the same. 
We have shown that $\t{ind}_n$ induces a deformation retraction of 
$|\mc{SEFT}_n|$ onto $|\mc{V}_n|$.
\qedsymbol

\vspace{0.1in}
\noindent
{\bf Comparing $|\mc{V}_n|$ and $Fred_n$.}
To finish the proof of theorem \ref{SEFT.Fred}, it remains to establish the 
following

{\prop \label{V.Fred}
There is a homotopy equivalence
\beqq
|\mc{V}_n|\simeq Fred_n.
\eeqq
}

{\it Remark.}
In \cite{Segal.Khlgy}, Segal outlined the proofs of the homotopy equivalences
\beqq
|\wh{Vect}|\simeq Fred_0,\qquad |QVect|\simeq Fred_1.
\eeqq
According to propositions \ref{V0.virVect} and \ref{V1.QVect}, these imply 
the cases $n=0,1$ of proposition \ref{V.Fred}. 
In essence, we are going to generalize Segal's arguments.

\proof
Define a topological category $\mc{C}$ as follows. 
An object $(V,i)$ of $\mc{C}$ consists of an object $V$ of 
$\mc{V}_n$ and an even, $Cl_{-n}$-linear, isometric 
embedding $i:V\hookrightarrow\mc{H}$. 
Let $\t{Emb}'(V,\mc{H})$ denote the space of such embeddings. 
A morphism $(V,V',f,A,i)$ in $\mc{C}$ 
consists of two objects $V,V'$ of $\mc{V}_n$, 
a morphism $(f,A)$ of $\mc{V}_n$ from $V$ to $V'$, 
and some $i\in\t{Emb}'(V',\mc{H})$. 
Such a $5$-tuple is a morphism from $(V,if)$ to $(V',i)$. 
A pair of composable morphisms must be of the form 
$(V,V',f,A,if')$ and $(V',V'',f',A',i)$, 
whose composition is $(V,V'',f'f,f'(A)\oplus A',i)$. 

{\lemma \label{V.C0} 
$|\mc{V}_n|\simeq|\mc{C}|$.
}

\proof
There is a functor from $\mc{C}$ to $\mc{V}_n$ defined by 
$(V,i)\mapsto V$ and $(V,V',f,A,i)\mapsto(f,A)$. 
One easily sees that, for any $p\in\bb{Z}_+$, the induced map from 
$N_p\mc{C}$ to $N_p\mc{V}_n$ is a fiber bundle whose fibers 
are $\t{Emb}'(V,\mc{H})$ for various $V$. 
\footnote{
As usual, $N_\bullet$ denotes the nerve functor on small categories. 
}
Since $\t{Emb}'(V,\mc{H})\simeq *$, 
the functor induces a homotopy equivalence of the classifying spaces. 
\qedsymbol

Now, we define another topological category $\mc{C}'$. 
The objects of $\mc{C}'$ are the same as those of $\mc{V}_n$, 
i.e. finite dimensional, graded $Cl_{-n}$-submodules of $\mc{H}$. 
However, $N_0\mc{C}'$ is topologized as a subspace of $\t{Gr}(\mc{H})$. 
A morphism of $\mc{C}'$ is a pair 
$(V,A)\in N_0\mc{C}'\times\t{Gr}(\mc{H})$, where 
$A$ is a finite dimensional $Cl_{-n}$-submodule such that 
$A\perp V$ and $A\perp\epsilon A$; 
it is a morphism from $V$ to $V\oplus A\oplus\epsilon A$. 
Two composable morphisms must be of the form 
$(V,A)$ and $(V\oplus A\oplus\epsilon A,A')$, 
whose composition is $(V,A\oplus A')$. 

{\lemma \label{C0.C1}
$|\mc{C}|\simeq|\mc{C}'|$.
}

\proof
Let $\mc{C}_{\t{iso}}\subset\mc{C}$ be the subcategory with the same objects 
but only those morphisms of the form $(V,V',f,0,i)$. 
A point $x$ in $|\mc{C}_{\t{iso}}|$ is uniquely represented by 
the following data, 
\beq \label{point.in.Ciso}
\left(\xymatrix{
     V_0\ar[r]^{f_1} & V_1\ar[r]^{f_2} &
     \cdots\ar[r]^{f_p} & V_p\ar[r]^{i} & \mc{H}
   },\;
  \vec{s}\in\Delta^p-\partial\Delta^p\right),
\eeq
where $f_i$ are isomorphisms. 
We use the following notations for some of the data associated to $x$: 
$p^x=p$,\; $V^x_i=V_i$,\; $f^x_{\t{comp}}=f_p\cdots f_2 f_1$,\;
$i^x=i$,\; and $\vec{s}^x=\vec{s}$.

We define a category $\mc{D}$ such that 
$N_0\mc{D}=|\mc{C}_{\t{iso}}|$ and $|\mc{D}|\cong|\mc{C}|$. 
For $x,y\in|\mc{C}_{\t{iso}}|$, a morphism in $\mc{D}$ from $x$ to $y$ 
is a $4$-tuple $(x,f,A,y)$, 
where $(f,A)$ is a morphism in $\mc{V}_n$ from $V^x_{p^x}$ to $V^y_0$, 
such that $i^x=i^y\circ f^y_{\t{comp}}\circ f$. 
Composition is defined as follows,
\beqq
(y,f',A',z)\circ(x,f,A,y)=
\big(x,\;f'\circ f^y_{\t{comp}}\circ f,\;
  f'\circ f^y_{\t{comp}}(A)\oplus A',\;z\big).
\eeqq
We construct a map $\phi:|\mc{D}|\rightarrow|\mc{C}|$. 
The restriction of $\phi$ on $N_0\mc{D}=|\mc{C}_{\t{iso}}|$ is the map 
induced by the inclusion $\mc{C}_{\t{iso}}\subset\mc{C}$. 
For $q\geq 1$, a point in $N_q\mc{D}\times\Delta^q$ has the form 
$\left(x_0,f_1,A_1,x_1,f_2,A_2,\ldots,f_q,A_q,x_q;\vec{t}\;\right)$,
where $\vec{t}\in\Delta^q$ and the $i$-th morphism of the $q$-chain is 
$(x_{i-1},f_i,A_i,x_i)$. 
Its image under $\phi$ is defined to be  
\beqq
\left(
\xymatrix{
  V^{x_0}_\bullet \ar[r]^{(f_1,A_1)} &
  V^{x_1}_\bullet \ar[r]^{(f_2,A_2)} &
  \cdots \ar[r]^{(f_q,A_q)} &
  V^{x_q}_\bullet \ar[r]^{i^{x_q}} & \mc{H}  
},\; 
(t_0\vec{s}^{x_0},\ldots,t_q\vec{s}^{x_q})\right)\;\in|\mc{C}|.
\eeqq
It is straightforward to check that $\phi$ is well-defined, 
i.e. compatible with faces and degeneracies, 
and is a homemorphism. 

There is a functor $G:\mc{D}\rightarrow\mc{C}'$, defined by 
\beqq
N_0 G(x)=i^x(V^x_{p^x}),\qquad
N_1 G(x,f,A,y)=\big(i^x(V^x_{p^x}),i^y\circ f^y_{\t{comp}}(A)\big).
\eeqq
First, we show that $N_0 G$ is a fiber bundle with contractible fibers. 
\footnote{
In fact, $|\mc{C}_{\t{iso}}|\simeq  
\coprod_{[V]}BU^\ev_{Cl_{-n}}(V)\simeq N_0\mc{C}'$, 
where $V$ is any finite dimensional, graded $Cl_{-n}$-module, 
and the disjoint union has one component for each isomorphism 
class of such modules.
}
Fix $v\in N_0\mc{C}'$. 
Let $U^\ev_{Cl_{-n}}(\mc{H})$ be the space of even, $Cl_{-n}$-linear, 
unitary operators on $\mc{H}$. 
Since the map from $U^\ev_{Cl_{-n}}(\mc{H})$ to $N_0\mc{C}'$ given by 
$T\mapsto T(v)$ is a fiber bundle, 
it admits a section $V\mapsto T_V$ for $V$ in a neighhorhood 
$U\subset N_0\mc{C}'$ of $v$. 
Notice that $T_V(v)=V$. 
There is a homeomorphism $(N_0 G)^{-1}(U)\cong U\times(N_0 G)^{-1}(v)$;
for any $V\in U$, the image of $x\in (N_0 G)^{-1}(V)$ is $(V,x')$, 
where the data comprising $x'\in(N_0 G)^{-1}(v)$ 
-- see (\ref{point.in.Ciso}) -- 
are the same as those of $x$ except with $i^{x'}=T_V^{-1}\circ\,i^x$. 
The fiber $(N_0 G)^{-1}(v)$ is the classifying space of 
$\mc{C}_v\subset\mc{C}_{\t{iso}}$, which is the full subcategory 
consisting of objects $(V,i)$ such that $i(V)=v$. 
Any object of $\mc{C}_v$ is both initial and terminal, 
hence $(N_0 G)^{-1}(v)=|\mc{C}_v|\simeq *$. 

Similarly, $N_1 G$ is also a fiber bundle with contractible fibers. 
Fix $(v,a)\in N_1\mc{C}'$. 
Let $U$ be a neighborhood of $(v,a)$ that admits a continuous map to 
$U^\ev_{Cl_{-n}}(\mc{H})$, given by $(V,A)\mapsto T_{V,A}$, 
such that $T_{V,A}(v)=V$, $T_{V,A}(a)=A$, $\forall\,(V,A)\in U$. 
There is a homeomorphism $(N_1 G)^{-1}(U)\cong U\times(N_1 G)^{-1}(v,a)$ 
defined as follows:
$\forall\,(V,A)\in U$, $(x,f,A,y)\in(N_1 G)^{-1}(V,A)$ is mapped to 
$(x',f,A,y')\in (N_1 G)^{-1}(v,a)$, where the data comprising 
$x'$ (resp. $y'$) are the same as those of $x$ (resp. $y$) except with 
$i^{x'}=T_{V,A}^{-1}\circ i^x$ (resp. $i^{y'}=T_{V,A}^{-1}\circ i^y$). 
Now, consider the fiber $(N_1 G)^{-1}(v,a)$. 
Notice that, for any $(x,f,A,y)\in(N_1 G)^{-1}(v,a)$,  
$x\in(N_0 G)^{-1}(v)$ and $y\in(N_0 G)^{-1}(v')$, 
where $v'=v\oplus a\oplus \epsilon a$. 
Let $x_0=(v,\t{inc})\in N_0\mc{C}_v$ and 
$y_0=(v',\t{inc})\in N_0\mc{C}_{v'}$. 
Since all objects of $\mc{C}_v$ and $\mc{C}_{v'}$ are initial, 
there are canonical contractions of $(N_0 G)^{-1}(v)=|\mc{C}_v|$ and  
$(N_0 G)^{-1}(v')=|\mc{C}_{v'}|$ to the points 
$x_0$ and $y_0$ respectively. 
These two contractions uniquely determine a contraction of 
$(N_1 G)^{-1}(v,a)$ to the point $(x_0,\t{inc},a,y_0)$, 
so that $(N_1 G)^{-1}(v,a)\simeq *$. 

We have shown that $N_0 G$ and $N_1 G$ are homotopy equivalences. 
Since the face maps $N_1\mc{D}\rightrightarrows N_0\mc{D}$, 
$N_1\mc{C}'\rightrightarrows N_0\mc{C}'$ are fibrations, 
$|G|:|\mc{D}|\rightarrow|\mc{C}'|$ is also a homotopy equivalence. 
\qedsymbol

Let us define yet another topological category $\mc{C}''$ as follows: 
\beqq
N_0\mc{C}''
  &=& \big\{(V,F)\in N_0\mc{C}'\times Fred_n|
  F(V)\subset V,\ker F\subset V\big\}, \\
N_1\mc{C}''
  &=& \big\{(V,A,F)\in N_1\mc{C}'\times Fred_n|
  F(V)\subset V,F(A)\subset A,\ker F\subset V, F|_A>0 \big\}, 
\eeqq
where $(V,A,F)\in N_1\mc{C}''$ is a morphism from 
$(V,F)$ to $(V\oplus A\oplus\epsilon A,F)$. 
Hence, two morphisms are composable iff they are of the form 
$(V,A,F)$ and $(V\oplus A\oplus\epsilon A,A',F)$, 
and their composition is defined to be $(V,A\oplus A',F)$. 

{\lemma \label{C1.C2}
$|\mc{C}'|\simeq|\mc{C}''|$. 
}

\proof
Define a functor $\pi_1:\mc{C}''\rightarrow\mc{C}'$ by 
$(V,F)\mapsto V$, $(V,A,F)\mapsto(V,A)$. 
It suffices to show that $N_p\pi_1$ is 
a homotopy equivalence for every $p\geq 0$. 
A $p$-chain of $\mc{C}''$ whose $i$-th morphism is 
$(V_{i-1},A_i,F)$, $i=1,\ldots,p$, 
is determined by the $(p+2)$-tuple $(V_0,A_1,\ldots,A_p,F)$; 
we represent any $p$-chain of $\mc{C}''$ by such a $(p+2)$-tuple. 
Analogously, we represent a $p$-chain of $\mc{C}'$ by a $(p+1)$-tuple 
$(V,A_1,\ldots,A_p)$. 

Fix $(v,a_1,\ldots,a_p)\in N_p\mc{C}'$. 
Since $T\mapsto(T(v),T(a_1),\ldots,T(a_p))$ defines a fiber bundle 
$U^\ev_{Cl_{-n}}(\mc{H})\rightarrow N_p\mc{C}'$, 
$(v,a_1,\ldots,a_p)$ has a neighhorhood $U$ in $N_p\mc{C}'$ 
such that, $\forall\,C=(V,A_1,\ldots,A_p)\in U$, 
$\exists\,T_C\in U^\ev_{Cl_{-n}}(\mc{H})$ which 
depends continuously on $C$ and satisfies $T_C(v)=V$, 
$T_C(a_i)=A_i$, $i=1,\ldots,p$. 
Denote the vector space 
$v\oplus a_1\oplus\epsilon a_1\oplus\cdots\oplus a_p\oplus\epsilon a_p$ 
by $v_p$. 
Observe there is a homeomorphism 
\beq \label{Np.pi1}
\begin{array}{ccc}
(N_p\pi_1)^{-1}(U) &\cong& 
  U\times \t{End}'(v)\times\prod_{i=1}^p\t{End}'_+(a_i)\times GL'(v_p^\perp)\\
(C,F) &\mapsto& 
  (C,\;T_C^{-1}FT_C|_v,\;T_C^{-1}FT_C|_{a_1},\ldots,T_C^{-1}FT_C|_{a_p},\;
  T_C^{-1}FT_C|_{v_p^\perp}). 
\end{array}
\eeq
Here, $\t{End}'(v)$ is the vector space of odd, self-adjoint, 
$Cl_{-n}$-linear endomorphisms of $v$; 
$\t{End}'_+(a_i)$ is the space of positive definite, 
$Cl_{-n}$-linear endomorphisms of $a_i$; 
$GL'(v_p^\perp)$ is the space of invertible operators on $v_p^\perp$ 
which are odd, self-adjoint, $Cl_{-n}$-linear and,  
if $n\equiv 1\pmod{4}$, condition (\ref{ess.pos.neg}) is satisfied. 
Notice that $\mc{H}=v\oplus a_1\oplus\epsilon a_1\oplus\cdots\oplus 
a_p\oplus\epsilon a_p\oplus v_p^\perp$, hence the odd operator 
$T_C^{-1}FT_C$, or $F$, is determined by its restrictions on 
$v$, $a_1,\ldots,a_p$ and $v_p^\perp$. 
Because of the homeomorphism (\ref{Np.pi1}), 
$N_p\pi_1$ is a fiber bundle whose fiber is the product of 
$\t{End}'(v)$, $\t{End}'_+(a_i)$, $i=1,\ldots,p$, and $GL'(v_p^\perp)$. 

It remains to show that the fibers of $N_p\pi_1$ are contractible. 
Each $\t{End}'_+(a_i)$ is a convex subspace of a vector space, 
thus contractible. 
Let us consider $GL'(\mc{H}')$, where $\mc{H}'$ is any Hilbert space 
with a stable graded $Cl_{-n}$-action, e.g. $v_p^\perp$. 
In the case $n=0$, $T\mapsto T|_{{\mc{H}'}^\ev}$ defines a homeomorphism 
$GL'(\mc{H}')\cong GL({\mc{H}'}^\ev,{\mc{H}'}^\od)$, and the latter space is 
contractible by Kuiper's theorem \cite{Kuiper}. 
In the case $n<0$, 
let $e_1,\ldots,e_{|n|}$ be an orthonormal basis of $\bb{R}^{|n|}$, 
regarded as elements of $Cl_{-n}$. 
There is an action of $Cl^\ev_{-n}\cong Cl_{-n-1}$ on ${\mc{H}'}^\ev$. 
The map $T\mapsto (e_1 T)|_{{\mc{H}'}^\ev}$ defines a homeomorphism from 
$GL'(\mc{H}')$ to $GL''({\mc{H}'}^\ev)$, the space of skew-adjoint, 
$Cl_{-n-1}$-antilinear, 
\footnote{
An operator $T$ on a $Cl_k$-module $V$ is $Cl_k$-antilinear if $Tv=-vT$ for 
any $v\in\bb{R}^{|k|}\subset Cl_k$.
}
invertible operators $S$ on ${\mc{H}'}^\ev$ with the property that, if 
$n\equiv 1\pmod 4$, $e_1 e_2\cdots e_{|n|-1}S$ is not essentially positive or 
negative. 
It follows from Kuiper's theorem again that $GL''({\mc{H}'}^\ev)\simeq *$. 
\cite{AS.Fred} 
Finally, if $n>0$, suppose $n=4k-\ell$ where $k>0$ and $\ell\geq 0$ are 
integers. 
Let $M$ be an irreducible graded $Cl_{\ell,\ell}$-module. 
Since $Cl_{0,n}\otimes Cl_{\ell,\ell}\cong Cl_{\ell,4k}\cong Cl_{4k+\ell}$ 
as a graded algebra, $\mc{H}'\otimes M$ admits a stable graded 
$Cl_{4k+\ell}$-action. 
The map in (\ref{map.btx.Fred}) defines a homeomorphism 
$GL'(\mc{H}')\cong GL'(\mc{H}'\otimes M)$, 
where the latter space has just been shown to be contractible. 
\qedsymbol

{\lemma \label{C2.Fred}
$|\mc{C}''|\simeq Fred_n$. 
}

\proof
Regarding $Fred_n$ as a topological category with only identity 
morphisms, there is a functor $\pi_2:\mc{C}''\rightarrow Fred_n$ 
defined by $(V,F)\mapsto F$, $(V,A,F)\mapsto\t{id}_F$. 
Our first step is to construct continuous local sections of 
$N_0\pi_2:N_0\mc{C}''\rightarrow Fred_n$. 
This will also give rise to local sections of $|\pi_2|$, 
via the inclusion $N_0\mc{C}''\rightarrow|\mc{C}''|$. 

Fix $F_0\in Fred_n$. 
Let $4c=||(F_0|_{\ker F_0^\perp})^{-1}||^{-1}
=\inf\big(\t{sp}\,F_0\cap(0,\infty)\big)$, which is a positive 
number, 
\footnote{
Since $\t{im}\,F_0=\ker F_0^\perp$, the restriction of $F_0$ to 
$\ker F_0^\perp$ has a bounded inverse by the open mapping theorem. 
}
and $U=\{F\in Fred_n:||F-F_0||<c\}$. 
For each $F\in U$, let $e_F$ be the resolution of $1_{\mc{H}}$ 
associated with $F$, and define 
\beq \label{VF}
p_F=e_F\big((-c,c)\big)\in B(\mc{H}),\qquad V_F=p_F(\mc{H}). 
\eeq 
Notice that 
\beq \label{VF.perp}
1-p_F=e_F\big((-\infty,-3c)\cup(3c,\infty)\big),\quad 
V_F^\perp=(1-p_F)(\mc{H}).
\eeq
(See \cite{Rudin}, \S 10.20.) 
We claim $F\mapsto(V_F,F)$ defines a continuous section of $N_0\pi_2$ 
over $U$. 
This means, for $F\in U$,  \\
\indent (A) $\dim V_F<\infty$, \\
\indent (B) $V_F\subset\mc{H}$ is a graded $Cl_{-n}$-submodule, \\
\indent (C) $F(V_F)\subset V_F$, \\
\indent (D) $\ker F\subset V_F$, \\
\indent (E) the map $F\mapsto p_F$ is continuous. \\
Among these, (B), (C) and (D) follow immediately from (\ref{VF}). 

Consider (A). 
It is true for $F=F_0$, because $V_{F_0}=\ker F_0$. 
To prove it for all $F\in U$, it suffices to show that, 
$\forall\,F_1,F_2\in U$, $p_{F_2}$ is injective on $V_{F_1}$. 
Let us assume the contrary, i.e. for some $F_1,F_2\in U$, 
$\exists\,u\in V_{F_1}$ such that $||u||=1$ and $p_{F_2}(u)=0$. 
Since $p_{F_1}(u)=u=(1-p_{F_2})(u)$, by (\ref{VF}) and (\ref{VF.perp}), 
we have $||F_1(u)||<c$ and $||F_2(u)||>3c$, hence 
$||(F_2-F_1)(u)||>2c$. 
However, $||F_i-F_0||<c$, $i=1,2$, imply $||F_2-F_1||<2c$, 
yielding a contradiction. 

To prove (E), we fix $F_1\in U$, $0<\epsilon<1$, and try to show  
$||p_F-p_{F_1}||\leq\epsilon$ for any  
$F\in U$ with $||F-F_1||<c\epsilon$. 
It suffices to obtain $||p_F(u)-p_{F_1}(u)||\leq\epsilon$ for the following 
two cases: \\
\indent (i) $u\in V_{F_1}$ is an eigenvector of $F_1$, $||u||=1$, \\
\indent (ii) $u\perp V_{F_1}$, $||u||=1$. \\
For case (i), suppose $F_1 u=\lambda u$ and $u=u'+u''$, where 
$u'=p_F(u)\in V_F$, $u''\perp V_F$. 
Notice that $|\lambda|<c$ and $p_{F_1}(u)-p_F(u)=u-u'=u''$. 
The desired inequality $||u''||\leq \epsilon$ follows from:
\beqq
\begin{array}{rclcl}
c\epsilon 
&>& ||F(u)-F_1(u)||
&=& ||F(u')+F(u'')-\lambda u'-\lambda u''|| \\
&\geq& ||F(u'')-\lambda u''||
&\geq& ||F(u'')||-|\lambda|\,||u''|| \\
&>& 3c||u''|| - c||u''||
&=& 2c||u''||, 
\end{array}
\eeqq
where the fourth inequality follows from (\ref{VF.perp}), and 
the second inequality from the fact that 
$F(u')-\lambda u'\in V_F$ and $F(u'')-\lambda u''\in V_F^\perp$ 
are orthogonal. 
For case (ii), again write $u$ as $u'+u''$, where $u'=p_F(u)$. 
Notice that $p_F(u)-p_{F_1}(u)=u'-0=u'\in V_F$. 
Let $v'$ be an arbitrary vector in $V_F$ with $||v'||=1$. 
By the proof of (A), $v'=p_F(v)$ for a unique $v\in V_{F_1}$. 
Write $v$ as $v'+v''$; notice that $v''\perp V_F$. 
From the proof of (i), we have 
$||v-v'||<\epsilon||v||/2<||v||/2$, implying $||v||<2$. 
Because $u\perp v$, we have $\langle u',v'\rangle+\langle u'',v''\rangle=0$, 
hence 
\beqq
|\langle u',v'\rangle|
\leq ||u''||\,||v''|| 
\leq ||p_{F_1}(v)-p_F(v)||
\leq \frac{\epsilon}{2}\,||v||
< \epsilon, 
\eeqq
where the third inequality follows from case (i). 
Since $||u'||$ equals the supremum of $|\langle u',v'\rangle|$ over all 
$v'\in V_{F_1}$ with $||v'||=1$, we have
$||p_F(u)-p_{F_1}(u)||=||u'||\leq \epsilon$, as desired. 
This finishes the proof of our claim that 
$F\mapsto(V_F,F)$ is a continuous section of $N_0\pi_2$, or $|\pi_2|$, 
over $U$. 

Next, keeping the same notations, 
we show that $|\pi_2|^{-1}(U)$ can be contracted to the subspace 
$s_U=\{(V_F,F)|F\in U\}$ along the fibers of $|\pi_2|$. 
Since open sets $U$ with such contractions form a base of 
$Fred_n$, this will prove $|\pi_2|$ is a homotopy equivalence 
by \cite{DoldThom}. 

Let $\mc{C}''_{U,s_U}\subset\mc{C}''_U$ be the full subcategories 
of $\mc{C}''$ such that, given $(V,F)\in N_0\mc{C}''$, 
\beqq
\begin{array}{lcl}
(V,F)\in N_0\mc{C}''_U &\Leftrightarrow& F\in U, \\
(V,F)\in N_0\mc{C}''_{U,s_U} &\Leftrightarrow& F\in U\t{ and }V\supset V_F. 
\end{array}
\eeqq
There is a functor 
$r:\mc{C}''_U\rightarrow\mc{C}''_{U,s_U}$, defined by 
\beqq
r(V,F)=(V+V_F,F),\qquad r(V,A,F)=(V+V_F,A\cap V_F^\perp,F).
\eeqq
Suppose $(V,F)\in N_0\mc{C}''_U$ and 
$p$ is the orthogonal projection onto $V$. 
Since $V$, $V^\perp$ and $V_F$ are invariant subspaces of $F$, 
the orthogonal projection onto $V+V_F=V\oplus(V_F\cap V^\perp)$ 
is given by $p+p_F(1-p)=p+p_F-pp_F$. 
Clearly, $(p,F)\mapsto p+p_F-pp_F$ is continuous, 
thus so is $(V,F)\mapsto V+V_F$. 
A similar argument applies to $(A,F)\mapsto A\cap V_F^\perp$. 
This shows $r$ is continuous. 
The composition $r\circ\t{inc}$ is the identity on 
$\mc{C}''_{U,s_U}$, 
where $\t{inc}$ is the inclusion of $\mc{C}''_{U,s_U}$ into 
$\mc{C}''_U$. 
On the other hand, for any $(V,F)\in N_0\mc{C}''$, 
let $A_{V,F}=e_F\big((0,\infty)\big)(V^\perp)$, 
which depends continuously on $(V,F)$. 
There is a natural transformation from the identity on $\mc{C}''_U$ to 
$\t{inc}\circ r$, given by $(V,F)\mapsto(V,V_F\cap A_{V,F},F)$. 
Indeed, for any morphism $(V,A,F)$ of $\mc{C}''_U$, the following diagram 
commutes:
\beqq
\xymatrix{
  (V,F) \ar[rrrr]^{(V,V_F\cap A_{V,F},F)} \ar[dd]_{(V,A,F)} &&&&
  (V+V_F,F) \ar[dd]^{(V+V_F,A\cap V_F^\perp,F)} \\ \\
  (V',F) \ar[rrrr]^{(V',V_F\cap A_{V',F},F)} &&&&
  (V'+V_F,F),
}
\eeqq
where $V'=V\oplus A\oplus\epsilon A$. 
Commutativity follows from the fact that there is at most one morphism 
between any two objects of $\mc{C}''$. 
Therefore, $|r|$ defines a deformation retraction of 
$|\pi_2|^{-1}(U)=|\mc{C}''_U|$ onto the subspace $|\mc{C}''_{U,s_U}|$. 

Finally, any object $(V,F)$ of $\mc{C}''_{U,s_U}$ admits 
a unique morphism $(V_F,V\cap A_{V_F,F},F)$ from $(V_F,F)$, 
which depends continuously on $(V,F)$. 
In other words, each $(V_F,F)\in s_U$ is an initial object of a preimage 
category $\pi_2^{-1}(F)$, and the associated contraction of 
$|\pi_2^{-1}(F)|=|\pi_2|^{-1}(F)$ to the point $(V_F,F)$ depends 
continuously on $F\in U$. 
Therefore, $|\mc{C}''_{U,s_U}|$ deformation retracts onto $s_U$. 
Notice that the fibers of $|\pi_2|$ are preserved during both of the 
deformations constructed above. 
\qedsymbol

This completes the proof of proposition \ref{V.Fred}.
\qedsymbol

\newpage
\section{ANNULAR FIELD THEORIES AND THE ELLIPTIC COHOMOLOGY
FOR THE TATE CURVE}
\label{chap.AFT}


\setcounter{equation}{0}
\subsection{Annular Field Theories}

{\bf Some super manifolds of dimensions $(1|1)$ and $(2|1)$.}
Consider the $(1|1)$-manifold 
\beqq
S=\bb{R}/\bb{Z}\times\bb{R}^{0|1}. 
\eeqq
Let $(s,\lambda)$ be the standard coordinates on $S$, 
where $s$ is defined only locally. 
The metric $ds+\lambda d\lambda$ on $\bb{R}^{1|1}$ defined right above 
(\ref{Hom.00}) descends to a metric 
\beqq
\omega^S=ds+\lambda d\lambda
\eeqq
on $S^{1|1}$. 
Also consider the $(2|1)$-manifolds 
\beqq
A=\bb{R}/\bb{Z}\times\bb{R}\times\bb{R}^{0|1},\qquad
A_t=\bb{R}/\bb{Z}\times[0,t]\times\bb{R}^{0|1},
\eeqq 
with standard coordinates $(s,s',\lambda)$. 
Let 
\beqq
\omega^A_1=ds+ids',\qquad \omega^A_2=ds-ids'+\lambda d\lambda,
\eeqq 
which are complex $1$-forms on $A$ or $A_t$. 
Define the following maps between these super manifolds:
\beqq
\begin{array}{rll}
& \tau_t:S\rightarrow S, & \tau_t^*(s,\lambda)=(s+t,\lambda), \\
& \nu_{t,t'}:S\rightarrow A, & \nu_{t,t'}^*(s,s',\lambda)=(s+t,t',\lambda), \\
& \kappa_{t,t'}:A\rightarrow A, 
  & \nu_{t,t'}^*(s,s',\lambda)=(s+t,s'+t',\lambda), \\
& \epsilon:S\rightarrow S, & \epsilon^*(s,\lambda)=(s,-\lambda), \\
\t{or} & \epsilon:A\rightarrow A, 
  & \epsilon^*(s,s',\lambda)=(s,s',-\lambda),
\end{array}
\eeqq
where $t\in\bb{R}/\bb{Z}$, $t'\in\bb{R}$.  
Simple calculations then show the following:
\beq
&\t{Hom}(S,\omega^S;S,\omega^S)=\{\tau_t,\tau_t\epsilon\},& 
  \label{Hom.SS} \\
&\t{Hom}(S,ds,\omega^S;A,\omega^A_1,\omega^A_2)
  =\{\nu_{t,t'},\nu_{t,t'}\epsilon\},& 
  \label{Hom.SA} \\
&\t{Hom}(A,\omega^A_1,\omega^A_2;A,\omega^A_1,\omega^A_2)
  =\{\kappa_{t,t'},\kappa_{t,t'}\epsilon\}.& 
  \label{Hom.AA}
\eeq
The notation in e.g. (\ref{Hom.SA}) refers the set of maps from $S$ to $A$ 
that pull back $\omega^A_1$ and $\omega^A_2$ to $ds=\omega^S|_{S_\r}$ and 
$\omega^S$ respectively. 

Let $B$ be a super manifold with $B_\r=$ a point, and $\Lambda=\mc{O}(B)$. 
The natural transformations $(s_B,\lambda_B)$ and $(s_B,s'_B,\lambda_B)$ 
--- see (\ref{fcn.B}) --- allow us to make the identifications
\beqq
S(B)=\Lambda^\ev/\bb{Z}\times\Lambda^\od, \quad
A(B)=\Lambda^\ev/\bb{Z}\times\Lambda^\ev\times\Lambda^\od,
\eeqq
where $\Lambda^\ev/\bb{Z}$ is a quotient of abelian groups. 
Defining
\beqq
\begin{array}{rlll}
& \tau_x:S(B)\rightarrow S(B), 
  & (w,\eta)\mapsto(w+x,\eta),\\
& \nu_{x,y,\theta}:S(B)\rightarrow A(B), 
  & (w,\eta)\mapsto(w+x-\theta\eta/2,y-i\theta\eta/2,\theta+\eta),\\
& \kappa_{x,y,\theta}:A(B)\rightarrow A(B), 
  & (w,w',\eta)\mapsto(w+x-\theta\eta/2,w'+y-i\theta\eta/2,\theta+\eta),\\
& \epsilon:S(B)\rightarrow S(B), & (w,\eta)\mapsto(w,-\eta), \\
\t{or} & \epsilon:A(B)\rightarrow A(B), 
  & (w,w',\eta)\mapsto(w,w',-\eta),
\end{array}
\eeqq
\footnote{
More precisely, $\nu_{x,y,\theta}$ and $\kappa_{x,y,\theta}$ are maps 
into $A(B)\otimes\bb{C}$.  
}
with $x\in\Lambda^\ev/\bb{Z}$, $y\in\Lambda^\ev$, $\theta\in\Lambda^\od$, 
we have 
\beq
&&\t{Hom}\big(S(B),ds_B,\omega^S_B;S(B),ds_B,\omega^S_B\big)
  =\{\tau_x,\tau_x\epsilon\}, \label{BHom.SS} \\
&&\t{Hom}\big(S(B),ds_B,\omega^S_B;A(B),(\omega^A_1)_B,(\omega^A_2)_B\big)
  =\{\nu_{x,y,\theta},\nu_{x,y,\theta}\epsilon\}, 
  \label{BHom.SA} \\
&&\t{Hom}\big(A(B),(\omega^A_1)_B,(\omega^A_2)_B,
  A(B),(\omega^A_1)_B,(\omega^A_2)_B\big)
  =\{\kappa_{x,y,\theta},\kappa_{x,y,\theta}\epsilon\}. 
  \label{BHom.AA} 
\eeq
Recalling (\ref{1form.B}), the notation in e.g.~(\ref{BHom.SA}) refers to 
the set of maps $\nu:S(B)\rightarrow A(B)$ such that 
$(\omega^A_1)_B\circ T\nu=ds_B$ and $(\omega^A_2)_B\circ T\nu=\omega^S_B$. 
The proofs of (\ref{BHom.SS})-(\ref{BHom.AA}) are similar to those 
of (\ref{BHom.00})-(\ref{BHom.11}). 

There are analogues of lemmas \ref{BHom.red} and \ref{BHom.Hom} for 
the maps in (\ref{BHom.SS})-(\ref{BHom.AA}). 
More precisely, if $\phi:M_1(B)\rightarrow M_2(B)$ is any of those maps, 
there is a map $\phi_\r:(M_1)_\r\rightarrow(M_2)_\r$ of ordinary manifolds 
such that $\phi(f)_\r=\phi_\r f_\r$, $\forall f\in M_1(B)$, and also a map 
$\ell(\phi):M_1\rightarrow M_2$ of super manifolds such that 
$\ell(\phi)_\r=\phi_\r$. 
The latter ones are defined as follows:
\beq \label{ell.AFT}
\begin{array}{ll}
\ell(\tau_x)=\tau_{x^0}, & 
\ell(\tau_x\epsilon)=\tau_{x^0}\epsilon, \\
\ell(\nu_{x,y,\theta})=\nu_{x^0,y^0}, &
\ell(\nu_{x,y,\theta}\epsilon)=\nu_{x^0,y^0}\epsilon, \\
\ell(\kappa_{x,y,\theta})=\kappa_{x^0,y^0}, &
\ell(\kappa_{x,y,\theta}\epsilon)=\kappa_{x^0,y^0}\epsilon, 
\end{array}
\eeq
where e.g $x^0$ is the unique real number so that $x-x^0$ is nilpotent.  
The induced maps $\phi_\r$ and $\ell(\phi)$ respect compositions, i.e. 
\beqq
(\phi'\phi)_\r=\phi'_\r\phi_r,\quad
\ell(\phi'\phi)=\ell(\phi')\ell(\phi),
\eeqq
whenever $\phi$ and $\phi'$ are composable. 

Define $C(S)$ to be the Clifford algebra on the space of constant odd 
functions on $S$, i.e. $\bb{R}\lambda\subset\mc{O}(S)^\od$, with 
respect to the quadratic form $c\lambda\mapsto c^2$. 
Clearly, $C(S)\cong Cl_1$. 
Also, for each $t>0$, let $F(A_t)$ be the exterior algebra on the space 
of constant odd functions on $A_t$. 
\footnote{
$C(S)$ and $F(A_t)$ are small parts of the corresponding Clifford algebra 
and Fock space defined in \cite{ST}. 
}
$C(S)$ and $F(A_t)$ have the similar properties as the Clifford algebras and 
Fock spaces in section \ref{sec.1SEFT}. 
Let us repeat them below in terms of our $(1|1)$- and $(2|1)$-manifolds. 
For details, refer to \cite{ST}. 

Consider a \emph{bordism}
\beqq
\xymatrix{
& (A_t,\omega^A_1,\omega^A_2) & \\
(S,ds,\omega^S) \ar[ru]^\alpha && (S,ds,\omega^S) \ar[lu]_\beta,
}
\eeqq
i.e. $\alpha_\r$ (resp. $\beta_\r$) maps $S_\r$ onto the $s'=0$ 
(resp. $s'=t$) boundary of $A_t$; 
the notations in the diagram mean that $\alpha$ and $\beta$ belong to 
(\ref{Hom.SA}). 
Any such bordism canonically induces a $C(S)$-$C(S)$-action on $F(A_t)$. 
In fact, $F(A_t)$ is irreducible and is generated by a special element 
$\Omega$. 
It also has the property that $\lambda\Omega=\Omega\lambda$, where 
$\lambda\in C(S)$. 
If there is a commutative diagram
\beqq
\xymatrix{
& (A_t,\omega^A_1,\omega^A_2) \ar[dd]_\kappa & \\
(S,ds,\omega^S) \ar[ru]^\alpha \ar[rd]_{\alpha'} 
  && (S,ds,\omega^S) \ar[lu]_\beta \ar[ld]^{\beta'} \\
& (A_t,\omega^A_1,\omega^A_2) & \\
}
\eeqq
where the map $\kappa$ is invertible, then each of the two bordisms in the 
diagram gives rise to a $C(S)$-$C(S)$-action on $F(A_t)$, 
and $\kappa$ induces an isomorphism 
\beqq
\kappa_*:F(A_t)\stackrel{\sim}{\rightarrow}F(A_t)
\eeqq
of the two $C(S)$-$C(S)$-bimodules. 
Notice that by (\ref{Hom.AA}), $\kappa=1$ or $\epsilon$. 
Finally, suppose there is a diagram (with the $1$-forms suppressed)
\beqq
\xymatrix{
&& A_{t_1+t_2} && \\
& A_{t_1} \ar[ru]^\iota && A_{t_2} \ar[lu]_{\iota'} & \\
S \ar[ru]^\alpha && S \ar[lu]_\beta \ar[ru]^{\alpha'} && 
  S, \ar[lu]_{\beta'}
}
\eeqq
where the four maps on the bottom define two bordisms and the square is a 
pushout. 
The maps $\iota\alpha$ and $\iota'\beta'$ then also define a bordism, and 
there is a canonical isomorphism of $C(S)$-$C(S)$-bimodules,
\beq \label{Fock.gluing.AFT}
F(A_{t_2})\otimes_{C(S)} F(A_{t_1})\cong F(A_{t_1+t_2}),
\eeq
which sends $\Omega\otimes\Omega$ to $\Omega$. 

\vspace{0.1in}
\noindent
{\bf Field theories defined using $S$ and $A_t$.}
We now construct a super category $\mc{SAB}_n$ for each $n\in\bb{Z}$. 
It is a `super' version of a $(1+1)$-dimensional bordism category whose 
only bordisms are annuli. 
This should be thought of as a subcategory of a super category that 
incorporates all compact surfaces with conformal structures. 

Fix a super manifold $B$ with $B_\r=$ a point, and let $\Lambda=\mc{O}(B)$. 
The only object of $\mc{SAB}_n(B)$ is the $(1|1)$-manifold $S$. 
Like $\mc{SEB}^1_n(B)$, the morphisms are divided into 
`decorated $B$-isomorphisms' and equivalence classes of 
`decorated $B$-bordisms.' 
A \emph{decorated $B$-isomorphism} $(\tau,c)$ consists of an invertible map 
\beqq
\tau:\big(S(B),\omega^S_B\big)\rightarrow\big(S(B),\omega^S_B\big),
\eeqq 
together with an element $c\in C(S)^{\otimes -n}$. 
(See section \ref{sec.1SEFT}, footnote \ref{neg.power.tensor}.)
A \emph{decorated $B$-bordism} is a $4$-tuple $(A_t,\alpha,\beta,\Psi)$ 
consisting of a $(2|1)$-manifold $A_t$, maps 
\beqq
\xymatrix{
& \big(A_t(B),(\omega^A_1)_B,(\omega^A_2)_B\big) & \\
\big(S(B),ds_B,\omega^S_B\big) \ar[ru]^\alpha &&
\big(S(B),ds_B,\omega^S_B\big), \ar[lu]_\beta 
}
\eeqq
and an element $\Psi\in F(A_t)^{\otimes -n}$. 
This $4$-tuple is \emph{equivalent} to those of the form 
$(A_t,\kappa\alpha,\kappa\beta,\ell(\kappa)_*\Psi)$, where 
\beqq
\kappa:\big(A_t(B),(\omega^A_1)_B,(\omega^A_2)_B\big)\rightarrow
  \big(A_t(B),(\omega^A_1)_B,(\omega^A_2)_B\big)
\eeqq
is an invertible map, and $\ell$ is as in (\ref{ell.AFT}). 
The equivalence class is denoted $[A_t,\alpha,\beta,\Psi]$. 
Composition in $\mc{SAB}^1_n(B)$ is defined as follows:
\beqq
\begin{array}{lcl}
(\tau_2,c_2)\circ(\tau_1,c_1)&=&(\tau_2\tau_1,\ell(\tau_1)^*c_2\cdot c_1) \\
(\tau,c)\circ[A_t,\alpha,\beta,\Psi]&=&
  [A_t,\alpha,\beta\tau^{-1},\ell(\tau^{-1})^*c\cdot\Psi] \\
\left[A_t,\alpha,\beta,\Psi\right]\circ(\tau,c)&=&
  [A_t,\alpha\tau,\beta,\Psi c] \\
\left[A_{t_2},\alpha_2,\beta_2,\Psi_2\right]\circ
  [A_{t_1},\alpha_1,\beta_1,\Psi_1] &=& 
  [A_{t_1+t_2},\iota_1\alpha_1,\iota_2\beta_2,\Psi_3]. 
\end{array}
\eeqq
In e.g. the second case, the 
$C(S)^{\otimes -n}$-$C(S)^{\otimes -n}$-action on $F(A_t)^{\otimes -n}$ 
is induced from the bordism (not $B$-bordism!) given by $\ell(\alpha)$ and 
$\ell(\beta\tau^{-1})$. 
For the last case, $\iota_1$, $\iota_2$ are maps such that the square 
in the following diagram ($1$-forms suppressed),
\beqq
\xymatrix{
&& A_{t_1+t_2}(B) && \\
& A_{t_1}(B) \ar[ru]^{\iota_1} && A_{t_2}(B) \ar[lu]_{\iota_2} & \\
S(B) \ar[ru]^{\alpha_1} && S(B) \ar[lu]_{\beta_1} \ar[ru]^{\alpha_2} && 
  S(B), \ar[lu]_{\beta_2}
}
\eeqq
is a pushout; 
$\Psi_3$ is the image of $\Psi_2\otimes\Psi_1$ 
under the isomorphism (\ref{Fock.gluing.AFT}), 
which is induced by applying $\ell$ to the above diagram. 

The decorated $B$-isomorphisms in $\mc{SAB}^1_n(B)$ are classified by 
(\ref{BHom.SS}) and the fact that $C(S)^{\otimes -n}\cong Cl_{-n}$. 
The relations among decorated $B$-isomorphisms are given by 
\beq 
(1,c)\circ(\epsilon,1)=(\epsilon,1)\circ(1,\epsilon^*c), 
  \label{SAB.Clifford.epsilon} \\
(1,c)\circ(\tau_x,1)=(\tau_x,1)\circ(1,c).
  \label{SAB.Clifford.rotation}
\eeq
Given a decorated $B$-bordism $(A_t,\alpha,\beta,\Psi)$, 
$\alpha$ and $\beta$ are elements of (\ref{BHom.SA}), and there exists a 
unique element $\kappa$ of (\ref{BHom.AA}) such that 
$\kappa\alpha=\nu_{0,0,0}$, while $\kappa\beta$ can be 
any element of (\ref{BHom.SA}) with $y^0>0$. 
Therefore, all the equivalence classes of decorated $B$-bordisms are 
precisely as follows: 
\beqq
A_{x,y,\theta,\Psi}=[A_{y^0},\nu_{0,0,0},\nu_{x,y,\theta},\Psi], \qquad
(\epsilon,1)\circ A_{x,y,\theta,\Psi}
  =[A_{y^0},\nu_{0,0,0},\nu_{x,y,\theta}\epsilon,\Psi],
\eeqq
where $(x,y,\theta)\in
\Lambda^\ev/\bb{Z}\times\Lambda^\ev_{(0,\infty)}\times\Lambda^\od$ 
and $\Psi\in F(A_{y^0})$. 
Let us define $A_{x,y,\theta}:=A_{x,y,\theta,\Omega^{\otimes -n}}$. 
Notice that any $A_{x,y,\theta,\Psi}$ can be rewritten as 
$(1,c_1)\circ A_{x,y,\theta}\circ(1,c_2)$, for some 
$c_1,c_2\in C(S)^{\otimes -n}$. 
The relations between decorated $B$-isomorphisms and decorated 
$B$-bordisms are given by 
\beq
(\tau_{x_1},1)\circ A_{x_2,y,\theta}
  &=& A_{x_2,y,\theta}\circ(\tau_{x_1},1)=A_{-x_1+x_2,\;y,\;\theta}, 
  \label{SAB.rotation.bordism}  \\
(\epsilon,1)\circ A_{x,y,\theta} &=& A_{x,y,-\theta}\circ(\epsilon,1), 
  \label{SAB.epsilon.bordism} \\
(1,c)\circ A_{x,y,\theta} &=& A_{x,y,\theta}\circ (1,c).
  \label{SAB.Clifford.bordism}
\eeq
The first two follow from simple calculations and the last one from the 
fact that $\lambda\Omega=\Omega\lambda$. 
Finally, a diagram similar to the one that proves (\ref{sbordism.compose}) 
shows 
\beq \label{SAB.bordism.bordism}
A_{x_2,y_2,\theta_2}A_{x_1,y_1,\theta_1} =
A_{x_1+x_2-\frac{1}{2}\theta_1\theta_2,\;
  y_1+y_2-\frac{i}{2}\theta_1\theta_2,\;\theta_1+\theta_2}.
\eeq

Now, given a map of super manifolds $f:B'\rightarrow B$, the functor 
$\mc{SAB}^1_n(f):\mc{SAB}^1_n(B)\rightarrow\mc{SAB}^1_n(B')$ 
is defined by
\beqq
S\mapsto S,\quad
(\epsilon,1)\mapsto(\epsilon,1),\quad
(\tau_x,c)\mapsto(\tau_{f^*x},c),\quad
A_{x,y,\theta}\mapsto A_{f^*x,f^*y,f^*\theta},
\eeqq
where $f^*:\Lambda\rightarrow\Lambda'=\mc{O}(B')$. 
We indeed have a functor because the morphisms of $\mc{SAB}^1_n(B)$ 
are generated by $(\epsilon,1)$ and those of the forms $(\tau_x,c)$ and 
$A_{x,y,\theta}$. 
This finishes the definition of the super category $\mc{SAB}^1_n$. 

\defn
An \emph{annular field theory, or AFT, of degree $n$} is a functor 
$E:\mc{SAB}^1_n\rightarrow\mc{S}\t{Hilb}^{\bb{C}}$ of super categories. 
The super category $\mc{S}\t{Hilb}^{\bb{C}}$ is defined similarly as 
$\mc{S}\t{Hilb}$, but with complex Hilbert spaces and complex-linear 
operators. 
Let $E_B:\mc{SAB}^1_n(B)\rightarrow\mc{S}\t{Hilb}(B)$ be the (ordinary) 
functors comprising $E$, where $B$ is any object of $\mc{S}$.  
Since the values of $E_B$ on $S$, $(\epsilon,1)$ and $(1,c)$ 
are independent of $B$,
\footnote{
See definition \ref{defn.SEFT0}.
}
we simply write $\mc{H}:=E(S)$, $E(\epsilon,1)$ and $E(1,c)$.
Fix $B$ and let $\Lambda=\mc{O}(B)$. 
$E_B$ is required to satisfy the following assumptions:
\begin{itemize}
\item[(i)]
$E_B(A_{x,y,\theta})\in\Lambda\otimes B(\mc{H})$ are Hilbert-Schmidt 
operators.
\item[(ii)]
$E_B(\tau_x,1)^\dagger=E_B(\tau_{-x},1)$ and 
$E_B(A_{x,y,\theta})^\dagger=E_B(A_{-x,y,-i\theta})$.
\footnote{
We make this assumption on $E_B(A_{x,y,\theta})^\dagger$ 
so that $E_B(A_{x,y,0})^\dagger=E_B(A_{-x,y,0})$, 
as expected according to the usual axioms of field theories, 
and taking adjoints respects (\ref{SAB.bordism.bordism}). 
To check the latter, one should beware that 
$[(\theta T)(\theta'T')]^\dagger=-(\theta'T')^\dagger(\theta T)^\dagger$, 
for $T,T'\in B(\mc{H})$. 
}
\item[(iii)]
$\epsilon=E(\epsilon,1)$ defines the $\bb{Z}/2$-grading in $\mc{H}$, i.e. 
$\epsilon|_{\mc{H}^\ev}=1$, $\epsilon|_{\mc{H}^\od}=-1$. 
\item[(iv)]
$E_B(\tau_x,c)$ is analytic in $x$ and linear in $c$.
$E_B(A_{x,y,\theta,\Psi})$ is analytic in $(x,y,\theta)$ and linear 
in $\Psi$. 
\end{itemize}

Let us analyze the consequences of the definition of an AFT $E$: 
\begin{itemize}
\item
By (iii) and (\ref{SAB.Clifford.epsilon}), 
the operators $E(1,c)$, for $c\in C(S)^{\otimes -n}\cong Cl_{-n}$, 
define a $\bb{Z}/2$-graded $Cl_{-n}$-action on $\mc{H}$. 
\item
By (ii) and $\tau_x\tau_{-x}=1$, 
$x\in\bb{R}/\bb{Z}\mapsto E_B(\tau_x,1)$ defines a unitary $S^1$-action 
on $\mc{H}$, so that $E_B(\tau_x,1)=e^{2\pi ixP}$, where $P$ is even, 
self-adjoint, $Cl_{-n}$-linear, because of (\ref{SAB.Clifford.rotation}).
This extends to all $x\in\Lambda^\ev/\bb{Z}$ because of (iv).  
\item
Since (i), (ii) and (\ref{SAB.bordism.bordism}) together 
imply $y\mapsto E_B(A_{0,y,0})$ is a commutative semigroup of 
self-adjoint, Hilbert-Schmidt operators, we have 
$E_B(A_{0,y,0})|_{\mc{H}'}=e^{-2\pi y H}$ and 
$E_B(A_{0,y,0})|_{\mc{H}'^\perp}=0$, where 
$\mc{H}'\subset\mc{H}$ is a subspace and 
$H$ is a self-adjoint operator on $\mc{H}'$ with compact resolvent. 
Again, this holds for all $y\in\Lambda^\ev_{(0,\infty)}$. 
Furthermore, due to (\ref{SAB.epsilon.bordism}) and 
(\ref{SAB.Clifford.bordism}), $\mc{H}'$ is in fact a $\bb{Z}/2$-graded 
$Cl_{-n}$-submodule, and $H$ is even and $Cl_{-n}$-linear. 
For simplicity, we just write $E_B(A_{0,y,0})=e^{-2\pi y H}$, 
thinking that `$H=+\infty$' on $\mc{H}'^\perp$. 
\item
The relation in (\ref{SAB.rotation.bordism}) implies that 
$P$ and $H$ commute.
\item
Let $E_B(A_{0,y,\theta})=e^{-2\pi y H}+\sqrt{2\pi i}\theta K(y)$. 
It follows from (\ref{SAB.epsilon.bordism}), (ii) and 
(\ref{SAB.Clifford.bordism}) that $K(y)$ is odd, self-adjoint and 
$Cl_{-n}$-linear. 
Applying $E_B$ to (\ref{SAB.bordism.bordism}) with $x_1=x_2=0$ and using 
(\ref{SAB.rotation.bordism}), we obtain 
\beq \label{eqn.K}
\begin{array}{rl}
& \sqrt{2\pi i}\,\big[\theta_1 e^{-2\pi y_2 H} K(y_1)
  + \theta_2 K(y_2) e^{-2\pi y_1 H}\big]
  + 2\pi i\theta_1\theta_2 K(y_2)K(y_1) \\
=& \sqrt{2\pi i}(\theta_1+\theta_2)K(y_1+y_2-\frac{i}{2}\theta_1\theta_2)
  + \pi i\theta_1\theta_2(H+P).
\end{array} 
\eeq
Putting $\theta_1=0$ and $\theta_2=0$ into (\ref{eqn.K}) separately shows 
\beqq
K(y_2)e^{-2\pi y_1 H}=K(y_1+y_2)=e^{-2\pi y_2 H}K(y_1).
\eeqq
In particular, all $K(y)$ commute with $H$. 
Let $\tilde{G}=K(1)e^{2\pi H}$. 
For $y>1$, we have $K(y)=K(1)e^{-2\pi(y-1)H}=\tilde{G}e^{-2\pi yH}$; 
this in fact extends to all $y\in\Lambda^\ev_{(0,\infty)}$ by (iv). 
Finally, putting this back into (\ref{eqn.K}) yields 
$2\tilde{G}^2=H+P$. 
\end{itemize}
To summarize, let $q=e^{2\pi i(x+iy)}$ and $\bar{q}=e^{-2\pi i(x-iy)}$; 
using (\ref{SAB.rotation.bordism}), we have 
\beq \label{AFT.generators}
E_B(A_{x,y,\theta})
  =q^L\bar{q}^{\tilde{G}^2}(1+\sqrt{2\pi i}\theta\tilde{G}),
\eeq
\footnote{
More conceptually, an AFT is determined by two operators because 
the $B$-bordisms, as $B$ varies, define a super semigroup by 
(\ref{SAB.bordism.bordism}), and its super Lie algebra is generated by 
two elements. 
}
$\forall (x,y,\theta)\in\Lambda^\ev/\bb{Z}\times
\Lambda^\ev_{(0,\infty)}\times\Lambda^\od$,  
where $L=-P+\tilde{G}^2$ is even, $\tilde{G}$ is odd, and 
both are self-adjoint, $Cl_{-n}$-linear and have compact resolvent. 
Furthermore, $L-\tilde{G}^2=-P$ has integer eigenvalues, 
and the eigenvalues of $L+\tilde{G}^2=H$ are bounded below. 
The pair $(L,\tilde{G})$ are called the \emph{generators} of $E$. 

\vspace{0.1in}
\noindent
{\bf A few words about AFTs.}
Our definitions of $\mc{SAB}_n$, and hence of AFTs, are unsatisfactory 
for a number of reasons. 
Obviously, we would like to include bordisms given by more general surfaces. 
Furthermore, these surfaces should be equipped with an appropriate 
notion of `super conformal structure.' 
The author has not been able to identify such a structure. 
The geometric structures attached to the super manifolds $S$ and $A_t$ 
in the definition of $\mc{SAB}_n$ should be particular representatives of 
certain `super conformal classes.' 
Also, our choices of the Clifford algebra $C(S)$ and the Fock spaces 
$F(A_t)$, though will give us the desired result below, are rather ad hoc. 


\setcounter{equation}{0}
\subsection{Categories $\mc{AFT}_n$}
For any AFT $E$, we denote its generators by $(L_E,\tilde{G}_E)$. 

\defn 
Let $n\in\bb{Z}$, and
\beqq
\mc{H}
=\bigoplus_{k\in\bb{Z}}\mc{H}_k
=\bigoplus_{k\in\bb{Z}}(\mc{H}_k^\ev\oplus\mc{H}_k^\od)
\eeqq
be a $\bb{Z}\times\bb{Z}/2$-graded Hilbert space, such that each 
$\mc{H}_k$ admits a stable (see footnote \ref{stable}) 
$\bb{Z}/2$-graded $Cl_{-n}$-action. 
We define a topological category $\mc{AFT}_n(\mc{H})$ as follows. 
An object of $\mc{AFT}_n(\mc{H})$ is an AFT $E$ of degree $n$ such that    
$E(S)=\mc{H}$,        
$\mc{H}_k=\{L_E-\tilde{G}_E^2=k\}$ and the given graded $Cl_{-n}$-action on 
$\mc{H}$ coincides with the one coming from $E$. 
In oder to define the morphisms in $\mc{AFT}_n(\mc{H})$, 
we need to first discuss the spectral decompositions of the generators of 
an AFT of degree $n$. 

\begin{figure}[t] 
\setlength{\unitlength}{0.037cm}
\begin{picture}(400,120)
  \put(0,20){\line(1,0){400}}
  \put(5,0){$\cdots$}
  \put(50,20){\circle*{3}}
  \put(48,0){\scriptsize $V^E_{-\mu,k}$}
  \put(150,20){\circle*{3}}
  \put(148,0){\scriptsize $V^E_{-\lambda,k}$}
  \put(200,20){\circle*{3}}
  \put(198,0){\scriptsize $V^E_{0,k}$}
  \put(250,20){\circle*{3}}
  \put(248,0){\scriptsize $V^E_{\lambda,k}$}
  \put(350,20){\circle*{3}}
  \put(348,0){\scriptsize $V^E_{\mu,k}$}
  \put(380,0){$\cdots$}

  \put(0,60){\line(1,0){400}}
  \put(5,40){$\cdots$}
  \put(130,60){\circle*{3}}
  \put(128,40){\scriptsize $V^E_{-\nu,k+1}$}
  \put(200,60){\circle*{3}}
  \put(198,40){\scriptsize $V^E_{0,k+1}$}
  \put(270,60){\circle*{3}}
  \put(268,40){\scriptsize $V^E_{\nu,k+1}$}
  \put(380,40){$\cdots$}

  \put(0,100){\line(1,0){400}}
  \put(5,80){$\cdots$}
  \put(150,100){\circle*{3}}
  \put(148,80){\scriptsize $V^E_{-\lambda,k+2}$}
  \put(250,100){\circle*{3}}
  \put(248,80){\scriptsize $V^E_{\lambda,k+2}$}
  \put(380,80){$\cdots$}
\end{picture}
\caption{The spectral decomposition of $(\tilde{G}_E,L_E-\tilde{G}_E^2)$ 
for an AFT $E$. \label{sp.L.G}}
\end{figure}
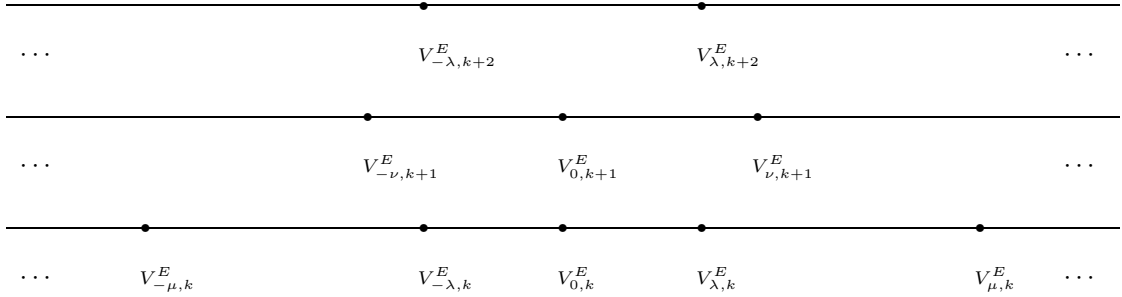

Recall from (\ref{AFT.generators}) that any AFT $E$ 
is determined by its generators 
$(L_E,\tilde{G}_E)$, which have the properties that  
(i) $L_E$ is even, $\tilde{G}_E$ is odd, 
(ii) $L_E$ and $\tilde{G}_E$ commute and are self-adjoint, 
(iii) they are $Cl_{-n}$-linear, 
(iv) they have compact resolvent, 
(v) $L_E-\tilde{G}_E^2$ has integer eigenvalues, and
(vi) $L_E$ is bounded below on $\ker\tilde{G}_E$.  
By (ii), the operators $L_E$ and $\tilde{G}_E$ have real spectra and 
a simultaneous spectral decomposition. 
Furthermore, (iv) implies their spectra are discrete and each eigenspace 
is finite dimensional. 
Let $V^E_{\lambda,k}$ be the common eigenspace of 
$\tilde{G}_E$, with eigenvalue $\lambda$, 
and $L_E-\tilde{G}_E^2$, with eigenvalue $k\in\bb{Z}$. 
By (iii), $V^E_{\lambda,k}$ are $Cl_{-n}$-submodules of $\mc{H}_k$, 
and (vi) says $V^E_{0,k}=0$ for $k<<0$. 
Finally, (i) implies $\epsilon V^E_{\lambda,k}=V^E_{-\lambda,k}$, where 
$\epsilon|_{\mc{H}^\ev}=1$ and $\epsilon|_{\mc{H}^\od}=-1$. 
Figure \ref{sp.L.G} shows a picture of these spectral data. 
In the picture, 
the lines represent the subset $\bb{R}\times\bb{Z}$ of $\bb{R}^2$. 
We regard $\t{sp}(\tilde{G}_E,L_E-\tilde{G}_E^2)$, 
the simultaneous spectrum of $\tilde{G}_E$ and $L_E-\tilde{G}_E^2$, 
as a subset of $\bb{R}\times\bb{Z}\subset\bb{R}^2$, and
label each of its points $(\lambda,k)$ with the corresponding 
eigenspace $V^E_{\lambda,k}$. 
As in the case of $1$-SEFTs, there is a symmetry about the $y$-axis in 
figure \ref{sp.L.G} due to the isomorphisms between 
$V^E_{\lambda,k}$ and $V^E_{-\lambda,k}$ defined by $\epsilon$. 

Now, we can define the morphisms in $\mc{AFT}_n(\mc{H})$. 
A morphism from an object $E$ to another object $E'$ is 
a `deformation' $(\alpha,f,A)$ from the spectral decomposition of 
$(\tilde{G}_E,L_E-\tilde{G}_E^2)$ to that of 
$(\tilde{G}_{E'},L_{E'}-\tilde{G}_{E'}^2)$, where 
\begin{itemize}
\item[-]
$\alpha:\t{sp}(\tilde{G}_E,L_E-\tilde{G}_E^2)\rightarrow
\t{sp}(\tilde{G}_{E'},L_{E'}-\tilde{G}_{E'}^2)$ is a proper map that 
is constant on the second coordinate, 
odd and order preserving in the first coordinate 
(regarding $\t{sp }(\tilde{G}_E,L_E-\tilde{G}_E^2)$, etc., as 
subsets of $\bb{R}^2$), 
\item[-] 
$f:\bigoplus_{(\lambda,k)}V^E_{\lambda,k}\rightarrow
\bigoplus_{(\lambda',k')}V^{E'}_{\lambda',k'}$ is an even map 
that embeds $V^E_{\lambda,k}$ isometrically into 
$V^{E'}_{\alpha(\lambda,k)}$, and 
\item[-] 
$A\subset\mc{H}$ is a $Cl_{-n}$-submodule, such that 
$\tilde{G}_{E'}\geq 0$ on $A$, $A\perp\epsilon A$ and 
$\bigoplus_{\lambda'}V^{E'}_{\lambda'}=f(\bigoplus_\lambda V^E_\lambda)
\oplus A\oplus\epsilon A$. 
\end{itemize}
Fixing $\alpha$ and 
$f\in\prod_{(\lambda,k)}\t{Hom}(V^E_{\lambda,k},V^{E'}_{\alpha(\lambda,k)})$, 
$A$ is determined by 
$A\cap V^{E'}_{0,k}\in\t{Gr}(V^{E'}_{0,k})$, $k\in\bb{Z}$. 
Indeed, $A\oplus\epsilon A$ is the orthogonal complement of the image of 
$f$, and the assumption that $\tilde{G}_{E'}\geq 0$ determines $A$ up to 
the subspace where $\tilde{G}_{E'}=0$. 
The set of morphisms $\{(\alpha,f,A)\}$ from $E$ to $E'$ is topologized 
as a subspace of 
\beqq
\coprod_\alpha\left(
\prod_{(\lambda,k)\in\t{sp}(\tilde{G}_E,L_E-\tilde{G}_E^2)}
\t{Hom}(V^E_{\lambda,k},V^{E'}_{\alpha(\lambda,k)})
\times\prod_{k\in\bb{Z}}\t{Gr}(V^{E'}_{0,k})\right). 
\eeqq
Composition of morphisms is defined by 
\beqq
(\alpha',f',A')\circ(\alpha,f,A)=(\alpha'\alpha,f'f,f'(A)\oplus A').
\eeqq 
This finishes the definition of $\mc{AFT}_n(\mc{H})$.

{\thm \label{AFT.TateK}
For each $n\in\bb{Z}$, we have 
\beqq
|\mc{AFT}_n(\mc{H})|\simeq (K_{\t{Tate}})_n. 
\eeqq
The right hand side is the $n$-th space of the elliptic spectrum 
$K_{\t{Tate}}$ associated with the Tate curve.
}

We briefly recall the Tate curve and $K_{\t{Tate}}$. 
For each $q\in\bb{C}$, $|q|<1$, there is an elliptic curve over $\bb{C}$ 
defined by a Weierstrass equation
\beqq
y^2+xy=x^3+a_4(q)x+a_6(q).
\eeqq
If $0<|q|<1$, the elliptic curve is isomorphic to 
$\bb{C}^\times/q^{\bb{Z}}$ as a Riemann surface, with the group structure 
induced from multiplication in $\bb{C}^\times$. 
As $q$ varies, $a_4$ and $a_6$ define analytic functions over the open disk 
$\{\,|q|<1\}$. 
Their Taylor series at $q=0$ have integer coefficients, hence the above 
equation also defines an elliptic curve $C_{\t{Tate}}$ over $\bb{Z}[[q]]$, 
the ring of integral formal power series. \cite{Silverman}
However, we will instead regard $C_{\t{Tate}}$ as an elliptic curve over 
$\bb{Z}[[q]][q^{-1}]$, the ring of integral formal Laurent series. 

An elliptic spectrum $(E,C,t)$ consists of an even periodic 
ring spectrum $E$, an elliptic curve $C$ over $\pi_0 E$, and an isomorphism 
$t:\t{Spf}\,E^0(\bb{C}P^\infty)\stackrel{\sim}{\rightarrow}\wh{C}$ 
of formal groups over $\pi_0 E$. \cite{AHS} 
On one hand, as a result of even periodicity, 
$E^0(\bb{C}P^\infty)$ is abstractly the ring of formal power series 
over $\pi_0 E$, and the formal scheme $\t{Spf}\,E^0(\bb{C}P^\infty)$ admits 
a formal group structure. 
On the other hand, the group structure on $C$ induces a formal group 
structure on $\wh{C}$, the formal neighborhood of $C$ at the 
identity. 
As an example, the elliptic spectrum associated with the Tate curve is 
\beqq
K_{\t{Tate}}=(K[[q]][q^{-1}],C_{\t{Tate}},t_{\t{Tate}}), 
\eeqq
where the cohomology theory associated with $K[[q]][q^{-1}]$ is periodic 
$K$-theory tensored with $\bb{Z}[[q]][q^{-1}]$.
The formal groups arising from $K[[q]][q^{-1}]$ and $C_{\t{Tate}}$ 
are both the multiplicative formal group. 

\proof[Proof of theorem \ref{AFT.TateK}]
Recall the categories $\mc{V}_n$ from section \ref{sec.Vn}. 
For each $k\in\bb{Z}$, define a functor 
$\t{ind}_{n,k}:\mc{AFT}_n(\mc{H})\rightarrow\mc{V}_n(\mc{H}_k)$ 
by $E\mapsto V^E_{0,k}$ and 
\beqq
(\alpha,f,A)\mapsto
  \Big(f|_{V^E_{0,k}},
  f(\bigoplus_{\substack{\lambda>0 \\ \alpha(\lambda,k)=(0,k)}} 
  V^E_{\lambda,k})\oplus (A\cap V^{E'}_{0,k})\Big).
\eeqq
Collecting $\t{ind}_{n,k}$ for all $k\in\bb{Z}$ yields another functor, 
\beqq
\t{ind}_n^{S^1}:\mc{AFT}_n(\mc{H})\rightarrow
{\prod_{k\in\bb{Z}}}'\mc{V}_n(\mc{H}_k). 
\eeqq
The symbol $\prod'$ denotes the full subcategory of the product category 
consisting of objects of the form $(V_k)_{k\in\bb{Z}}$, where 
$V_k\in N_0\mc{V}_n(\mc{H}_k)$ and $V_k=0$ for all $k<<0$. 
(Recall: $V^E_{0,k}=0$ for $k<<0$ for any AFT $E$.)

On the other hand, there is a functor in the opposite direction 
\beqq
\mc{E}_n:{\prod_{k\in\bb{Z}}}'\mc{V}_n(\mc{H}_k)\rightarrow
\mc{AFT}_n(\mc{H}).
\eeqq
An object of $\prod'_k\mc{V}_n(\mc{H}_k)$ is a sequence 
${\bf V}=(V_k)_{k\in\bb{Z}}$ of $\bb{Z}/2$-graded $Cl_{-n}$-modules. 
The degree $n$ AFT $\mc{E}_n({\bf V})$ is defined by 
$\mc{E}_n({\bf V})(S)=\mc{H}$ and 
\beqq        
\mc{E}_n({\bf V})_B(A_{x,y,\theta})=
  \begin{cases} 
    q^k & \t{on }V_k \\
    0 & \t{on }(\bigoplus_k V_k)^\perp
  \end{cases},
\eeqq
where $B$ is any object of $\mc{S}$ and $q=e^{2\pi i(x+iy)}$. 
For a morphism $({\bf f,A})=(f_k,A_k)_{k\in\bb{Z}}$ in 
$\prod'_k\mc{V}_n(\mc{H}_k)$, say from ${\bf V}=(V_k)$ to 
${\bf V}'=(V'_k)$, we have 
\beqq
\mc{E}_n(({\bf f,A}))=\left(\alpha_{\bf V,V'},
  \bigoplus_{k\in\bb{Z}}f_k,\bigoplus_{k\in\bb{Z}}A_k\right),
\eeqq
where $\alpha_{\bf V,V'}$ is the inclusion of 
$\{(0,k)|V_k\neq 0\}$ into 
$\{(0,k)|V'_k\neq 0\}$. 
Notice that 
(i) by definition, $f_k:V_k\rightarrow V'_k$ are embeddings, 
so $V'_k\neq 0$ if $V_k\neq 0$, and 
(ii) we have 
$\t{sp}(\tilde{G}_{\mc{E}_n({\bf V})},
L_{\mc{E}_n({\bf V})}-\tilde{G}_{\mc{E}_n({\bf V})}^2)=
\{(0,k)|V_k\neq 0\}$, etc.

The composition $\t{ind}_n^{S^1}\circ\mc{E}_n$ is the identity on 
$\prod'_k\mc{V}_n(\mc{H}_k)$. 
Also, for each object $E$ of $\mc{AFT}_n(\mc{H})$, there is a morphism  
\beqq
(\alpha_E,i_E,A_E):\mc{E}_n(\t{ind}_n^{S^1}(E))\rightarrow E,
\eeqq
where $\alpha_E$ is the inclusion 
$\{(0,k):V_k\neq 0\}\subset\t{sp}(\tilde{G}_E,L_E-\tilde{G}_E^2)$, 
$i_E$ is the inclusion 
$\bigoplus_k V^E_{0,k}\subset\bigoplus_{(\lambda,k)}V^E_{\lambda,k}$, 
and $A_E=\bigoplus_{\lambda>0} V^E_{\lambda,k}$.  
One easily checks that every diagram 
\beqq
\xymatrix{
  \mc{E}_n(\t{ind}_n^{S^1}(E))
    \ar[rrr]^{\phantom{aaaa}(\alpha_E,i_E,A_E)}
    \ar[d]_{\mc{E}_n(\t{ind}_n^{S^1}(\alpha,f,A))} &&& 
    E \ar[d]^{(\alpha,f,A)} \\
  \mc{E}_n(\t{ind}_n^{S^1}(E'))
    \ar[rrr]^{\phantom{aaaa}(\alpha_{E'},i_{E'},A_{E'})} &&& E'
}
\eeqq
commutes.
Therefore, $E\mapsto(\alpha_E,i_E,A_E)$ is a natural transformation 
from $\mc{E}_n\circ\t{ind}_n^{S^1}$. 
This proves the functors $\t{ind}_n^{S^1}$ and $\mc{E}_n$ induce a
homotopy equivalence
\beqq
|\mc{AFT}_n(\mc{H})|\simeq
\left|{\prod_k}'\mc{V}_n(\mc{H}_k)\right|=
{\prod_k}'|\mc{V}_n(\mc{H}_k)|,
\eeqq
where the last term is the space consisting of points 
$(x_k)_{k\in\bb{Z}}$ with $x_k\in|\mc{V}_n(\mc{H}_k)|$ and 
$x_k=0$ (the zero submodule) for $k<<0$.
By proposition \ref{V.Fred}, we have 
\beqq
\big[X,|\mc{AFT}_n(\mc{H})|\big]={\prod_k}'K^n(X)=K^n(X)[[q]][q^{-1}], 
\eeqq
where the second term consists of sequences $(a_k)_{k\in\bb{Z}}$ with 
$a_k\in K^n(X)$ and $a_k=0$ for $k<<0$, 
and the second equality identifies such a sequence with the formal 
Laurent series $\sum_k a_k q^k$. 
\qedsymbol

\newpage
\appendix

\titleformat{\section}[block]
  {\large\sffamily\filcenter\bfseries}
  {Appendix:}{.5em}{}

\renewcommand{\theequation}{\thesection.\arabic{equation}}

\section{Proofs of Two Lemmas} \label{lemmas}
In this appendix, we give the proofs of lemmas \ref{BHom.red} and 
\ref{BHom.Hom}. 
For any map $f:M\rightarrow N$ of super manifolds, we will denote 
by $f_B:M(B)\rightarrow N(B)$ the natural transformation 
$g\mapsto f\circ g$, where $B$ is any super manifold. 

{\renewcommand{\theequation}{\ref{BHom.red}}
\lemma
\addtocounter{equation}{-1}
\renewcommand{\theequation}{\thechapter.\arabic{equation}}
Suppose $(M_i,\chi_i)$, $i=1,2$, are Euclidean $(0|1)$- or $(1|1)$-manifolds, 
and $B$ is a super manifold with $B_\r=\{b\}$. 
Given a map
\beqq
\phi:\big(M_1(B),(\chi_1)_B\big)\rightarrow\big(M_2(B),(\chi_2)_B\big), 
\eeqq
there exists a unique map $\phi_\r:(M_1)_\r\rightarrow(M_2)_\r$ of 
ordinary manifolds so that the diagram 
\beqq
\xymatrix{
  \big(M_1(B),(\chi_1)_B\big)\ar[r]^\phi\ar[d] & 
    \big(M_2(B),(\chi_2)_B\big)\ar[d] \\
  (M_1)_\r\ar[r]^{\phi_\r} & (M_2)_\r \; , 
} 
\eeqq
with the vertical arrows given by $f\mapsto f_\r(b)$, commutes. 
In other words, $\phi(f)_\r=\phi_\r\circ f_\r$, $\forall f\in M_1(B)$. 
Furthermore, $(\phi'\phi)_\r=\phi'_\r\circ\phi_\r$ whenever 
$\phi$ and $\phi'$ are composable.
}

\proof
Assume that $(M_i)_\r$, $i=1,2$, are connected. 
This is sufficient because, for any super manifold $N$, 
the map $f\mapsto f_\r(b)$ sends different components of $N(B)$ to different 
components of $N_\r$.

The lemma is clear for $\dim M_1=\dim M_2=(0|1)$. 

Suppose $\dim M_1=(0|1)$ and $\dim M_2=(1|1)$. 
Since $(M_1)_\r$ is a point, the existence of $\phi_\r$ is equivalent to 
the map $f\mapsto\phi(f)_\r(b)$, $f\in M_1(B)$, being constant, or just 
locally constant, as $M_1(B)$ is connected. 
Fix $f_0\in M_1(B)$ and let $y=\phi(f_0)_\r(b)$. 
If $\lambda_1\in\mc{O}(M_1)^\od$ is an element such that 
$\chi_1=\lambda_1 d\lambda_1$, sending $\lambda$ to $\lambda_1$ 
defines an invertible map 
$\mu:(M_1,\chi_1)\stackrel{\sim}{\rightarrow}(\bb{R}^{0|1},\lambda d\lambda)$. 
On the other hand, there exists a neighborhood $U$ of $y\in(M_2)_\r$, 
together with an invertible map 
$\tau:((M_2)_U,\chi_2|_{(M_2)_U)})\stackrel{\sim}{\rightarrow}
(\bb{R}^{1|1}_{(a_1,a_2)},ds+\lambda d\lambda)$, for some $a_1<a_2$. 
(See the note below (\ref{Hom.11}).)
Since $\phi(f_0)\in(M_2)_U(B)\subset M_2(B)$, 
let $V\subset M_1(B)$ be a neighborhood of $f_0$ such that 
$\phi(V)\subset(M_2)_U(B)$. 
Also, let $V'=\mu_B(V)\subset\bb{R}^{0|1}(B)$. 
We have the following commutative diagram: 
\beqq
\xymatrix{
  \big(V,(\chi_1)_B\big) \ar[rr]^{\phi|_V\phantom{aaa}} \ar[d]_{\mu_B}^\cong 
    && \big((M_2)_U(B),(\chi_2|_{(M_2)_U})_B\big) \ar[d]_{\tau_B}^\cong \\
  \big(V',(\lambda d\lambda)_B\big) \ar[rr]^{\tilde{\phi}\phantom{aaa}} 
    && \big(\bb{R}^{1|1}_{(a,a')}(B),(ds+\lambda d\lambda)_B\big). 
}
\eeqq
By (\ref{BHom.01}), $\tilde{\phi}=\gamma_{z,\theta}$ or 
$\gamma_{z,\theta}\epsilon$, 
for some $(z,\theta)\in\mc{O}(B)^\ev\times\mc{O}(B)^\od$. 
For any $\eta\in V'$, 
\beqq
\tilde{\phi}(\eta)_\r(b)=\gamma_{z,\theta}(\pm\eta)_\r(b)=z^0\in(a_1,a_2),
\eeqq
which is independent of $\eta$. 
(Recall that $z^0$ denotes the unique real number such that 
$z-z^0$ is nilpotent.)
This means $\phi(f)_\r(b)=y$ for any $f\in V$, as desired. 

Now, suppose $\dim M_1=\dim M_2=(1|1)$. 
Fix $x\in(M_1)_\r$. 
Let 
\beqq
M_1(B)_x=\{f\in M_1(B)|f_\r(b)=x\}.
\eeqq
We need to show that the map $f\mapsto\phi(f)_\r(b)$ is constant on 
$M_1(B)_x$. 
Let $U\subset(M_1)_\r$ now be a neighborhood of $x$ and 
$s_U\in\mc{O}((M_1)_U)^\ev$, $\lambda_U\in\mc{O}((M_1)_U)^\od$ be local 
functions such that $\chi_1|_{(M_1)_U}=ds_U+\lambda_U d\lambda_U$. 
Define a map 
$\gamma:(\bb{R}^{0|1},\lambda d\lambda)\rightarrow(M_1,\chi_1)$ by 
$\gamma_\r(\bb{R}^0)=x$, $\gamma^*(\lambda_U)=\lambda$. 
Notice that any $f\in M_1(B)_x$ is in the image of 
$\gamma_B:\bb{R}^{0|1}(B)\rightarrow M_1(B)$;  
indeed, viewing $f^*(\lambda_U)$ as an element of $\bb{R}^{0|1}(B)$ using 
(\ref{R01B.R11B}), we have $\gamma_B(f^*(\lambda_U))=f$. 
As proved above, the composition
\beqq
\big(\bb{R}^{0|1}(B),(\lambda d\lambda)_B\big) 
  \stackrel{\gamma_B}{\rightarrow}
\big(M_1(B),(\chi_1)_B\big) \stackrel{\phi}{\rightarrow}
\big(M_2(B),(\chi_2)_B\big) \rightarrow (M_2)_\r\;,
\eeqq
where the last map is $g\mapsto g_\r(b)$, is constant. 
This proves $\phi$ is constant on $M_1(B)_x$. 

Given 
\beqq
\phi:\big(M_1(B),(\chi_1)_B\big)\rightarrow\big(M_2(B),(\chi_2)_B\big), \\
\phi':\big(M_2(B),(\chi_1)_B\big)\rightarrow\big(M_3(B),(\chi_3)_B\big), 
\eeqq
where $(M_i,\chi_i)$, $i=1,2,3$, are Euclidean $(0|1)$- or 
$(1|1)$-manifolds, we have 
\beqq
(\phi'\phi)_\r f_\r=(\phi'\phi)(f)_\r=\phi'(\phi(f))_\r
=\phi'_\r \phi(f)_\r=\phi'_\r \phi_\r f_\r, 
\eeqq
for any $f\in M_1(B)$. 
Any $x\in(M_1)_\r$ is the image of $f_\r$ for some $f\in M_1(B)$, 
hence we have 
$(\phi'\phi)_\r(x)=\phi'_\r\phi_\r(x)$, $\forall x\in(M_1)_\r$. 
\qedsymbol

{\it Note.}
We need not consider $(\dim M_1,\dim M_2)=((1|1),(0|1))$, 
because in this case,
\beqq
\t{Hom}\big(M_1(B),(\chi_1)_B;M_2(B),(\chi_2)_B\big)=\emptyset.
\eeqq 
It suffices and is easy to check this for 
$(M_1,\chi_1)=(\bb{R}^{1|1}_{(a_1,a_2)},ds+\lambda d\lambda)$, where 
$a_1<a_2$, and 
$(M_2,\chi_2)=(\bb{R}^{0|1},\lambda d\lambda)$. 

{\renewcommand{\theequation}{\ref{BHom.Hom}}
\lemma
\addtocounter{equation}{-1}
\renewcommand{\theequation}{\thechapter.\arabic{equation}}
Suppose $(M_i,\chi_i)$, $i=1,2$, are Euclidean $(0|1)$- or $(1|1)$-manifolds 
and $B$ is a super manifold with $B_\r=$ a point. 
There is a map 
\beqq
\ell:\t{Hom}\big(M_1(B),(\chi_1)_B;M_2(B),(\chi_2)_B\big)\rightarrow
\t{Hom}(M_1,\chi_1;M_2,\chi_2),
\eeqq
with the following properties:
\begin{itemize}
\item[(i)]
$\ell(\phi)_\r=\phi_\r$ for any 
$\phi\in\t{Hom}\big(M_1(B),(\chi_1)_B;M_2(B),(\chi_2)_B\big)$, 
\item[(ii)]
$\ell(\varphi_B)=\varphi$ for any 
$\varphi\in\t{Hom}(M_1,\chi_1;M_2,\chi_2)$, and 
\item[(iii)]
$\ell(\phi'\phi)=\ell(\phi')\ell(\phi)$ whenever $\phi$ and $\phi'$ 
are composable. 
\end{itemize}
}

\proof
First, consider the case in which each $(M_i,\chi_i)$, $i=1,2$, is either 
$(\bb{R}^{0|1},\lambda d\lambda)$ or 
$(\bb{R}^{1|1}_{(a_1,a_2)},ds+\lambda d\lambda)$, 
where $-\infty\leq a_1<a_2\leq +\infty$. 
Recall (\ref{Hom.00})-(\ref{Hom.11}) and (\ref{BHom.00})-(\ref{BHom.11}). 
We define
\beq \label{ell.basic}
\begin{array}{c}
\ell(1)=1,\quad \ell(\epsilon)=\epsilon, \\
\ell(\gamma_{z,\theta})=\gamma_{z^0},\quad
  \ell(\gamma_{z,\theta}\epsilon)=\gamma_{z^0}\epsilon, \\
\ell(\tau_{z,\theta})=\tau_{z^0},\quad
  \ell(\tau_{z,\theta}\epsilon)=\tau_{z^0}\epsilon,
\end{array}
\eeq
for any $(z,\theta)\in\mc{O}(B)^\ev\times\in\mc{O}(B)^\od$. 
(Again, $z^0$ is the unique real number such that $z-z^0$ is nilpotent.)
One easily verifies (\ref{ell.basic}) satisfies (i)-(iii). 

In general, consider a map
$\phi:\big(M_1(B),(\chi_1)_B\big)\rightarrow\big(M_2(B),(\chi_2)_B\big)$. 
Let $U_i\subset(M_i)_\r$, $i=1,2$, be open sets, such that 
$\phi_\r(U_1)\subset U_2$ and there are invertible maps 
$\varphi_i$ from $\big((M_i)_{U_i},\chi_i|_{(M_i)_{U_i}}\big)$ to 
either $(\bb{R}^{0|1},\lambda d\lambda)$ or 
$(\bb{R}^{1|1}_{(a_1,a_2)},ds+\lambda d\lambda)$ for some $a_1<a_2$. 
According to the definition of $\phi_\r$ in lemma \ref{BHom.red}, 
$\phi$ maps $(M_1)_{U_1}(B)$ into $(M_2)_{U_2}(B)$. 
Using (\ref{ell.basic}), we define 
\beq \label{ell.local}
\ell(\phi|_{(M_1)_{U_1}(B)})=
  \varphi_2^{-1}\circ
  \ell\big((\varphi_2)_B\circ\phi|_{(M_1)_{U_1}(B)}\circ
    (\varphi_1)_B^{-1}\big)
  \circ\varphi_1.
\eeq
It follows from properties (ii), (iii) for (\ref{ell.basic}) that the right 
hand side is independent of the choices of $\varphi_i$. 
This implies the maps $\ell(\phi|_{(M_1)_{U_1}(B)})$ for various 
$U_1$ agree on overlaps, and thus glue into a single map 
$\ell(\phi):(M_1,\chi_1)\rightarrow(M_2,\chi_2)$. 

It suffices to verify (i)-(iii) locally, 
i.e. for (\ref{ell.local}). 
For example, 
$\ell(\phi|_{(M_1)_{U_1}(B)})_\r=(\phi|_{(M_1)_{U_1}(B)})_\r$ 
follows from lemma \ref{BHom.red}, 
the fact that (\ref{ell.basic}) satifies (i), 
and $((\varphi_i)_B)_\r=(\varphi_i)_\r$, $i=1,2$. 
The proofs of (ii) and (iii) are similar. 
\qedsymbol


\newpage

{\footnotesize
\begin{thebibliography}{MMMM99}

\bibitem[ABS64]{ABS}
M.~F.~Atiyah, R.~Bott and A.~Shapiro,
\emph{Clifford modules},
Topology, 3:3--38, 1964.

\bibitem[AHS01]{AHS} 
M.~Ando, M.~Hopkins and N.~Strickland,
\emph{Elliptic spectra, the Witten genus and the theorem of the cube},
Invent. Math. 146 (2001), no.3, 595--687.

\bibitem[AS69]{AS.Fred}
M.~F.~Atiyah and I.~Singer,
\emph{Index theory for skew-adjoint Fredholm operators},
Inst. Hautes \'Etudes Sci. Publ. Math., 37:5--26, 1969.

\bibitem[BT89]{BottTaubes}
R.~Bott and C.~Taubes,
\emph{On the rigidity theorems of Witten},
J.~Amer. Math. Soc., vol.~2, 1:137--186, 1989.

\bibitem[DM99]{QFS.susy}
P.~Deligne and J.~Morgan,
\emph{Notes on supersymmetry}
in: Quantum fields and strings: a course for mathematicians, vol.~1,
41--97, Amer. Math. Soc., Providence, 1999.

\bibitem[DT58]{DoldThom}
A.~Dold and R.~Thom,
\emph{Quasifaserungen und unendliche symmetrische produkts},
Annals of Math., 61, 239--281, 1958.

\bibitem[Hop02]{Hopkins.ICM}
M.~Hopkins,
\emph{Algebraic topology and modular forms},
in: Proceedings of the International Congress of Mathematicians, 
Vol. I, 291--317, Higher Ed. Press, 2002.

\bibitem[Kui65]{Kuiper}
N.~H.~Kuiper,
\emph{The homotopy type of the unitary group of Hilbert space},
Topology, 3:19--30, 1965.

\bibitem[LRS93]{LRS}
P.~Landweber, D.~Ravenel and R.~Stong,
\emph{Periodic cohomology theories defined by elliptic curves},
in: The \v Cech centennial, 317--337,
Contemp. Math., 181, Amer. Math. Soc., 1995.

\bibitem[LM89]{LM}
H.~B.~Lawson and M.-L.~Michelsohn,
\emph{Spin geometry},
Princeton Univ. Press, 1989.

\bibitem[Och87]{Ochanine}
S.~Ochanine,
\emph{Sur les genres multiplicatifs d\'efinis par les int\'egrales 
elliptiques},
Topology, 26:143--151, 1987.

\bibitem[Qui73]{Quillen.algK}
D.~Quillen,
\emph{Higher algebraic K-theory I}
in: Algebraic K-theory I: higher K-theories, 85--147,
Springer, Berlin, 1973.

\bibitem[Rui91]{Rudin}
W.~Rudin,
\emph{Functional analysis},
McGraw-Hill, New York, 1991.

\bibitem[Seg74]{Segal.Gamma}
G.~Segal,
\emph{Categories and cohomology theories},
Topology, 13:293--312, 1974.

\bibitem[Seg77]{Segal.Khlgy}
G.~Segal,
\emph{K-homology theory and algebraic K-theory}
in: K-theory and operator algebras, 113--127,
Springer, Berlin, 1977.

\bibitem[Seg88]{Segal.Ell}
G.~Segal,
\emph{Elliptic cohomology},
S\'eminaire Bourbaki, Vol.~1987/88.
Ast\'erisque No.~161-162 (1988), Exp. No. 695, 4, 187--201.

\bibitem[Seg04]{Segal.CFT}
G.~Segal,
\emph{The definition of conformal field theory}
in: Topology, geometry and quantum field theory, 423--577,
Cambridge Univ. Press, 2004.

\bibitem[Sil94]{Silverman}
J.~Silverman,
\emph{Advanced topics in the arithmetic of elliptic curves},
Springer-Verlag, New York, 1994.

\bibitem[ST04]{ST}
S.~Stolz and P.~Teichner,
\emph{What is an elliptic object?},
in: Topology, geometry and quantum field theory, 247--343,
Cambridge Univ. Press, 2004.

\bibitem[Wit87]{Witten.ell.genus}
E.~Witten,
\emph{Elliptic genera and quantum field theory},
Comm. Math. Phys. 109 (1987), no.~4, 525--536.

\bibitem[Wit99]{Witten.indD}
E.~Witten,
\emph{Index of Dirac operators}
in: Quantum fields and strings: a course for mathematicians, vol.~1,
475--511, Amer. Math. Soc., Providence, 1999.

\end{thebibliography}
}
\end{document}